\documentclass[twoside]{amsart}
\usepackage[centertags]{amsmath}
\usepackage{amsfonts}
\usepackage{amssymb}
\usepackage{dsfont}
\usepackage{amsthm}
\usepackage[ansinew]{inputenc}
\usepackage[cmtip,poly,all,arc]{xy}
\usepackage{graphicx}
\newtheorem{thm}[equation]{Theorem}
\newtheorem{cor}[equation]{Corollary}
\newtheorem{lem}[equation]{Lemma}
\newtheorem{prop}[equation]{Proposition}
\theoremstyle{definition}
\newtheorem{defn}[equation]{Definition}
\theoremstyle{remark}
\newtheorem{rem}[equation]{Remark}

\numberwithin{equation}{section}

\newcommand{\abs}[1]{\left\vert#1\right\vert}
\newcommand{\set}[1]{\left\{#1\right\}}
\newcommand{\Real}{\mathbb R}

\newcommand{\To}{\longrightarrow}

\newcommand{\I}{\mathbb{I}}

\newcommand{\C}[1]{\mathbf{#1}} 

\def\nill{{nil}}
\def\abb{{ab}}

\def\r{\rightarrow} 
\def\l{\leftarrow} 
\def\rr{\Rightarrow} 

\def\hom{\operatorname{Hom}}

\newcommand{\sym}[1]{\operatorname{Sym}(#1)}

\newcommand{\symt}[1]{\operatorname{Sym}_\vc(#1)}
\newcommand{\symtt}{\operatorname{Sym}_\vc}

\def\sign{\operatorname{sign}}

\def\aut{\operatorname{Aut}}

\def\cc{\ul{\circledcirc}}
\def\ca{\circledcirc}
\def\co{\bar{\otimes}}
\def\wc{\ul{\wedge}}
\def\wa{\wedge}
\def\sc{\ul{\#}}
\def\sa{\#}

\def\st{\stackrel} 

\def\ul{\underline} 
\def\ol{\overline} 

\newcommand{\hopf}{\mathit{Hopf}}

\def\coker{\operatorname{Coker}}
\renewcommand{\ker}{\operatorname{Ker}}
\renewcommand{\aleph}{\mathbb{H}}

\def\Z{\mathbb{Z}}

\def\N{\mathds{N}}

\def\S{\Sigma}

\def\L{\Omega}

\newcommand{\ad}{\mathsf{Ad}}

\newcommand{\grupo}[1]{\langle #1\rangle}

\newcommand{\vc}{\Box}    
\newcommand{\vi}{\boxminus}  

\newcommand{\Hg}{\mathrm{G}}

\begin{document}

\title[Smash products for secondary homotopy groups]{Smash products for secondary homotopy groups}%
\author{Hans-Joachim Baues and Fernando Muro}%
\address{Max-Planck-Institut f\"ur Mathematik, Vivatsgasse 7, 53111 Bonn, Germany}%
\email{baues@mpim-bonn.mpg.de, muro@mpim-bonn.mpg.de}%

\thanks{The second author was partially supported
by the project MTM2004-01865 and the MEC postdoctoral fellowship EX2004-0616.}%

\begin{abstract}
We construct a smash product operation on secondary homotopy
groups yielding the structure of a lax symmetric monoidal functor. Applications on cup-one products, Toda
brackets and Whitehead products are considered.
\end{abstract}
\maketitle
\begin{footnotesize}
\tableofcontents
\end{footnotesize}

\section*{Introduction}

The classical homotopy groups $\pi_nX$, $n\geq 0$, of a pointed space $X$
give rise to a graded abelian group $\Pi_*X$
obtained by additivization in low dimensions. In particular
$\Pi_nX=\pi_n X$ for $n\geq 2$,
$\Pi_1X=(\pi_1X)_\abb$ is the abelianized fundamental group, and
$\Pi_0X=\Z[\pi_0X]$ is the free abelian group on the pointed set of
path components of $X$. The smash product on homotopy groups induces a natural
homomorphism of graded abelian groups
\begin{equation*}\tag{1}
\Pi_*X\otimes\Pi_*Y\st{\wedge}\To\Pi_*(X\wedge Y),
\end{equation*}
which carries $f\otimes g$ with $f\colon S^n\r X$ and $g\colon S^m\r Y$ to $f\wedge g\colon S^{n+m}\r X\wedge Y$. 
This shows that $\Pi_*$ is a lax symmetric monoidal functor from pointed spaces to graded abelian groups.

The smash product (1) can be used for example to define the Whitehead product on homotopy groups, compare
Section \ref{swp}.

The purpose of this paper is to generalize these properties of primary homotopy groups on the level of secondary
homotopy theory.

Secondary homotopy operations like Toda brackets \cite{toda} or cup-one products \cite{K1p}, \cite{c1Tb}, are defined by pasting
tracks, where tracks are homotopy classes of homotopies. Since secondary homotopy operations play a crucial role in homotopy
theory it is of importance to develop the algebraic theory of tracks. We do this by introducing secondary
homotopy groups of a pointed space $X$
$$\Pi_{n,*}X=\left(\Pi_{n,1}X\st{\partial}\r\Pi_{n,0}X\right)$$
which have the structure of a quadratic pair module, see Sections \ref{sqpm} and \ref{el1}. 
Here $\partial$ is a group homomorphism
with cokernel $\Pi_nX$ for $n\geq 0$ and kernel $\Pi_{n+1}X$ for $n\geq 3$. 

We define $\Pi_{n,*}X$ for $n\geq 2$ directly in terms of maps 
$S^n\r X$ and tracks from such maps to the trivial map. For $n\geq 0$ the functor $\Pi_{n,*}$ is an additive
version of the functor $\pi_{n,*}$ studied in \cite{2hg1}. 

We introduce and study
the smash product morphism for additive secondary homotopy
groups
\begin{equation*}\tag{2}
\Pi_{*,*}X\odot\Pi_{*,*}Y\st{\wedge}\To\Pi_{*,*}(X\wedge Y).
\end{equation*}
Here one needs the symmetric monoidal structure $\odot$ of the category
of quadratic pair modules $\C{qpm}$ which is based on
the symmetric monoidal structure on the category of square groups constructed in \cite{qaI}. 
The smash product morphism (2) is compatible with the associativity isomorphisms
but it is not directly compatible with the commutativity isomorphisms.

In order to deal with commutativity we need the action
of the symmetric track group $\symt{n}$ on $\Pi_{n,*}X$ in \cite{2hg2}. We show that $\wedge$ in (2) is
equivariant with respect to this action, and  is commutative up to
the action of a shuffle permutation. This leads to the definition of the
symmetric monoidal category $\C{qpm}_0^{\symtt}$ with objects given by symmetric 
sequences of quadratic pair modules with extra structure. Then the morphism
(2) induces a morphism in $\C{qpm}_0^{\symtt}$ for which the
associativity and commutativity isomorphisms are compatible with the
symmetric monoidal structure of $\C{qpm}_0^{\symtt}$. Therefore $\Pi_{*,*}$ considered as a functor to  the category 
$\C{qpm}_0^{\symtt}$ is, in fact, a lax symmetric monoidal functor. 

The smash product (2) is used to define the
Whitehead product on secondary homotopy groups, compare Section \ref{swp}.

As an illustrating application of the results in this paper we prove a formula of Barratt-Jones-Mahowald
on unstable cup-one products, see Section \ref{c1}. This formula was stated in \cite{K1p}, but a proof did not
appear in the literature. A further application yields a formula for a triple Toda bracket which generalizes a
well-known formula in \cite{toda}, see Section \ref{t11}.

In a sequel of this paper we generalize the theory of secondary homotopy groups to symmetric spectra. There we
show that the smash
product operation defined in this paper endows the secondary homotopy groups $\Pi_{*,*}R$ of a fibrant ring spectrum $R$
with a graded algebra structure in the category of quadratic pair modules. The graded algebra $\Pi_{*,*}R$
determines all Toda brackets in $\pi_*R$, which can be regarded as Massey products in $\Pi_{*,*}R$. Moreover,
$\Pi_{*,*}R$ determines the universal matrix Toda bracket in the category of finitely generated free $R$-modules.
If $R$ is an $E_\infty$-ring spectrum then $\Pi_{*,*}R$ is a commutative algebra up to coherent homotopies in
$\C{qpm}_0^{\symtt}$ which encodes not only Toda brackets, but also cup-one products in a purely algebraic way.

The paper consists of three parts. The first part is concerned with the algebra needed for the statements of the
main theorems. In Part \ref{II} we present our main results and we give applications. Part \ref{III} contains the
construction of the smash product operation for additive secondary homotopy groups. There we prove all the
properties which imply our main results.

\part{Quadratic pair modules and their tensor product}\label{I}

In this part we describe the algebraic concepts needed for the structure of secondary homotopy groups. We
introduce the category of quadratic pair modules and we show that this category is symmetric monoidal. The tensor
product of quadratic pair modules is related to the exterior cup-products in the category $\C{Top}^*$.

\section{Square groups and quadratic pair modules}\label{sqpm}

We first recall the notion of square group, see \cite{ecg} and \cite{qaI}.

\begin{defn}\label{primen}
A \emph{square group} $X$ is a diagram
$$X=(X_e\mathop{\leftrightarrows}\limits^P_HX_{ee})$$
where $X_e$ is a group with an additively written group law, $X_{ee}$
is an abelian group, $P$ is a homomorphism, $H$ is a quadratic map, i. e. a function such
that the cross effect
$$(a|b)_H=H(a+b)-H(b)-H(a)$$
is linear in $a,b\in X_e$, and the following relations are
satisfied, $x,y\in X_{ee}$,
\begin{enumerate}
\item $(Px|b)_H=0$, $(a|Py)=0$,
\item $P(a|b)_H=[a,b]$,
\item $PHP(x)=P(x)+P(x)$.
\end{enumerate}
Here $[a,b]=-a-b+a+b$ is the commutator bracket.
The cross effect induces a homomorphism
\begin{equation}\label{ce}
(-|-)_H\colon\otimes^2\coker P\To X_{ee}.
\end{equation}
Here $\otimes^2A=A\otimes A$ is the tensor square of an
abelian group $A$. If (\ref{ce}) is an isomorphism we will say that the square group $X$ is \emph{good}.
The function $$T=HP-1\colon X_{ee}\To X_{ee}$$ is an involution, i.
e. a homomorphism with $T^2=1$. Moreover,
$$\Delta\colon X_e\To X_{ee}\colon x\mapsto (x|x)_H-H(x)+TH(x)$$
is a homomorphism which satisfies $T\Delta=-\Delta$.

A morphism of square groups $f\colon X\r Y$ is given by
homomorphisms $f_e\colon X_e\r Y_e$,
$f_{ee}\colon X_{ee}\r Y_{ee},$
commuting with $P$ and $H$. 

As an example of square group we can consider
$$\Z_\nill=(\Z\mathop{\leftrightarrows}\limits^P_H\Z)$$
with $P=0$ and $H(n)=\binom{n}{2}=\frac{n(n-1)}{2}$. This is the unit object of the symmetric monoidal structure
defined in the next section.

\end{defn}

We refer the reader to \cite{qaI} where the quadratic algebra of
square groups is developed. We need square groups for the definition
of quadratic pair modules as follows.

\begin{defn}\label{qpm}
A \emph{quadratic pair module} $C$ is a morphism $\partial\colon
C_{(1)}\r C_{(0)}$ between square groups
\begin{eqnarray*}
C_{(0)}&=&(C_{0}\mathop{\leftrightarrows}^{P_0}_HC_{ee}),\\
C_{(1)}&=&(C_{1}\mathop{\leftrightarrows}^P_{H_1}C_{ee}),
\end{eqnarray*}
such that $\partial_{ee}=1\colon C_{ee}\r C_{ee}$ is the identity
homomorphism. In particular $C$ is completely determined by
the diagram
\begin{equation}\label{ya}
\xymatrix{&C_{ee}\ar[ld]_P&\\C_{1}\ar[rr]_\partial&&C_{0}\ar[lu]_H}
\end{equation}
where $\partial=\partial_e$, $H_1=H\partial$ and $P_0=\partial P$.

Morphisms of quadratic pair modules $f\colon C\r D$ are therefore
given by group homomorphisms $f_{0}\colon C_{0}\r D_{0}$,
$f_{1}\colon C_{1}\r D_{1}$, $f_{ee}\colon C_{ee}\r D_{ee}$,
commuting with $H$, $P$ and $\partial$ in (\ref{ya}). 
They form a
category denoted by $\C{qpm}$.

A quadratic pair module $C$ is \emph{$0$-good} if the square group $C_{(0)}$ is good. The full subcategory of
$0$-good quadratic pair modules will be denoted by $\C{qpm}_0\subset\C{qpm}$.
\end{defn}

\begin{rem}\label{nse}
The category $\C{squad}$ of stable quadratic modules is described in \cite{ch4c} IV.C and \cite{2hg1}. Quadratic
modules in general are discussed in \cite{ch4c} and \cite{csch}, they are special $2$-crossed modules in the
sense of \cite{2cm}.
There is a faithful forgetful functor 
$\C{qpm}\r \C{squad}$
sending $C$ to the stable quadratic
module
\begin{equation*}
\otimes^2(C_{0})_\abb\st{P(-|-)_H}\To C_{1}\st{\partial}\To C_{0}.
\end{equation*}
Here $G_\abb$ denotes the abelianization of a group $G$.
Such a stable quadratic module gives rise to a crossed module $\partial\colon C_1\r C_0$ where $C_0$ acts on the
right of $C_1$ by the formula, $x\in C_1$, $y\in C_0$,
$$x^y=x+P(\partial(x)|y)_H,$$ so that we also get a forgetful functor $\C{qpm}\r \C{cross}$ to the track
category of crossed modules. Tracks in $\C{qpm}$ (i. e. invertible $2$-morphisms) will be
considered in Section \ref{tcqpm} below by using this forgetful functor.
\end{rem}





\section{The tensor product of square groups}

We now recall the notion of tensor product of square groups which is essential for this paper. This tensor
product, first defined in \cite{qaI}, originates from properties of the exterior cup-products in the
next section and in Part \ref{III}.

\begin{defn}\label{tpsg}
The \emph{tensor product} $X\odot Y$ of square groups $X, Y$ is
defined as follows. The group $(X\odot Y)_e$ is generated by the
symbols $x\cc y$, $x\ca y$, $a\co b$ for $x\in X_e$, $y\in Y_e$,
$a\in X_{ee}$ and $b\in Y_{ee}$, subject to the following relations
\begin{enumerate}
\item the symbol $a\co b$ is bilinear and central,
\item the symbol $x\cc y$ is right linear, $x\cc(y_1+y_2)=x\cc y_1+x\cc y_2$,
\item the symbol $x\ca y$ is left linear, $(x_1+x_2)\ca y=x_1\ca y+x_2\ca y$,
\item $P(a)\cc y=a\co \Delta(y)$,
\item $x\ca P(b)=\Delta(x)\co b$,
\item $T(a)\co T(b)=-a\co b$,
\item $x\cc y-x\ca y=H(x)\co TH(y)$.
\end{enumerate}
The abelian group $(X\odot Y)_{ee}$ is defined as the tensor product
$X_{ee}\otimes Y_{ee}$. The homomorphism $$P\colon(X\odot Y)_{ee}\To
(X\odot Y)_{e}$$ is $P(a\otimes b)=a\co b$, and
$$H\colon(X\odot Y)_{e}\To (X\odot Y)_{ee}$$
is the unique quadratic map satisfying
\begin{eqnarray*}
H(x\cc y)&=&(x|x)_H\otimes H(y)+H(x)\otimes\Delta(y),\\
H(x\ca y)&=&\Delta(x)\otimes H(y)+ H(x)\otimes (y|y)_H,\\
H(a\co b)&=& a\otimes b-T(a)\otimes T(b),\\
(a\co b|-)_H&=&(-|a\co b)_H=0,\\
(a\cc b|c\cc d)_H&=&(a\cc b|c\ca d)_H\\&=&(a\ca b|c\cc
d)_H\\&=&(a\ca b|c\ca d)_H\\&=&(a|c)_H\otimes (b|d)_H.
\end{eqnarray*}

Relation (7) above shows that $(X\odot Y)_e$ is actually generated just
by $x\cc y$ and $a\co b$. A complete list of relations for these generators is
given by (1), (2), (4) and (6) above together with
\begin{enumerate}\setcounter{enumi}{7}
\item $(x_1+x_2)\cc y=x_1\cc y+ x_2\cc y+(x_2|x_1)_H\co H(y)$,
\item $x\cc P(b)=(x|x)_H\co b$.
\end{enumerate}

Similarly $(X\odot Y)_e$ is also generated by just $x\ca y$ and $a\co b$ with
relations (1), (3), (5) and (6) above together with
\begin{enumerate}\setcounter{enumi}{9}
\item $x\ca (y_1+y_2)=x\ca y_1+ x\ca y_2+H(x)\co(y_2|y_1)_H$,
\item $P(a)\ca y=a\co (y|y)_H$.
\end{enumerate}
\end{defn}

As proved in \cite{qaI} the tensor product of square groups is a symmetric monoidal
structure on the category of square groups with unit $\Z_\nill$ in Definition \ref{primen}. The associativity isomorphism
$$(X\odot Y)\odot Z\cong X\odot (Y\odot Z)$$
is given by $(x\cc y)\cc z\mapsto x\cc(y\cc z)$, $(a\co b)\cc z\mapsto a\co(b\otimes\Delta(z))$ and
$(a\otimes b)\co c\mapsto a\co (b\otimes c)$ at the $e$-level and by the associativity isomorphism 
$$(X_{ee}\otimes Y_{ee})\otimes Z_{ee}\cong X_{ee}\otimes (Y_{ee}\otimes Z_{ee})$$
for the tensor
product of abelian groups at the $ee$-level. 
The symmetry isomorphism
$$\tau_\odot\colon X\odot Y\cong Y\odot X$$ is defined on $e$-groups
by $x\cc y\mapsto y\ca x$, $x\ca y\mapsto y\cc x$, and $a\co b\mapsto b\co a$, and on
${ee}$-groups by the standard symmetry isomorphism 
$$\tau_\otimes\colon X_{ee}\otimes Y_{ee}\cong Y_{ee}\otimes X_{ee}$$
for the tensor product of abelian groups. The unit isomorphism 
$$\Z_\nill\odot X\cong X$$
is defined on $e$-groups by the formulas $n\ca x\mapsto n\cdot x$ and $n\co a\mapsto n\cdot P(a)$.

\section{Exterior cup-products}\label{ecp}

We will work with the track category $\C{Top}^*$ of compactly generated pointed spaces. A track category is a
category enriched in groupoids, i. e. a $2$-category where all $2$-morphisms are vertically invertible. A
$2$-morphism in a track category is also termed a \emph{track}, and a $2$-functor between track categories 
is called a \emph{track functor}. 
Tracks in $\C{Top}^*$ are homotopy classes of homotopies between pointed maps. The
identity track on a pointed map $f\colon X\r Y$, also called the trivial track, will be denoted by
$0^\vc_f$. We
use the symbol $\vc$ for the vertical composition, and $F^\vi$
denotes the vertical inverse of a track $F$. Horizontal composition is denoted by juxtaposition. 

The smash product of pointed spaces $X, Y$ in $\C{Top}^*$ is the
quotient space 
\begin{equation}\label{smp}
X\wedge Y=X\times Y/X\vee Y
\end{equation} 
where the coproduct
$X\vee Y$ in $\C{Top}^*$ admits the canonical inclusion $X\vee
Y\subset X\times Y$ to the product.
If $\sigma$ is a permutation of $\set{1,\dots,n}$ the map
\begin{equation}\label{permuconv}
\sigma\colon X_1\wedge\dots\wedge X_n\To X_{\sigma^{-1}(1)}\wedge\dots\wedge X_{\sigma^{-1}(n)},
\end{equation}
induced by the permutation of coordinates according to $\sigma$, is
also denoted by $\sigma$. For the sake of simplicity we will ocasionally omit the permutation $\sigma$ in the diagrams
and equations where it is understood. 
Given a subset
$\set{i_1,\dots,i_k}\subset\set{1,\dots,n}$ we denote by
$\sigma=(i_1\,\dots\, i_k)$ the permutation defined by $\sigma
(i_s)=i_{s+1}$ if $1\leq s <k$, $\sigma(i_k)=i_1$ and $\sigma (m)=m$
otherwise.

The smash product is a track functor
$$\wedge \colon\C{Top}^*\times\C{Top}^*\To\C{Top}^*.$$
It is defined as usually at the level of pointed spaces and pointed
maps. The pointed space $I_+$ is the disjoint union of the interval $I=[0,1]$ with an outer base point $*$. Let $F\colon f\rr g$ be a track between maps $f,g\colon A\r B$
represented by a homotopy $F\colon I_+\wedge A\r B$,  and
let $h\colon X\r Y$ be another map. Then the track $F\wedge h\colon
f\wedge h\rr g\wedge h$ is represented by a homotopy $F\wedge h\colon
I_+\wedge A \wedge X\r B\wedge Y$ and $h\wedge F\colon h\wedge f\rr
h\wedge g$ is represented by the composite
$$I_+\wedge X\wedge A\cong X\wedge I_+\wedge A\st{h\wedge F}\To Y\wedge B.$$ 
If $G\colon h\rr k$ is another track then
$$F\wedge G=(g\wedge G)\vc(F\wedge h)=(F\wedge k)\vc(f\wedge G)\colon f\wedge h\rr g\wedge k.$$

The smash product defines a symmetric monoidal structure in the
category $\C{Top}^*$. The unit object is the
$0$-sphere $S^0$.

\begin{defn}
Given maps $f\colon\S A\r \S B$ and $g\colon\S X\r\S Y$ the
\emph{left exterior cup-product} $f\sa  g$ is the composite
\begin{equation*}
S^1\wedge A\wedge X\st{f\wedge X}\To S^1\wedge B\wedge X
\cong B\wedge S^1\wedge X\st{B\wedge g}\To B\wedge
S^1\wedge Y\cong S^1\wedge B\wedge Y.
\end{equation*}
Similarly the \emph{right exterior cup-product} $f\sc g$ is the composite
\begin{equation*}
S^1\wedge A\wedge X\cong A\wedge S^1\wedge X\st{A\wedge
g}\To A\wedge S^1\wedge Y\cong S^1\wedge A\wedge Y
\st{f\wedge Y}\To S^1\wedge B\wedge Y.
\end{equation*}
The equality
\begin{eqnarray}\label{tee}
f\sa g=(2\; 3)(g\sc f)(2\;3)
\end{eqnarray}
is always satisfied. 
\end{defn}

These constructions give rise to homotopy operations called
\emph{exterior cup-products}
$$\sa ,\sc \colon [\S A,\S B]\times[\S X,\S Y]\To[\S A\wedge X,\S B\wedge Y].$$
see \cite{ccghc} II.1.14. The operation $\sa $ is left-linear and $\sc $ is
right-linear,
\begin{eqnarray*}
(f_1+f_2)\sa g&=&f_1\sa  g+f_2\sa g,\\
f\sc (g_1+g_2)&=&f\sc g_1+f\sc  g_2.
\end{eqnarray*}

Given a pointed discrete set $E$ we denote $\vee_ES^n=\S^n E$. Then $$\pi_1(\vee_ES^1)=\grupo{E}$$ is the free group with
basis $E-*$, and $$\pi_n(\vee_ES^n)=\Z[E]$$ is the free abelian group with basis $E-*$ for $n\geq 2$. We write
$$\grupo{E}_\nill$$ for the free group of nilpotency class $2$ (\emph{nil-group} for short) with basis $E-*$, which is obtained from
$\grupo{E}$ by dividing out triple commutators.

If $A=X=S^0$, and $B=E$, $Y=\bar{E}$ are pointed sets then the
exterior cup-products are functions
$$\sa ,\sc \colon\grupo{E}\times\grupo{\bar{E}}\To\grupo{E\wedge\bar{E}}.$$
These functions factor in a unique way through the natural
projection onto the nilization in the following way
\begin{equation}\label{edi}
\xymatrix{\grupo{E}\times\grupo{\bar{E}}\ar[r]^-{\sa ,\sc
}\ar@{->>}[d]&\grupo{E\wedge
\bar{E}}\ar@{->>}[d]\\
\grupo{E}_{nil}\times\grupo{\bar{E}}_{nil}\ar[r]^-{\sa ,\sc
}&\grupo{E\wedge \bar{E}}_{nil}}
\end{equation}

A free nil-group $\grupo{E}_\nill$ on a pointed set $E$ gives rise to
a square group 
\begin{equation}\label{znil}
\Z_\nill[E]=(\grupo{E}_\nill \mathop{\leftrightarrows}\limits^P_H \otimes^2\Z[E])
\end{equation}
defined by $P(a\otimes b)=[b,a]$, $H(e)=0$ for any $e\in E$ and
$(s|t)_H=t\otimes s$ so that $\Z_\nill[S^0]=\Z_\nill$ in Definition \ref{primen}. These square groups are the
main examples of good square groups in the sense of Definition \ref{primen}. For (\ref{znil}) the involution $T$ is up to sign the
interchange of factors in the tensor square $T(a\otimes b)=-b\otimes a$ and $\Delta$ is defined
by $\Delta(e)=e\otimes e$ for $e\in E$. Recall that we denote $$\tau_\otimes\colon A\otimes B\cong B\otimes A$$ to the symmetry
isomorphism for the tensor product of abelian groups, which should not be confused with $T=-\tau_\otimes$ in this case.

The next proposition is essentially \cite{qaI} 34. It shows the connection between the tensor product of square groups and
the exterior cup-products.

\begin{prop}\label{facil}
Given two pointed sets $E$ and $\bar{E}$ there is a square group isomorphism
$$\Z_\nill[{E}]\odot\Z_\nill[\bar{E}]\st{\cong}\To\Z_\nill[E\wedge\bar{E}]$$
defined on the $e$-groups
by $x\cc y\mapsto x\sc y$, $x\ca y\mapsto x\sa  y$, 
and on the ${ee}$-groups by
\begin{eqnarray*}
1\otimes\tau_\otimes\otimes1\colon\Z[E]\otimes\Z[E]\otimes\Z[\bar{E}]\otimes\Z[\bar{E}]&\cong&
\Z[E]\otimes\Z[\bar{E}]\otimes\Z[E]\otimes\Z[\bar{E}],\\
a\otimes b\otimes c\otimes d&\mapsto& a\otimes c\otimes b\otimes d.
\end{eqnarray*}
\end{prop}

This is the quadratic analogue of the well-known fact
that free abelian groups have the tensor product
$$\Z[E]\otimes\Z[\bar{E}]\cong\Z[E\wedge\bar{E}].$$

\section{The symmetric monoidal category of quadratic pair modules}\label{ape}

A \emph{pair} in a category $\C{C}$ is a morphism $f\colon X\r Y$ in
$\C{C}$. Let $\C{Pair}(\C{C})$ be the category of such pairs.
Morphisms $(\alpha,\beta)\colon f\r f'$ in $\C{Pair}(\C{C})$ are given by morphisms $\alpha\colon X\r
X'$ and $\beta\colon Y\r Y'$ in $\C{C}$ satisfying $\beta f=f'\alpha$. A
quadratic pair module is a special pair
in the category $\C{SG}$ of square groups and the category
$\C{qpm}$ of quadratic pair modules is a full reflective subcategory of
$\C{Pair}(\C{SG})$. The left adjoint to the inclusion $\C{qpm}\subset\C{Pair}(\C{SG})$, i. e. the reflection
functor,
is denoted by
$$\Phi\colon\C{Pair}(\C{SG})\To\C{qpm}.$$
Given a pair $f\colon D\r C$ in $\C{SG}$ we have $\Phi(f)_{(0)}=C$, so that $\Phi(f)_0=C_e$ and
$\Phi(f)_{ee}=C_{ee}$. Moreover, $\Phi(f)_1$ is the quotient group
$$\Phi(f)_1=D_e\times C_{ee}/\sim,$$
\begin{eqnarray*}
(0,f_{ee}(d))&\sim&(P(d),0),\;\;d\in D_{ee},\\
(0,HP(c))&\sim&(0,2c),\;\;c\in C_{ee}.
\end{eqnarray*}
The operators $P$ and $H$ for $\Phi(f)_{(1)}$ and the homomorphism $\partial\colon\Phi(f)_1\r\Phi(f)_0=C_e$ are
defined by the formulas, $c\in C_{ee}$, $d\in D_{e}$,
\begin{eqnarray*}
P(c)&=&(0,c),\\
H(d,c)&=&f_{ee}H(d)+HP(c),\\
\partial(d,c)&=&f_e(d)+P(c).
\end{eqnarray*}
The unit of this adjunction is a natural morphism in
$\C{Pair}(\C{SG})$ 
\begin{equation}\label{uni}
\xymatrix{D\ar[r]^{f}\ar[d]_{\upsilon}&C\ar@{=}[d]\\
\Phi(f)_{(1)}\ar[r]^\partial&\Phi(f)_{(0)}}
\end{equation}
and is given by $\upsilon_e(d)=(d,0)$ for $d\in D_e$ and $\upsilon_{ee}(d')=f_{ee}(d')$ for $d'\in D_{ee}$.
We use the functor $\Phi$ for the following definition of the tensor
product in $\C{qpm}$.

The category $\C{qpm}$ is a symmetric monoidal category. This
structure is inherited from the tensor product in $\C{SG}$ described
above. More precisely, the \emph{tensor product} $C\odot D$ of two
quadratic pair modules $\partial\colon C_{(1)}\r
C_{(0)},\partial\colon D_{(1)}\r D_{(0)}$ is given 
as follows. Consider the push-out diagram 
\begin{equation}\label{tp}
\xymatrix{C_{(1)}\odot
D_{(1)}\ar@{}[rd]|{\text{push}}\ar[d]_{\partial\odot1}\ar[r]^{1\odot\partial}&
C_{(1)}\odot D_{(0)}\ar[d]^\zeta\ar@/^20pt/[rdd]^{\partial\odot 1}&\\
C_{(0)}\odot D_{(1)}\ar[r]_\xi\ar@/_20pt/[drr]_{1\odot\partial}&C\ol{\odot} D\ar[rd]^{\ol{\partial}}&\\
&&C_{(0)}\odot D_{(0)}}
\end{equation}
in the category $\C{SG}$.
Here $\ol{\partial}$ is a pair in $\C{SG}$ for which we derive the
tensor product in $\C{qpm}$ by the functor $\Phi$ above, that is,
$$C\odot D=\Phi(\ol{\partial}\colon C\ol{\odot}D\r C_{(0)}\odot D_{(0)}),$$
is particular $(C\odot D)_{(0)}=C_{(0)}\odot D_{(0)}$ and $(C\odot
D)_{ee}=C_{ee}\otimes D_{ee}$. Moreover, notice that $\upsilon\zeta$ and $\upsilon\xi$ are both the identity on $ee$-groups.
The unit element
for this tensor product is the quadratic pair module $\overline{\Z}_\nill=\Phi(0\r\Z_\nill)$ given
by
\begin{equation}\label{el1.}
\xymatrix{&\Z\ar@{->>}[ld]&\\\Z/2\ar[rr]^0&&\Z\ar[lu]_H}
\end{equation}
where $H(n)=\binom{n}{2}$.

\section{The track category of quadratic pair modules}\label{tcqpm}

Given a quadratic pair module $C$ the group $C_0$ acts on $C_1$ by the formula, $x\in C_1$, $y\in C_0$,
$$x^y=x+P(\partial(x)|y)_H,$$ 
see Remark \ref{nse}.
We define tracks in $\C{qpm}$ as follows.

\begin{defn}\label{tracks}
A \emph{track} $\alpha\colon f\rr g$ between two morphisms
$f,g\colon C\r D$ in $\C{qpm}$ is a function
$$\alpha\colon C_{0}\To D_{1}$$
satisfying the equations, $x,y\in C_{0}$, $z\in C_{1}$,
\begin{enumerate}
\item $\alpha(x+y)=\alpha(x)^{f_{0}(y)}+\alpha(y)$,
\item $g_{0}(x)=f_{0}(x)+\partial\alpha(x)$,
\item $g_{1}(z)=f_{1}(z)+\alpha\partial(z)$.
\end{enumerate}
\end{defn}

These tracks are pulled back from the track category of crossed modules along the forgetful functor in Remark
\ref{nse}. The track structure for crossed modules is described in \cite{2hg1} 7. In particular we obtain the
following result.

\begin{prop}
The category is a $\C{qpm}$ track category.
\end{prop}

The vertical composition of tracks $\alpha$, $\beta$ is defined by addition
$(\alpha\vc\beta)(x)=\beta(x)+\alpha(x)$. The horizontal
composition of a track $\alpha$ and a map $f$, $g$ is defined as $(f\alpha)(x)=f_1\alpha(x)$ and $(\alpha
g)(x)=\alpha g_0(x)$. A trivial track $0^\vc_f\colon f\rr f$ is always defined as $0^\vc_f(x)=0$.

One can use the \emph{interval quadratic pair module} $\I$ to characterize tracks in $\C{qpm}$ and
$\C{qpm}$ in some cases. This quadratic pair module $\I$ is defined as follows. 
$$\I=\left(\begin{array}{c}\xymatrix{&\otimes^2\Z[i_0,i_1]\ar[ld]_P&\\
\Z\bar{\imath}\oplus\Z/2P(i_0|i_0)_H\oplus\Z P(i_0|i_1)_H
\ar[rr]_-\partial&&\grupo{i_0,i_1}_\nill\ar[lu]_H}\end{array}\right)$$
The quadratic map $H$ is defined as in (\ref{znil}).
The structure homomorphisms $P$ and $\partial$ are completely determined by the laws of a quadratic pair module and the equality
$\partial(\bar{\imath})=-i_0+i_1$. There are two obvious inclusions $i_0,i_1\colon\overline{\Z}_\nill\r\I$ and a projection
$p\colon\I\r\overline{\Z}_\nill$ defined by $p(i_0)=p(i_1)=1$ and $p(\bar{\imath})=0$.

\begin{lem}\label{cilin}
Let $f,g\colon C\r D$ be morphisms in $\C{qpm}$. Assume that $C$ is $0$-good. Then tracks $\alpha\colon f\rr g$ are in bijection with
morphisms $\bar{\alpha}\colon \I\odot C\r D$ with $f=\bar{\alpha} i_0$ and $g=\bar{\alpha} i_1$. The equivalence is given by the formula
$\alpha(c)=\bar{\alpha}\upsilon_e\zeta_e(\bar{\imath}\ca c)$ for $c\in C_0$. Here we use the square group
morphisms $\upsilon$ and $\zeta$ in (\ref{uni}) and (\ref{tp}).
\end{lem}

This lemma can be derived from the definition of the tensor product of square groups.



\begin{prop}\label{seext2}
The tensor product functor
$$\odot\colon\C{qpm}\times\C{qpm}\To\C{qpm}$$ is a track functor.
\end{prop}

\begin{proof}
Given tracks $\alpha\colon f\rr g\colon C\r D$ and $\beta\colon h\rr k\colon X\r Y$ in
$\C{qpm}$, the track $\alpha\odot\beta\colon f\odot h\rr g\odot k\colon C\odot X\r D\odot Y$ in $\C{qpm}$ with
$$\alpha\odot\beta\colon(C\odot X)_0=(C_{(0)}\odot X_{(0)})_e\To(D\odot Y)_1$$
is
defined as follows. Given $c\in C_0$ and $x\in X_0$
\begin{eqnarray*}
(\alpha\odot\beta)(c\cc x)&=&{\upsilon}_e\xi_e(f_0(c)\cc\beta(x))+{\upsilon}_e\zeta_e(\alpha(c)\cc
k_0(x))\\&&+(-f_0+g_0|f_0(c))_H\co Hk_0(x)\\
&\st{\text{(a)}}=&{\upsilon}_e\zeta_e(\alpha(c)\cc h_0(x))+{\upsilon}_e\xi_e(g_0(c)\cc
\beta(x))\\&&+(-f_0(c)+g_0(c)|f_0(c))_H\co Hh_0(x),
\end{eqnarray*}
and given $a\in C_{ee}$ and $b\in C_{ee}$
\begin{eqnarray*}
(\alpha\odot\beta)(a\co b)&=&-f_{ee}(a)\co h_{ee}(b)+g_{ee}(a)\co k_{ee}(b).
\end{eqnarray*}
Here we use the square group morphisms ${\upsilon}$, $\zeta$, and $\xi$ in (\ref{uni}) and (\ref{tp}). For
(a) we use 
\begin{eqnarray*}
{\upsilon}_e\zeta_e(\partial\alpha(c)\cc\beta(x))&=&{\upsilon}_e\xi_e(\alpha(c)\cc\partial\beta(x)).
\end{eqnarray*}
This equality follows from the fact that the square in (\ref{tp}) commutes.
We leave the reader to check that $\alpha\odot\beta$ is indeed a track $f\odot h\rr g\odot k$ and that the axioms
of a track functor are satisfied.
\end{proof}

The following commutativity property for the tensor product of tracks holds.

\begin{lem}
Given tracks $\alpha$, $\beta$ in $\C{qpm}$ the equation
$\tau_\odot(\alpha\odot\beta)=(\beta\odot\alpha)\tau_\odot$ holds.
\end{lem}

The proof is a straightforward but somewhat lengthy computation. One can also use the track
functor $\odot$ in Proposition \ref{seext2} to show that $\C{qpm}$ is indeed a symmetric monoidal $2$-category, compare
\cite{hdaI} and \cite{gcbsm2c}.

\part{Secondary homotopy groups as a lax symmetric monoidal functor}\label{II}

In this part we introduce the additive secondary homotopy group as a quadratic pair module and we formulate our main
results on the smash product for additive secondary homotopy groups leading to a lax symmetric monoidal functor.
We also give applications to unstable cup-one products, Toda brackets, and secondary Whitehead
products.

\section{Homotopy groups and secondary homotopy groups}\label{el1}

Let $\C{Ab}$ be the category of abelian groups. Using classical homotopy groups $\pi_nX$ we obtain for
$n\geq 0$ the functor 
$$\Pi_n\colon\C{Top}^*\To\C{Ab}$$
with
\begin{equation}\label{Clas}
\Pi_nX=\left\{
\begin{array}{ll}
\pi_nX, & n\geq2,\\
(\pi_1X)_\abb, & n=1,\\
\Z[\pi_0X], & n=0,
\end{array}\right.
\end{equation}
termed \emph{additive homotopy group}. 

One readily checks that the smash product
$f\wedge g\colon S^{n+m}\To X\wedge Y$
of maps $\set{f\colon S^n\r X}\in\pi_nX$ and $\set{g\colon S^m\r
Y}\in\pi_mY$ induces a well-defined homomorphism
\begin{equation}\label{ClasS}
\wedge\colon\Pi_nX\otimes\Pi_mY\To\Pi_{n+m}(X\wedge Y).
\end{equation}
This homomorphism is symmetric in the sense that the symmetry isomorphism
$\tau_\wedge\colon X\wedge Y\r Y\wedge X$ yields the equation in
$\Pi_{n+m}(Y\wedge X)$
\begin{equation}\label{inter}
(\tau_\wedge)_*(f\wedge g)=(-1)^{nm}g\wedge f.
\end{equation}
The sign $(-1)^{nm}$ is given by the degree of the symmetry isomorphism
\begin{equation}\label{inter2}
\tau_{n,m}=\tau_\wedge\colon S^{n+m}=S^n\wedge S^m\To S^m\wedge S^n=S^{m+n}.
\end{equation}
Here 
\begin{equation}\label{chuf}
\tau_{n,m}\in\sym{n+m}
\end{equation} 
is the shuffle permutation of $n+m$ elements which
exchanges the blocks $\set{1,\dots,n}$ and $\set{n+1,\dots,n+m}$. For this we recall that the symmetric group of $k$ letters
$\sym{k}$ acts on the $k$-sphere $$S^k=S^1\wedge\st{k}\cdots\wedge S^1$$ by permutation of coordinates according to
(\ref{permuconv}).

The main purpose of this paper is the generalization of the smash product operator (\ref{ClasS}) for
additive secondary homotopy groups.

\begin{defn}
Let $n\geq 2$. For a pointed space $X$ we define the \emph{additive
secondary homotopy group} $\Pi_{n,*}X$ which is the quadratic pair
module given by the diagram
$$\Pi_{n,*}X=\left(\begin{array}{c}\xymatrix@C=0pt{&**[r]\Pi_{n,ee}X=\otimes^2\Z[\L^nX]\ar[ld]_{P}&\\\Pi_{n,1}X\ar[rr]_\partial&&**[r]\Pi_{n,0}X=\grupo{\L^nX}_\nill\ar[lu]_H}\end{array}\right)$$
Here $\L^nX$ is the discrete pointed set of maps $S^n\r X$ in $\C{Top}^*$ and $H$ is defined as in (\ref{znil}). 

We describe the group $\Pi_{n,1}X$ and the homomorphisms $P$ and
$\partial$ as follows. The group $\Pi_{n,1}X$  is given by the set of equivalence classes $[f,F]$ 
represented
by a map $f\colon S^1\r \vee_{\L^n X}S^1$ and a track
$$\xymatrix{S^n\ar[r]_{\S^{n-1}f}^<(.98){\;\;\;\;\;}="a"\ar@/^25pt/[rr]^0_{}="b"&S^n_X\ar[r]_{ev}&X.\ar@{=>}"a";"b"_F}$$
Here the pointed space $$S^n_X=\vee_{\L^n X}S^n=\S^n\L^n X$$ is the $n$-fold suspension of the discrete pointed set
$\L^nX$, which is
the coproduct of $n$-spheres indexed by the set of non-trivial maps $S^n\r X$. The map $ev\colon S^n_X\r X$
is the obvious evaluation map. Given a map
$f\colon S^1\r \vee_{\L^n X}S^1$ we will denote $f_{ev}=ev(\S^{n-1}f)$, so that $F$ in the previous diagram is a track $F\colon f_{ev}\rr 0$.

The equivalence relation $[f,F]=[g,G]$ holds provided there is a 
track $$N\colon\Sigma^{n-1}f\rr\S^{n-1}g$$ with $\overline{\hopf}(N)=0$, see (\ref{elhopf}) and (\ref{sibar}) below, such that 
the composite track in the following diagram is the trivial track.
\begin{equation}\label{algo}
\xymatrix@C=50pt{S^n\ar@/^40pt/[rr]^0_{\;}="a"\ar@/_40pt/[rr]_0^{\;}="f"\ar@/^15pt/[r]|{\S^{n-1}f}^<(.935){\;}="b"_{\;}="c"\ar@/_15pt/[r]|{\S^{n-1}g}^{\;}="d"_<(.93){\;}="e"
&S^n_X\ar[r]^{ev}&X\ar@{=>}"a";"b"^{F^\vi}\ar@{=>}"c";"d"^{N}\ar@{=>}"e";"f"^G}
\end{equation}
That is $F=G\vc(ev\,N)$. The map
$\partial$ is defined by the  formula
$$\partial[f,F]=(\pi_1 f)_\nill(1),$$
where $1\in\pi_1S^1=\Z$.

The \emph{Hopf invariant} $\hopf(N)$ of a track $N\colon\Sigma^{n-1}f\rr\S^{n-1}g$ between maps as above is defined in \cite{2hg1} 3.3 by the
homomorphism 
\begin{equation}\label{elhopf}
H_2(IS^1,S^1\vee S^1)\st{ad(N)_*}\To H_2(\L^{n-1}S^n_X,\vee_{\L^ nX}S^1)\cong\left\{
\begin{array}{ll}
\hat{\otimes}^2\Z[\L^nX],& n\geq 3,\\&\\
{\otimes}^2\Z[\L^nX],& n=2.
\end{array}
\right.
\end{equation}
which carries the generator $1\in\Z\cong H_2(IS^1,S^1\vee S^1)$ to $\hopf(N)$. Here the isomorphism is induced by
the Pontrjagin product and $ad(N)_*$ is the homomorphism
induced in homology by the adjoint of $$S^{n-1}\wedge I_+\wedge S^1\st{(1\; 2)}\cong I_+\wedge S^n\st{N}\r S^n_X.$$ The
reduced tensor square $\hat{\otimes}^2$ in (\ref{elhopf}) is given by
$$\hat{\otimes}^2A=\frac{A\otimes A}{a\otimes b+b\otimes a\sim 0},$$ 
and $\bar{\sigma}\colon\otimes^2A\twoheadrightarrow\hat{\otimes}^2A$
is the natural projection. 
We define
\begin{equation}\label{sibar}
\overline{\hopf}=\left\{\begin{array}{ll}
\hopf,&\text{for $n\geq 3$},\\
&\\
\bar{\sigma}\hopf,&\text{for $n=2$}.
\end{array}\right.
\end{equation}
We refer the reader to \cite{2hg1} 3 for the elementary properties of $\hopf$ which will be used in this paper.

This completes the definition
of $\Pi_{n,1}X$, $n\geq 2$, as a set. The
group structure of $\Pi_{n,1}X$ is induced by the comultiplication $\mu\colon S^1\r S^1\vee S^1$, compare
\cite{2hg1} 4.4.

We now define the homomorphism $P$. Consider the diagram
$$\xymatrix{S^n\ar[rr]_{\S^{n-1}\beta}^{\;}="a"\ar@/^30pt/[rr]^0_{\;}="b"&&S^n\vee S^n\ar@{=>}"a";"b"_{B}}$$
where $\beta\colon S^1\r S^1\vee S^1$ is given such that $(\pi_1\beta)_\nill(1)=[a,b]\in\grupo{a,b}_\nill$ 
is the commutator of the generators. The track $B$ 
is any track with
$\overline{\hopf}(B)=-\bar{\sigma}(a\otimes b) \in\hat{\otimes}^2\Z[a,b]$. 
Given $x\otimes y\in\otimes^2\Z[\L^nX]$ let
$\tilde{x},\tilde{y}\colon S^1\r\vee_{\L^nX}S^1$ be maps with
$(\pi_1\tilde{x})_\abb(1)=x$ and $(\pi_1\tilde{y})_\abb(1)=y$. Then
the diagram
\begin{equation}\label{omedia}
\xymatrix{S^n\ar[rr]_{\S^{n-1}\beta}^{\;}="a"\ar@/^30pt/[rr]^0_{\;}="b"&&S^n\vee
S^n\ar@{=>}"a";"b"_{B}\ar[rr]_{\S^{n-1}(\tilde{y},\tilde{x})}&&S^n_X\ar[r]_{ev}&X}
\end{equation}
represents an element
\begin{equation*}
P(x\otimes y)=[(\tilde{y},\tilde{x})\beta,(\tilde{y}_{ev},\tilde{x}_{ev})B]\in\Pi_{n,1}X.
\end{equation*}
This completes the definition of the quadratic pair module $\Pi_{n,*}X$ for $n\geq 2$. For $n=0,1$ we define the
additive secondary homotopy groups $\Pi_{n,*}X$ by Remark \ref{dever} below. In this way we get for $n\geq 0$ a
functor
$$\Pi_{n,*}\colon\C{Top}^*\To\C{qpm}$$
which is actually a track functor.
\end{defn}

For a quadratic pair
module $C$ we denote $h_0C=\coker\partial$, $h_1C=\ker\partial$. There are natural
isomorphisms 
\begin{eqnarray}
\label{h0} h_0\Pi_{n,*}X&\cong&\pi_nX,\;\;n\geq2,\\
\nonumber h_1\Pi_{n,*}X&\cong&\pi_{n+1}X,\;\;n\geq 3,\\
\nonumber h_1\Pi_{2,*}X&\cong&\pi_{3}X/[\pi_2X,\pi_2X],
\end{eqnarray}
where $[-,-]$ denotes here the Whitehead product. Here we use \cite{2hg1} 5.1 and also \cite{2hg1} 6.11 for the case $n=2$.
Furthermore, the following property is crucial.

\begin{prop}\label{qlo}
The homomorphism $$h_0\Pi_{n,*}X\r h_1\Pi_{n,*}X\colon x\mapsto P(x|x)_H$$ coincides via (\ref{h0}) with the
homomorphism $\eta^*\colon\pi_nX\otimes\Z/2\r\pi_{n+1}X$ if $n\geq 3$ and
$\eta^*\colon\pi_2X\otimes\Z/2\r\pi_{3}X/[\pi_2X,\pi_2X]$ if $n=2$. Here $\eta^*$ is induced by precomposition with $\S^{n-2}\eta$
where $\eta\colon S^3\r
S^2$ is the Hopf map. 
\end{prop}

This follows from \cite{2hg1} 8.2.

\begin{rem}\label{dever}\label{remel}
Considering maps $f\colon S^n\r X$ together with tracks of
such maps to the trivial map, we introduced in \cite{2hg1} the secondary
homotopy group $\pi_{n,*}X$, which is a groupoid for $n=0$, a crossed module for $n=1$, a reduced quadratic
module for $n=2$, and a stable quadratic module for $n\geq 3$. 

Then using the adjoint functors $\ad_n$ of the forgetful functors $\phi_n$ as discussed in \cite{2hg1} 6 we get the 
\emph{additive secondary homotopy group} track functor
\begin{equation*}
\Pi_{n,*}\colon\C{Top}^*\To\C{squad}
\end{equation*}
given by
$$\Pi_{n,*}X=\left\{\begin{array}{ll}
\pi_{n,*}X,&\text{for }n\geq3,\\
\ad_3\pi_{2,*}X,&\text{for }n=2,\\
\ad_3\ad_2\pi_{1,*}X,&\text{for }n=1,\\
\ad_3\ad_2\ad_1\pi_{0,*}X,&\text{for }n=0.\\
\end{array}\right.$$
This is the secondary analogue of (\ref{Clas}).

Here the category $\C{squad}$ of stable quadratic modules
is not appropriate to study the smash product of
secondary homotopy groups since we do not have a symmetric monoidal structure in $\C{squad}$. Therefore 
we introduced above the category
$\C{qpm}$ of quadratic pair modules and we observe that $\Pi_{n,*}X$ in $\C{squad}$ yields a functor
to the category $\C{qpm}$ in the following way. As a stable quadratic module $\Pi_{n,*}X$ looks as follows
$$\otimes^2\Z[\L^nX]=\otimes^2(\Pi_{n,0}X)_\abb\st{P(-|-)_H}\To\Pi_{n,1}X\st{\partial}\To\Pi_{n,0}X=\grupo{\L^nX}_\nill.$$
In the quadratic pair module $\Pi_{n,*}X$ the quadratic map $H$ is defined as in (\ref{znil}) and 
the homomorphism $(-|-)_H\colon\otimes^2(\Pi_{n,0}X)_\abb=\Pi_{n,ee}X$ is the identity.
A map $f\colon X\r Y$ in $\C{Top}^*$ induces a homomorphism
$\Pi_{n,0}f\colon\Pi_{n,0}X\r\Pi_{n,0}Y$ between free nil-groups
which carries generators in $\Pi_{n,0}X$ to generators in
$\Pi_{n,0}Y$ and therefore $\Pi_{n,*}f$ is compatible with $H$. This shows that there is 
a canonical lift
\begin{equation*}
\xymatrix{&\C{qpm}\ar[d]\\\C{Top}^*\ar[ru]^{\Pi_{n,*}}\ar[r]_{\Pi_{n,*}}&\C{squad}}
\end{equation*}
with $\Pi_{n,(0)}X=\Z_\nill[\L^nX]$ for all $n\geq 0$, in particular $\Pi_{n,*}X$ is always a $0$-good quadratic
pair module.
Compare \cite{2hg2} 1.15.

The definition of $\Pi_{2,*}X$ given above coincides with the lifting of $\ad_3\pi_{2,*}X$ to $\C{qpm}$ by the
claim (*) in the proof of \cite{2hg1} 4.9.

Generalizing (\ref{h0}) we have for all $n\geq 0$ a natural isomorphism
\begin{eqnarray}
h_0\Pi_{n,*}X&\cong\Pi_nX,
\end{eqnarray}
see \cite{2hg1} 6.10.
\end{rem}

In this paper we are concerned with the properties of the track
functor $\Pi_{n,*}$ mapping to the category $\C{qpm}$. In order to simplify notation given a map $f\colon X\r Y$ in $\C{Top}^*$ we will just denote
$$f_*=\Pi_{n,i}f\colon\Pi_{n,i}X\To\Pi_{n,i}Y,\;\;\text{for $i=0,1,ee$ and $n\geq 0$}.$$
Moreover, given a track $\alpha\colon f\rr g$ between maps $f,g\colon X\r Y$ we denote by
$$\alpha_*=\Pi_{n,*}\alpha\colon\Pi_{n,0}X\To\Pi_{n,1}Y,\;\;n\geq 0,$$
the induced track $\alpha_*\colon f_*\rr g_*$ in $\C{qpm}$.

\section{Smash product for secondary homotopy groups}\label{LA}

In this section we describe our main results connecting the tensor product of quadratic pair modules and the smash
product of pointed spaces. The smash product operator in the next theorem is the canonical analogue of the smash
product for classical homotopy groups in (\ref{ClasS}) above.

\begin{thm}\label{lamel}
The functor of additive secondary homotopy groups
$\Pi_{n,*}\colon\C{Top}^*\r\C{qpm}$ admits a well-defined smash
product operator
$$\wedge\colon\Pi_{n,*}X\odot\Pi_{m,*}Y\To\Pi_{n+m,*}(X\wedge Y),$$
which is a morphism in $\C{qpm}$, inducing the smash product of classical homotopy groups in (\ref{ClasS}) on
$h_0$.
This operator is natural in X and Y with respect to maps and tracks.
\end{thm}

This smash product operator is given in the $(0)$-level by the following morphism of square groups
$$\Z_\nill[\L^nX]\odot\Z_\nill[\L^mY]\cong\Z_\nill[\L^nX\wedge \L^mY]\st{\Z_\nill[\wedge]}\To\Z_\nill[\L^{n+m}(X\wedge Y)].$$
Here \begin{equation*}
\wedge\colon(\L^nX)\wedge(\L^m Y)\To\L^{m+n}(X\wedge
Y),\;\;n,m\geq0,
\end{equation*}
is the map between discrete pointed sets defined by
$$(f\colon S^n \r X)\wedge (g\colon S^m\r Y)\mapsto (f\wedge g\colon
S^{n+m}\r X\wedge Y).$$
On the $(1)$-level the definition of the smash product operator in Theorem \ref{lamel} is more elaborate, see
Definition \ref{ka} below. 

\begin{proof}[Proof of Theorem \ref{lamel}]
The first part of Theorem \ref{lamel} follows from Lemmas \ref{P1}, \ref{P2}, \ref{P3} and \ref{Pinf} in Part \ref{III}. 
For the naturality we use Lemma \ref{inri} in Part \ref{III} and Theorems \ref{lmf} and \ref{nosym}.
\end{proof}

We will use the following notation for the image of generators in the tensor product by the smash product
morphism in Theorem \ref{lamel}. Given $x\in\Pi_{n,i}X$ and $y\in\Pi_{m,j}Y$ for $0\leq i,j,i+j\leq 1$ we denote
by
$$x\wa y\in\Pi_{n+m,i+j}(X\wedge Y)$$
to the image by $\wedge$ of the element $x\ca y\in\Pi_{n+m,0}(X\wedge Y)$ if $i=j=0$, of
${\upsilon}_e\zeta_e(x\ca y)\in\Pi_{n+m,1}(X\wedge Y)$ if $i=1$ and $j=0$, and of 
${\upsilon}_e\xi_e(x\ca y)\in\Pi_{n+m,1}(X\wedge Y)$ if $i=0$ and $j=1$. Here we use the square group morphisms
${\upsilon}$, $\zeta$ and $\xi$ in (\ref{uni}) and (\ref{tp}). Similarly for $x\cc y$ and $x\wc y$. Moreover,
given $a\in\Pi_{n,ee}X$ and $b\in\Pi_{m,ee}Y$ we denote by $$a\wedge b\in\Pi_{n+m,ee}(X\wedge Y)$$
to the image of $a\otimes b\in(\Pi_{n,*}X\odot\Pi_{m,*}Y)_{ee}$ by $\wedge$.

\begin{thm}\label{lmf}
The smash product operator endows $\Pi_{*,*}$ with the structure of
a lax monoidal functor from $\C{Top}^*$ to the category of graded $0$-good quadratic pair modules. That is, the following
diagram commutes
$$\xymatrix{\Pi_{n,*}X\odot(\Pi_{m,*}Y\odot\Pi_{l,*}Z)\ar[dd]^{1\odot\wedge}\ar[r]^\cong&(\Pi_{n,*}X\odot\Pi_{m,*}Y)\odot\Pi_{l,*}Z
\ar[d]^{\wedge\odot 1}\\
&\Pi_{n+m,*}(X\wedge
Y)\odot\Pi_{l,*}Z\ar[d]^\wedge\\
\Pi_{n,*}X\odot\Pi_{m+l,*}(Y\wedge
Z)\ar[r]^\wedge&\Pi_{n+m+l,*}(X\wedge Y\wedge Z)}$$ and for the unit
$S^0$ of the symmetric monoidal category $(\C{Top}^*,\wedge)$ and the unit
$\overline{\Z}_\nill$ of $(\C{qpm},\odot)$ there is an isomorphism
$u\colon\overline{\Z}_\nill\cong\Pi_{0,*}S^0$ such that the following diagrams commute.
$$\xymatrix{\Pi_{0,*}S^0\odot\Pi_{n,*}X\ar[r]^\wedge&\Pi_{n,*}(S^0\wedge X)\ar[d]^\cong\\
\overline{\Z}_\nill\odot\Pi_{n,*}X\ar[r]^\cong\ar[u]^{u\odot 1}_\cong&\Pi_{n,*}X}$$
$$\xymatrix{\Pi_{n,*}X\odot\Pi_{0,*}S^0\ar[r]^\wedge&\Pi_{n,*}(X\wedge S^0)\ar[d]^\cong\\
\Pi_{n,*}X\odot \overline{\Z}_\nill\ar[r]^\cong\ar[u]^{1\odot u}_\cong&\Pi_{n,*}X}$$
\end{thm}

The isomorphism $u\colon\overline{\Z}_\nill\cong\Pi_{0,*}S^0$ is the unique one sending
$1\in\Z=(\overline{\Z}_\nill)_0$ to $u(1)=1_{S^0}\colon S^0\r S^0$ in $\Pi_{0,0}S^0$. 

\begin{proof}[Proof of Theorem \ref{lmf}]
In dimensions $\geq 1$ the
associativity property in Theorem \ref{lmf} follows from Lemma \ref{P4} in Part \ref{III}. In case dimension $0$
is involved we use the more algebraic Lemma
\ref{Pinf}. The commutativity of the squares with the isomorphism $u$ is easy to check. This is left to the
reader.
\end{proof}

The graded commutativity equation (\ref{inter}) for the smash product for classical homotopy groups has a secondary
analogue as follows.

\begin{thm}\label{nosym}
The following diagram commutes in
$\C{qpm}$.
$$\xymatrix{\Pi_{n,*}X\odot\Pi_{m,*}Y\ar[dd]^\cong_{\tau_\odot}\ar[r]^{\wedge}&\Pi_{n+m,*}(X\wedge
Y)\ar[d]_\cong^{(\tau_{X,Y})_*}\\&\Pi_{n+m,*}(Y\wedge X)\ar@{<-}[d]_\cong^{\tau_{n,m}^*}\\
\Pi_{m,*}Y\odot\Pi_{n,*}X\ar[r]^{\wedge}&\Pi_{m+n,*}(Y\wedge X)}$$
Here $\tau_\odot$ is the symmetry isomorphism of
the tensor product $\odot$ in $\C{qpm}$, and $\tau_{n,m}^*$ is given by the right action of the shuffle
permutation $\tau_{n,m}\in\sym{n+m}$ in (\ref{chuf}) on
$\Pi_{m+n,*}(Y\wedge X)$, see the next section.
\end{thm}

This follows from Lemma \ref{P5} in Part \ref{III}.

\section{The symmetric action on smash products}

Secondary homotopy groups, regarded
as a functor from pointed spaces to graded $0$-good quadratic pair modules
$$\Pi_{*,*}\colon\C{Top}^*\To\C{qpm}_0^\N$$
is a lax monoidal functor, see Theorem \ref{lmf}. The monoidal structure in $\C{qpm}^\N_0$ is
the usual graded extension of the tensor product $\odot$ of
quadratic pair modules, see the appendix, and $\C{Top}^*$ has the monoidal structure
given by the smash product $\wedge$. Both monoidal structures are
symmetric, however Theorem \ref{nosym} shows that $\Pi_{*,*}$ is not lax symmetric monoidal since the action of
the shuffle permutation
$\tau_{n,m}$ in (\ref{chuf}) is involved. This inconvenience is solved by enriching the
structure of secondary homotopy groups with the sign action of the
symmetric track groups constructed in \cite{2hg2}.


\begin{defn}\label{signg}
Let $\set{\pm1}$ be the multiplicative group of order $2$. A
\emph{sign group} $G_\vc$ is a diagram of group homomorphisms
$$\set{\pm1}\st{\imath}\hookrightarrow G_\vc\st{\partial}\twoheadrightarrow G\st{\varepsilon}\To \set{\pm1}$$
where the first two morphisms $\imath$ and $\partial$ form an extension. Here all groups
have a multiplicative group law and the composite
$\varepsilon\partial$ is also denoted by $\varepsilon\colon G_\vc\r
\set{\pm1}$.

A morphism $f_\vc\colon G_\vc\r K_\vc$ of sign groups is a
commutative diagram of group homomorphisms
$$\xymatrix{\set{\pm1}\ar@{=}[d]\ar@{^{(}->}[r]&G_\vc\ar[d]^{f_\vc}\ar@{->>}[r]&G\ar[d]^f\ar[r]&\set{\pm1}\ar@{=}[d]\\
\set{\pm1}\ar@{^{(}->}[r]&K_\vc\ar@{->>}[r]&K\ar[r]&\set{\pm1}}$$
This defines the category $\C{Gr}_\pm$ of sign groups. The initial object $1_\vc$ in this category given by $G=\set{1}$ will be termed the
\emph{trivial sign group}.
\end{defn}

\begin{rem}\label{cmes}
Recall from \cite{2hg2} 3.6 that a sign group $G_\vc$ gives rise to
a crossed module
$$\partial_\vc=(\varepsilon,\partial)\colon G_\vc\To\set{\pm1}\times G.$$
where $\set{\pm1}\times G$ acts  on $G_\vc$ by the formula
$$g^{(x,h)}=\bar{h}^{-1}g\bar{h}\imath\left(\varepsilon(g)^{\binom{x}{2}}\right).$$
Here $g\in G_\vc$, $x\in\set{\pm1}$, $h\in G$, and
$\bar{h}\in G_\vc$ is any element with $\partial(\bar{h})=h$.
This is a well-defined crossed module since $G_\vc$ is a central extension of $G$ by $\set{\pm1}$.
\end{rem}

The main examples of sign groups are the symmetric track groups
$$\set{\pm1}\hookrightarrow\symt{n}\st{\delta}\twoheadrightarrow\sym{n}\st{\sign}\To\set{\pm1}$$
defined as follows. The symmetric group $\sym{n}$ acts on the
left of the $n$-sphere $$S^n=S^1\wedge\st{n}\cdots\wedge S^1,$$ see (\ref{permuconv}). The
elements of the symmetric track group $\symt{n}$ for $n\geq 2$ are tracks
$\alpha\colon \sigma\rr(\cdot)_n^{\sign(\sigma)}$ between maps
$\sigma,(\cdot)_n^{\sign(\sigma)}\colon S^n\r S^n$, with
$\sigma\in\sym{n}$ and
$(\cdot)_n^{\sign(\sigma)}=\S^{n-1}(\cdot)^{\sign(\sigma)}$, where
$$(\cdot)^k\colon S^1\To S^1\colon z\mapsto z^k,\;\;k\in\Z,$$
is given by the (multiplicative) topological abelian group structure
of $S^1$. The group law in $\symt{n}$ is given by the horizontal
composition of tracks. For $n=0,1$, let $\symt{n}$ be the trivial sign group. Compare \cite{2hg2} 5 and 6.

The smash product $S^m\wedge -$ induces a sign group morphism
\begin{equation}\label{s-}
S^m\wedge-\colon\symt{n}\To\symt{m+n}
\end{equation} sending
a track $\alpha\colon \sigma\rr(\cdot)_n^{\sign(\sigma)}$ as above
to $S^m\wedge \alpha$. This morphism is given on symmetric groups by
the usual inclusion $\sym{n}\subset\sym{m+n}$ obtained by regarding
$\sym{n}$ as the subgroup of permutations of $m+n$ elements fixing
the first $m$ elements.



One can not directly define a sign group morphism $$-\wedge
S^n\colon\symt{m}\r\symt{m+n}$$ in a geometric way as above since
$(\cdot)^k_m\wedge S^n\neq (\cdot)^k_{m+n}$. With the help of the
crossed module structure for sign groups described in Remark
\ref{cmes} and the shuffle permutation $\tau_{n,m}$ in (\ref{chuf}) we define $-\wedge S^n$ as the following composite
\begin{equation}\label{-s}
-\wedge S^n\colon\symt{m}\st{S^n\wedge
-}\To\symt{n+m}\st{(\cdot)^{(1,\tau_{n,m})}}\To\symt{m+n}
\end{equation} 
This morphism is given on symmetric groups by
the inclusion $\sym{m}\subset\sym{m+n}$ obtained by identifying
$\sym{m}$ with the group of permutations of $m+n$ elements which
only permute the first $m$ ones.

\begin{defn}
A \emph{twisted bilinear morphism} of sign groups
$$(f_\vc,g_\vc)\colon G_\vc\times L_\vc\To K_\vc$$
is given by a pair of sign group morphisms
$f_\vc\colon G_\vc\r K_\vc$, $g_\vc\colon L_\vc\r K_\vc$,
such that given $a\in G$ and $b\in L$ the equality
$$f(a)g(b)=g(b)f(a)$$ holds in $K$, and therefore the group homomorphism
$$(f,g)\colon G\times L\To K\colon (a,b)\mapsto f(a)g(b)$$
is defined, and given $x\in G_\vc$, $y\in L_\vc$ the following
equality is satisfied in $K_\vc$
$$f_\vc(x)g_\vc(y)=g_\vc(y)f_\vc(x)\imath\left((-1)^{\binom{\varepsilon(x)}{2}\binom{\varepsilon(y)}{2}}\right).$$

The \emph{twisted product} $G_\vc\tilde{\times}L_\vc$ of sign groups
$G_\vc$, $L_\vc$ is a sign group
$$\set{\pm1}\st{\imath}\hookrightarrow G_\vc\tilde{\times}L_\vc\twoheadrightarrow G\times L\st{\varepsilon}\twoheadrightarrow\set{\pm1}$$
together with a universal twisted bilinear morphism
\begin{equation}\label{universal}
(i_{G_\vc},i_{L_\vc})\colon G_\vc\times L_\vc\To
G_\vc\tilde{\times}L_\vc.
\end{equation}
The group $G_\vc\tilde{\times}L_\vc$ is generated by the symbols
$\overline{s}$, $\overline{t}$ and $\omega$, for $t\in G_\vc$ and
$s\in L_\vc$, with the following relations:
\begin{enumerate}
\item $\omega$ is central,
\item $\overline{\imath(-1)}=\omega$ for both
$\imath\colon\set{\pm1}\hookrightarrow G_\vc$ and
$\imath\colon\set{\pm1}\hookrightarrow L_\vc$,
\item $\overline{r_1r_2}=\overline{r_1}\cdot\overline{r_2}$ when
$r_1,r_2$ lie both in the group $G_\vc$ or both in $L_\vc$,
\item $\overline{t}\overline{s}=\overline{s}\overline{t}\omega^{\binom{\varepsilon(t)}{2}\binom{\varepsilon(s)}{2}}$ for
$t\in G_\vc$ and $s\in L_\vc$.
\end{enumerate}
The homomorphism $\partial\colon
G_\vc\tilde{\times}L_\vc\twoheadrightarrow G\times L$ is defined by
$\partial(\overline{t})=(\partial(t),0)$ for $t\in G_\vc$, $\partial(\overline{s})=(0,\partial(s))$ for 
$s\in L_\vc$,
and $\partial(\omega)=0$. The universal bilinear morphism is given
by $i_{G_\vc}(t)=\overline{t}$ and $i_{L_\vc}(s)=\overline{s}$.
The twisted product is a non-symmetric monoidal structure in the
category $\C{Gr}_\pm$ of sign groups where rhe unit is the trivial sign group $1_\vc$.
\end{defn}

\begin{prop}\label{inclusym}
The morphisms in (\ref{s-}) and (\ref{-s}) induce a morphism of sign
groups
$$\symt{m}\tilde{\times}\symt{n}\To\symt{m+n}.$$
\end{prop}

This proposition can be derived easily from the presentation of the symmetric track groups given in \cite{2hg2}
6.11.

We now introduce the action of a sign group on a quadratic pair module. In \cite{2hg2} we show that the sign
group $\symt{n}$ acts in this sense on $\Pi_{n,*}X$.

\begin{defn}\label{laac}
A sign group $G_\vc$ \emph{acts on the right} of a quadratic pair module
$C$ if $G$ acts on the right of $C$ by morphisms $g^*\colon C\r
C$, $g\in G$, in $\C{qpm}$, and there is a bracket
$$\grupo{-,-}=\grupo{-,-}_{G}\colon C_0\times G_\vc\To C_1$$
satisfying the following properties, $x,y\in C_0$, $z\in C_1$,
$s,t\in G_\vc$.
\begin{enumerate}
\item
$\grupo{x+y,t}=\grupo{x,t}+\grupo{y,t}+P(-{\partial(t)^*(x)}+\varepsilon(t)^*(x)|{\partial(t)^*(y)})_H$,
\item $\varepsilon(t)^*(x)=\partial(t)^*(x)+\partial\grupo{x,t}$,
\item $\varepsilon(t)^*(z)=\partial(t)^*(z)+\grupo{\partial(z),t}$,
\item $\grupo{x,s\cdot
t}=\grupo{\partial(s)^*(x),t}+\grupo{\varepsilon(t)^*(x),s}$,
\item for the element $\omega=\imath(-1)\in G_\vc$ we have the
\emph{$\omega$-formula}:
$$\grupo{x,\omega}=P(x|x)_H.$$
Notice that the $\omega$-formula corresponds to the homomorphism in Proposition \ref{qlo}.
\end{enumerate}
Here we use the notation
\begin{itemize}
\item $(-1)^*(x)=-x+\partial PH(x)$,
\item $(-1)^*(z)=-y+HP\partial(y)$,
\item $1^*=$ identity,
\end{itemize}
introduced in \cite{2hg2}. Notice that $(-1)^*(-1)^*=1^*$.
The trivial sign group acts on any quadratic pair module in a unique way.
\end{defn}

This definition of a sign group action, as given in \cite{2hg2}, can be reinterpreted in terms of the tensor
product $\odot$ of quadratic pair modules by using the following ``group ring'' construction for sign groups.

\begin{defn}\label{gring}
A \emph{quadratic pair algebra} $R$ is a monoid in the category $\C{qpm}$ of quadratic pair modules. The image of
a right-linear generator $r\cc s$ in the tensor product by the monoid structure morphism $R\odot R\r R$ will be
denoted by $r\cdot s$. This notation will also be used below (in the proof of Lemma \ref{rein}) for right modules over a quadratic pair algebra.
Given a sign group $G_\vc$ the quadratic pair algebra $A(G_\vc)$ has generators
\begin{itemize}
\item $[g]$ for any $g\in G$ on the $0$-level,
\item $[t]$ for any $t\in G_\vc$ on the $1$-level,
\item no generators on the $ee$-level,
\end{itemize}
and relations
\begin{itemize}
\item $H[g]=0$ for $g\in G$,
\item $[1]=1$ the unit element,
\item $[gh]=[g]\cdot[h]$ for $g,h\in G$,
\item $\partial[t]=-[\partial(t)]+\varepsilon(t)$,
\item $[st]=[\partial(s)]\cdot[t]+\varepsilon(t)\cdot[s]
+\binom{\varepsilon(s)}{2}\binom{\varepsilon(t)}{2}P(1|1)_H$ for $s,t\in G_\vc$,
\item $[\omega]=P(1|1)_H$.
\end{itemize}
In these equations $-1$ can appear as a value of the homomorphism $\varepsilon$. This
$-1$ denotes the additive opposite of the unit element $1\in A_{0}(G_\vc)$, 
except when it appears as part of a cominatorial number, where it is regarded just as an integer. 
The relations above show that $A_{(0)}(G_\vc)=\Z_\nill[G_+]$ where $G_+$ is the group $G$ 
together with an outer base point, so $A(G_\vc)$ is $0$-good. 
\end{defn}

The ``group ring'' of a sign group defines a functor onto the category $\C{qpa}$ of quadratic pair algebras
$$A\colon\C{Gr}_\pm\To\C{qpa}.$$

\begin{prop}\label{essmon}
The functor $A$ is strict monoidal.
\end{prop}

\begin{proof}
The isomorphism $A(G_\vc)\odot A(K_\vc)\cong A(G_\vc\tilde{\times} K_\vc)$
is given by the composite
$$A(G_\vc)\odot A(K_\vc)\st{A(i_{G_\vc})\odot A(i_{K_\vc})}\To A(G_\vc\tilde{\times} K_\vc)\odot A(G_\vc\tilde{\times} K_\vc)
\st{\text{product}}\To A(G_\vc\tilde{\times} K_\vc).$$
For the trivial sign group $1_\vc$ the isomorphism $\Z_\nill\cong A(1_\vc)$ is the unit of the quadratic pair algebra $A(1_\vc)$.
\end{proof}

The following lemma gives a reinterpretation of sign group actions in terms of algebras and modules in the
monoidal category $\C{qpm}$.

\begin{lem}\label{rein}
Let $G_\vc$ be a sign group and let $C$ be a quadratic pair module. A sign action of $G_\vc$ on $C$ in the sense
of Definition \ref{laac} corresponds to a right $A(G_\vc)$-module structure on $C$.
\end{lem}

\begin{proof}
With the notation in in Definitions \ref{laac} and \ref{gring} the correspondence is given by the formulas, $g\in G$, $t\in G_\vc$,
\begin{eqnarray*}
g^*x&=&x\cdot[g],\\
\grupo{x,t}&=&x\cdot[t]+\binom{\varepsilon(t)}{2}PH(x).
\end{eqnarray*}
The technical details of this proof are left to the reader.
\end{proof}

In \cite{2hg2} we define a sign group action of $\symt{n}$ on $\Pi_{n,*}X$,
hence combining Proposition \ref{essmon} and Lemma \ref{rein} we readily obtain the following result.

\begin{thm}\label{actdi}
The sign group $\symt{n}\tilde{\times}\symt{m}$ acts on
$\Pi_{n,*}X\odot\Pi_{m,*}Y$.
\end{thm}

We will now consider the compatibility of the smash product operation in Theorem \ref{lamel} with the sign group
actions. 

\begin{thm}\label{yi}
The smash product morphism
$$\wedge\colon\Pi_{n,*}X\odot\Pi_{m,*}Y\To\Pi_{n+m,*}(X\wedge Y)$$
in Theorem \ref{lamel} is equivariant with respect to the action of
$\symt{n}\tilde{\times}\symt{m}$ on $\Pi_{n,*}X\odot\Pi_{m,*}Y$
defined by Theorem \ref{actdi} and the sign group morphism
$$\symt{n}\tilde{\times}\symt{m}\To\symt{n+n}$$
in Proposition \ref{inclusym}.
\end{thm}

This theorem follows from Lemma \ref{P7} in Part \ref{III}.

Since the secondary homotopy groups $\Pi_{n,*}X$ have a canonical action of the sign group $\symt{n}$ we are led
to consider the following category of symmetric sequences in $\C{qpm}$ (this is similar to the treatment of
symmetric spectra in \cite{se}).

\begin{defn}
An object $X$ in the category $\C{qpm}_0^{\symtt}$ of symmetric
sequences is a sequence of $0$-good
quadratic pair modules $X_n$ endowed with a sign group action of the
symmetric track group $\symt{n}$, $n\geq 0$. A morphism $f\colon X\r
Y$ in $\C{qpm}_0^{\symtt}$ is a sequence of $\symt{n}$-equivariant
morphisms $f_n\colon X_n\r Y_n$ in $\C{qpm}_0$. The
results in \cite{2hg2} show that secondary homotopy groups yield a
functor $$\Pi_{*,*}\colon\C{Top}^*\To\C{qpm}_0^{\symtt}.$$ The
category $\C{qpm}_0^{\symtt}$ has a symmetric monoidal structure
denoted by $\odot_{\symtt}$. The tensor product $X\odot_{\symtt} Y$
of two symmetric sequences $X$, $Y$ of $0$-good quadratic pair modules is characterized by the
following universal property: for any symmetric sequence $Z$ of quadratic pair modules
there is a natural bijection
$$\hom_{\C{qpm}_0^{\symtt}}(X\odot_{\symtt}
Y,Z)\cong\prod_{p,q\in\N}\hom_{\symt{p}\tilde{\times}\symt{q}}(X_p\odot
Y_q, Z_{p+q}).$$ Here $\hom_{\symt{p}\tilde{\times}\symt{q}}$
denotes the set of morphisms in $\C{qpm}$ which are equivariant with
respect to the sign group morphism in Proposition \ref{inclusym}. The explicit construction of $X\odot_{\symtt}
Y$ is indicated in the appendix.
The symmetry isomorphism
$$X\odot_{\symtt}Y\cong Y\odot_{\symtt}X$$
is induced by the morphisms in $\C{qpm}$
$$X_p\odot Y_q\cong Y_q\odot X_p\To (Y\odot_{\symtt} X)_{q+p}\st{\tau_{p,q}^*}\To(Y\odot_{\symtt} X)_{p+q}.$$
Here the first morphism is the symmetry isomorphism for $\odot$, the
second one is induced by the universal property of $Y\odot_{\symtt}
X$ and in the third morphism we use the sign group action of $\symt{p+q}$ and the shuffle permutation
$\tau_{p,q}\in\sym{p+q}$ in (\ref{chuf}). The associativity isomorphism is defined by using the universal property of the
$3$-fold tensor product, which is analogous to the $2$-fold case above. The unit element is $\overline{\Z}_\nill$
concentrated in degree $0$. 
\end{defn}

Now Theorems \ref{lamel}, \ref{yi}, \ref{lmf} and \ref{nosym} can be restated as follows.

\begin{thm}\label{yii}
The smash product operator induces a natural
morphism in the category $\C{qpm}_0^{\symtt}$
$$\wedge\colon\Pi_{*,*}X\odot_{\symtt}\Pi_{*,*}Y\To\Pi_{*,*}(X\wedge Y)$$
which is compatible with the associativity, commutativity and unit
isomorphisms for the symmetric monoidal structures $\wedge $ and
$\odot_{\symtt}$ in $\C{Top}^*$ and $\C{qpm}_0^{\symtt}$,
respectively. Equivalently the functor
$$\Pi_{*,*}\colon\C{Top}^*\To\C{qpm}_0^{\symtt}$$
given by secondary homotopy groups is lax symmetric monoidal.
\end{thm}

\section{Secondary Whitehead products}\label{swp}

The smash product may be used for the definition of Whitehead products in ordinary homotopy groups. In fact, any
path connected space $X$ is homotopy equivalent to the classifying space of a topological group $G$ so that
$\pi_nG=\pi_{n+1}X$. We consider the additive homotopy groups $\Pi_nG$ which satisfy 
$$\Pi_nG=\pi_nG\text{ for }n\geq1.$$ Using the smash product operator $\wedge$ in (\ref{ClasS}) for the functor $\Pi_*$ and 
the commutator map $c\colon G\wedge G\r G$ with $c(a\wedge b)=a^{-1}b^{-1}ab$ we obtain the composite
$$[-,-]\colon\Pi_*G\otimes\Pi_*G\st{\wedge}\To\Pi_*(G\wedge G)\st{\Pi_*c}\To\Pi_*G$$
which corresponds to the Whitehead product in $\pi_*X$. It is well known that $(\Pi_*G,[-,-])$ has the structure
of a graded Lie algebra if $X$ is simply connected.

In a similar way we now define the \emph{secondary Whitehead product} for the additive secondary homotopy groups
$\Pi_{n,*}G$ by the composite 
$$[-,-]\colon\Pi_{*,*}G\otimes\Pi_{*,*}G\st{\wedge}\To\Pi_{*,*}(G\wedge G)\st{\Pi_{*,*}c}\To\Pi_{*,*}G.$$

Marcum defines in \cite{ombc} the partial Whitehead product of a map $\alpha$ and a track $F$ as in the following
diagram
$$\xymatrix{\S A\ar[r]_\alpha&X&B\ar[l]^\beta&
\tilde{B}.\ar[l]^e_<(1){\;}="a"\ar@/_30pt/[ll]_0^{\;}="b"\ar@{=>}"a";"b"_F}$$
Marcum's partial Whitehead product lives in the group of homotopy classes $[\S A\wedge \tilde{B},X]$. It can be
obtained from the secondary Whitehead product for additive secondary homotopy groups in case $A$ and $\tilde{B}$
are spheres. 

We will explore this connection in a sequel of this paper where we shall discuss the algebraic properties of the structure $(\Pi_{*,*}G,[-,-])$ which
leads to the notion of a secondary Lie algebra. This should be compared with the notion of secondary Hopf
algebra discussed in \cite{asco}.

\section{Cup-one products}\label{c1}

Let $n\geq m> 1$ be even integers. The cup-one product operation
$$\pi_nS^m\To\pi_{2n+1}S^{2m}\colon\alpha\mapsto\alpha\cup_1\alpha=Sq_1(\alpha)$$
is defined in the following way, compare \cite{c1Tb} 2.2.1. Let $k$ be any positive integer and let 
$\tau_k\in\sym{2k}$ be the permutation exchanging the first and
the second block of $k$ elements in $\set{1,\dots,2k}$. If $k$ is even then $\sign \tau_k=1$. 
We choose for any even integer $k>1$ a track
$\hat{\tau}_k\colon\tau_k \rr1_{S^{2k}}$ in $\symt{2k}$. Consider the following diagram in the track category
$\C{Top}^*$ of pointed spaces where 
$a\colon S^n\r S^m$ represents $\alpha$.
\begin{equation}\label{cup1}
\xymatrix{S^{2n}\ar[r]^{a\wedge a}\ar[d]^{{\tau}_n}_{\;}="b"\ar@/_30pt/[d]_{1_{S^{2n}}}^{\;}="a"&
S^{2m}\ar[d]_{{\tau}_m}^{\;}="c"\ar@/^30pt/[d]^{1_{S^{2m}}}_{\;}="d"\\
S^{2n}\ar[r]^{a\wedge a}&S^{2m}\ar@{=>}"a";"b"^{\hat{\tau}_n^\vi}\ar@{=>}"c";"d"^{\hat{\tau}_m}}
\end{equation}
By pasting this diagram we obtain a self-track of $a\wedge a$
\begin{equation}\label{cup2}
(\hat{\tau}_m(a\wedge a))\vc((a\wedge a)\hat{\tau}_n^\vi)\colon a\wedge a\rr a\wedge a.
\end{equation}
The set of self-tracks $a\wedge a\rr a\wedge a$ is the automorphism group of the map $a\wedge a$ in the track
category $\C{Top}^*$. 
The element $\alpha\cup_1\alpha\in\pi_{2n+1}S^{2m}$ is given by the track (\ref{cup2})
via the well-known Barcus-Barratt-Rutter isomorphism
$$\aut(a\wedge a)\cong\pi_{2n+1}S^{2m},$$
see \cite{hcefm}, \cite{hcmifs} and also \cite{ch4c} VI.3.12 and \cite{slort} for further details.

The following proposition yields a description of the cup-one product in terms of the structure of additive
secondary homotopy groups.

\begin{prop}
Let $n$ and $m$ be even positive integers. For $\alpha\in\Pi_nS^m$ we choose $a\in\Pi_{n,0}S^m$ representing
$\alpha$ and we define in $\Pi_{2n,1}S^{2m}$
$$Sq_1(\alpha)=-\grupo{a\wa a,\hat{\tau}_{n}}+(\hat{\tau}_{m})_*(a\wc a)-P(H(a)\wedge TH(a)).$$
Then $\partial Sq_1(\alpha)=0$ so that $Sq_1(\alpha)\in h_1\Pi_{2n,1}S^{2m}=\pi_{2n+1}S^{2m}$. Moreover,
$Sq_1(\alpha)=\alpha\cup_1\alpha$.
\end{prop}

\begin{proof}
The track $\hat{\tau}_{m}\colon\tau_{m}\rr 1$ induces a track
$(\hat{\tau}_{m})_*\colon(\tau_{m})_*\rr 1$ in $\C{qpm}$ satisfying
$$\partial(\hat{\tau}_{m})_*(a\wc a)=-(\tau_m)_*(a\wc a)+a\wc a.$$
By the symmetric action we have the element $\grupo{a\wa a,\hat{\tau}_{n}}\in\Pi_{2n,1}S^{2m}$ satisfying
$$-\partial\grupo{a\wa a,\hat{\tau}_{n}}=-a\wa a+\tau_{n}^*(a\wa a).$$
Hence we get
$$\partial Sq_1(\alpha)=-a\wa a+(\tau_{n})^*(a\wa a)-(\tau_{m})_*(a\wc a)+a\wc a-\partial 
P(H(a)\wedge TH(a))$$
where $(\tau_{n})^*(a\wa a)=(\tau_{m})_*(a\wc a)$ and $a\wc a -a\wa a=\partial P(H(a)\wedge TH(a))$.
This shows $\partial Sq_1(\alpha)=0$. Using the definition of secondary homotopy groups as track functors in
\cite{2hg1} and the
symmetric actions in \cite{2hg2} we see that $Sq_1(\alpha)$ coincides with the track definition of $\alpha\cup_1\alpha$.
\end{proof}

\begin{thm}\label{sq}
Let $n$ and $m$ be even positive integers and $\alpha,\beta\in\pi_nS^m$. Then
$$Sq_1(\alpha|\beta)=Sq_1(\alpha+\beta)-Sq_1(\alpha)-Sq_1(\beta)=\left(\frac{n-m}{2}+1\right)(\alpha\wedge\beta)(\S^{2n-3}\eta).$$
\end{thm}

This result is stated in \cite{K1p}, but a proof did not appear in the literature.

\begin{proof}[Proof of Theorem \ref{sq}]
We choose representatives $a,b\in\Pi_{n,0}S^m$ of $\alpha,\beta$ with $$H(a)=0=H(b).$$ Then we have $$a\wa b=a\wc
b$$ and we get 
\begin{eqnarray*}
(a+b)\wa(a+b)&=&a\wa (a+b)+b\wa(a+b)\\
&=&a\wa a+a\wa b+b\wa a+b\wa b\\
&=&x+u+y,\\
(a+b)\wc(a+b)&=& (a+b)\wc a+(a+b)\wc b\\
&=&a\wa a+b\wa a+a\wa b+b\wa b\\
&=&x+v+y.
\end{eqnarray*}
Here we set $x=a\wa a$, $y=b\wa b$ and
\begin{eqnarray*}
u&=&a\wa b+b\wa a,\\
v&=&b\wa a+a\wa b,
\end{eqnarray*}
so that 
\begin{eqnarray*}
u&=&v+\partial P(a\wa b|b\wa a)_H.
\end{eqnarray*}
Now the formula for $Sq_1(\alpha)$ yields:
\begin{eqnarray*}
Sq_1(\alpha)&=&-\grupo{a\wa a,\hat{\tau}_{n}}+(\hat{\tau}_{m})_*(a\wa a),\\
Sq_1(\beta)&=&-\grupo{b\wa b,\hat{\tau}_{n}}+(\hat{\tau}_{m})_*(b\wa b).
\end{eqnarray*}
Moreover for $\gamma=\alpha+\beta$ represented by $c=a+b$ we have 
\begin{eqnarray*}
Sq_1(\gamma)&=&-\grupo{c\wa c,\hat{\tau}_{n}}+(\hat{\tau}_{m})_*(c\wc c)-P(H(c)\wa TH(c)).
\end{eqnarray*}
The summands of $Sq_1(\gamma)$ satisfy the formulas:
\begin{eqnarray*}
\grupo{c\wa c,\hat{\tau}_{n}}&=&\grupo{x+u+y,\hat{\tau}_{n}}\\
&=&\grupo{x,\hat{\tau}_{n}}^{\tau_{n}^*(u+y)}+\grupo{u,\hat{\tau}_{n}}^{\tau_{n}^*(y)}
+\grupo{y,\hat{\tau}_{n}},\\
(\hat{\tau}_{m})_*(c\wc c)&=&(\hat{\tau}_{m})_*(x+v+y)\\
&=&(\hat{\tau}_{m})_*(x)^{(\tau_{m})_*(v+y)}\\
&&+(\hat{\tau}_{m})_*(v)^{(\tau_{m})_*(y)}+(\hat{\tau}_{m})_*(y).
\end{eqnarray*}
Here we have
\begin{eqnarray*}
(\tau_{n})^*(a\wa b)&=&(\tau_{m})_*(b\wc a)\\
&=&(\tau_{m})_*(b\wa a),\\
(\tau_{n})^*(u)&=&(\tau_{m})_*(v).
\end{eqnarray*}
Hence we get
\begin{eqnarray*}
Sq_1(\gamma)&=&-\grupo{y,\hat{\tau}_{n}}-\grupo{u,\hat{\tau}_{n}}^{\tau_{n}^*(y)}\\
&&+\left(-\grupo{x,\hat{\tau}_{n}}+(\hat{\tau}_{m})_*(x)\right)^{(\tau_{m})_*(v+y)}\\
&&+(\hat{\tau}_{m})_*(v)^{(\tau_{m})_*(y)}+(\hat{\tau}_{m})_*(y)-P(H(c)\wedge
TH(c)).
\end{eqnarray*}
Since the action on $\ker\partial$ is trivial and since the image of $P$ and $\ker\partial$ are both central we
thus get
\begin{eqnarray*}
Sq_1(\gamma)-Sq_1(\beta)&=&-\grupo{u,\hat{\tau}_{n}}^{\tau_{n}^*(y)}+Sq_1(\alpha)
+(\hat{\tau}_{m})_*(v)^{(\tau_{m})_*(y)}\\&&-P(H(c)\wedge TH(c)).
\end{eqnarray*}
Therefore we have
\begin{eqnarray*}
Sq_1(\alpha|\beta)&=&-\grupo{u,\hat{\tau}_{n}}^{\tau_{n}^*(y)}
+(\hat{\tau}_{m})_*(v)^{(\tau_{m})_*(y)}-P(H(c)\wedge TH(c))\\
&=&\left(-\grupo{u,\hat{\tau}_{n}}+(\hat{\tau}_{m})_*(v)\right)^{(\tau_{m})_*(y)}
+P(a\wa b|b\wa a)_H
\end{eqnarray*}
since $\partial(-\grupo{u,\hat{\tau}_{n}}+(\hat{\tau}_{m})_*(v))$ is a commutator, and hence in the
image of $\partial P$. 
Here we have
\begin{eqnarray*}
\grupo{u,\hat{\tau}_{n}}&=&\grupo{a\wa b,\hat{\tau}_{n}}+\grupo{b\wa a,\hat{\tau}_{n}}+P(-\tau_{n}^*(a\wa
b)+a\wa b|\tau_{n}^*(b\wa a))_H,\\
(\hat{\tau}_{m})_*(v)&=&(\hat{\tau}_{m})_*(b\wa a)+(\hat{\tau}_{m})_*(a\wa b)\\
&&+P(-(\tau_{m})_*(b\wa a)+b\wa a|(\tau_{m})_*(a\wa b))_H.
\end{eqnarray*}
Hence we obtain 
\begin{eqnarray*}
Sq_1(a|b)&=&-\grupo{b\wa a,\hat{\tau}_{n}}-\grupo{a\wa b,\hat{\tau}_{n}}+(\hat{\tau}_{m})_*(b\wa
a)+(\hat{\tau}_{m})_*(a\wa b)+\text{(c)},\\
\text{(c)}&=&P(b\wa a-a\wa b|(\hat{\tau}_{m})_*(a\wa b))+P(a\wa b|b\wa a)_H.
\end{eqnarray*}
Now we consider the following formulas with $\partial\text{(a)}=0=\partial\text{(b)}$.
\begin{eqnarray*}
\text{(a)}&=&\grupo{a\wa b,\hat{\tau}_{n}\hat{\tau}_{n}}\\
&=&\grupo{\tau_{n}^*(a\wa b),\hat{\tau}_{n}}+\grupo{a\wa b,\hat{\tau}_{n}},\\
\text{(b)}&=&(\hat{\tau}_{m}\hat{\tau}_{m})_*(a\wa b)\\
&=&(\hat{\tau}_{m})_*((\tau_{m})_*(a\wa b))+(\hat{\tau}_{m})_*(a\wa b).
\end{eqnarray*}
For (b) we use that $\Pi_{2n,*}$ is a track functor and 
\begin{eqnarray*}
\hat{\tau}_{m}\hat{\tau}_{m}&=&\hat{\tau}_{m}\vc(\tau_{m}\hat{\tau}_{m})\\
&=&\hat{\tau}_{m}\vc(\hat{\tau}_{m}\tau_{m}).
\end{eqnarray*}
Now we have the following equations.
\begin{eqnarray*}
\grupo{a\wa b,\hat{\tau}_{n}}+\grupo{b\wa a,\hat{\tau}_{n}}&=&
\text{(a)}-\grupo{\tau_{n}^*(a\wa b),\hat{\tau}_{n}}+\grupo{b\wa a,\hat{\tau}_{n}}\\
&=&\text{(a)}+\grupo{-(\tau_{m})_*(b\wa a)+b\wa a,\hat{\tau}_{n}}-\text{(d)}\\
&=&\text{(a)}+\grupo{\partial(\hat{\tau}_{m})_*(b\wa a),\hat{\tau}_{n}}-\text{(d)}\\
&=&\text{(a)}-\tau_{n}^*(\hat{\tau}_{m})_*(b\wa a)+(\hat{\tau}_{m})_*(b\wa a)-\text{(d)}\\
&=&\text{(a)}-(\hat{\tau}_{m})_*\tau_{n}^*(b\wa a)+(\hat{\tau}_{m})_*(b\wa a)-\text{(d)}\\
&=&\text{(a)}-(\hat{\tau}_{m})_*(\tau_{m})_*(a\wa b)+(\hat{\tau}_{m})_*(b\wa
a)-\text{(d)}.
\end{eqnarray*}
Here (d) is given by the following formula.
\begin{eqnarray*}
\text{(d)}&=&-P(\tau_{n}^*(\tau_{m})_*(b\wa a)-(\tau_{m})_*(b\wa
a))|\tau_{n}^*(\tau_{m})_*(b\wa a))_H\\
&&+P(\tau_{n}^*(\tau_{m})_*(b\wa a)-(\tau_{m})_*(b\wa
a))|\tau_{n}^*(b\wa a))_H\\
&=&-P(a\wa b-(\tau_{m})_*(b\wa a)|a\wa b)_H\\
&&+P(a\wa b-(\tau_{m})_*(b\wa a)|\tau^*_{n}(b\wa a))_H\\
&=&-P(a\wa b|a\wa b)_H+P((\tau_{m})_*(b\wa a)|a\wa b)_H\\
&&+P(a\wa b|(\tau_{m})_*(b\wa a))_H-P((\tau_{m})_*(b\wa a)|(\tau_{m})_*(a\wa b))_H.
\end{eqnarray*}
On the other hand we get
\begin{eqnarray*}
(\hat{\tau}_{m})_*(b\wa a)+(\hat{\tau}_{m})_*(a\wa b)&=&(\hat{\tau}_{m})_*(b\wa a)\\
&&-(\hat{\tau}_{m})_*((\tau_{m})_*(a\wa b))+\text{(b)}.
\end{eqnarray*}
Now we get
\begin{eqnarray*}
Sq_1(\alpha|\beta)&=&-\text{(a)}-\text{(d)}-(\hat{\tau}_{m})_*(b\wa a)
+(\hat{\tau}_{m})_*((\tau_{m})_*(a\wa b))\\
&&+(\hat{\tau}_{m})_*(b\wa a)-(\hat{\tau}_{m})_*((\tau_{m})_*(a\wa
b))+\text{(b)}+\text{(c)}\\
&=&-\text{(a)}+\text{(d)}+\text{(e)}+\text{(b)}+\text{(c)}
\end{eqnarray*}
where (e) is the commutator:
\begin{eqnarray*}
\text{(e)}&=&P((\hat{\tau}_{m})_*(b\wa a)|-(\hat{\tau}_{m})_*((\tau_{m})_*(a\wa
b))_{H\partial}\\
&=&P(-(\tau_{m})_*(b\wa a)+b\wa a|(\tau_{m}\tau_{m})_*(a\wa b)\\
&&-(\tau_{m})_*(a\wa
b))_H\\
&=&-P((\tau_{m})_*(b\wa a)|a\wa b)_H+P(b\wa a|a\wa b)_H\\
&&+P((\tau_{m})_*(b\wa a)|(\tau_{m})_*(a\wa b))_H\\
&&-P(b\wa a|(\tau_{m})_*(a\wa b))_H.
\end{eqnarray*}

\begin{eqnarray*}
\text{(d)}+\text{(e)}+\text{(c)}&=&-P(a\wa b|a\wa b)_H+P((\tau_{m})_*(b\wa a)|a\wa b)_H\\
&&+P(a\wa b|(\tau_{m})_*(a\wa b))_H\\
&&-P((\tau_{m})_*(b\wa a)|(\tau_{m})_*(a\wa
b))_H\\
&&-P((\tau_{m})_*(b\wa a)|a\wa b)_H+P(b\wa a|a\wa b)_H\\
&&+P((\tau_{m})_*(b\wa a)|(\tau_{m})_*(a\wa b))_H\\
&&-P(b\wa a|(\tau_{m})_*(a\wa
b))_H\\
&&+P(b\wa a|(\tau_{m})_*(a\wa b))_H-P(a\wa b|(\tau_{m})_*(a\wa b))_H\\
&&+P(a\wa b|b\wa a)_H\\
&=&-P(a\wa b|a\wa b)_H.
\end{eqnarray*}
Hence we get
\begin{eqnarray*}
Sq_1(\alpha|\beta)=-\text{(a)}+\text{(b)}-P(a\wa b|a\wa b)_H,
\end{eqnarray*}
and this implies the result by use of Proposition \ref{qlo}.
\end{proof}

\section{Toda brackets}\label{t11}

For a pointed space $X$ we use the suspension $\S X=S^1\wedge X$ and the \emph{$E$-suspension} $E X=X\wedge S^1$.
Here $\S$ and $E$ are isomorphic endofunctors of $\C{Top}^*$. The $E$-suspension is for example used by Toda in
his book \cite{toda}.

\begin{defn}\label{tb}
Let $n\geq k\geq 0$ and consider morphisms in $\C{Top}^*/\simeq$
$$Z\st{\alpha}\l E^kY\st{E^k\beta}\l E^k X,\;\; Y\st{\beta}\l X\st{\gamma}\l S^{n-k}$$
with $\alpha(E^k\beta)=0$ and $\beta\gamma=0$. Then the \emph{Toda bracket}
$$\set{\alpha,E^k\beta,E^k\gamma}_k\subset\pi_{n+1}Z$$
is the subset of all elements in $\pi_{n+1}Z$ obtained by pasting tracks as in the diagram 
$$\xymatrix{Z&E^kY\ar[l]_a&E^kX\ar[l]_{E^kb}_<(.905){\;}="d"\ar@/_40pt/[ll]_0^{\;}="c"&
S^n\ar[l]_{E^kc}^<(1.05){\;}="a"\ar@/^40pt/[ll]^0_{\;}="b"\ar@{=>}"a";"b"^{E^kB}\ar@{=>}"c";"d"^A}$$
where $a$, $b$, $c$ represent $\alpha$, $\beta$, $\gamma$ and $B\colon bc\rr 0$ and $A\colon 0\rr a(E^kb)$.
\end{defn}

Let $\imath_k\in\pi_k(S^k)$ be the element represented by the identity of $S^k$, $k\geq 0$. Moreover, let 
$$\imath_{k,0}\in\L^kS^k\subset\Pi_{k,0}(S^k)=\grupo{\L^kS^k}_\nill$$
be given by the identity of $S^k$. 
This element yields a quadratic pair module morphism
$$\imath_{k,0}\colon\overline{\Z}_\nill\To\Pi_{k,*}S^k.$$
We define the morphisms $\bar{E}^k$ in $\C{qpm}$, $n\geq k$,
\begin{equation}\label{ebar}
\bar{E}^k\colon\Pi_{n-k,*}X\cong\Pi_{n-k,*}X\odot\overline{\Z}_\nill\st{1\otimes
\imath_{k,0}}\r\Pi_{n-k,*}X\odot\Pi_{k,*}S^k\st{\wedge}\r\Pi_{n,*}(E^kX).
\end{equation}

Let $\alpha$, $\beta$, $\gamma$ be given as in Definition \ref{tb} with $\alpha(E^k\beta)=0$ and
$\beta\gamma=0$. We choose maps $a$, $b$ representing $\alpha$, $\beta$ and we choose a track $A\colon 0\rr
a(E^kb)$ as in Definition \ref{tb}. Moreover, let $$\bar{c}\in\Pi_{n-k,0}X$$ be an element representing
$\gamma\in\pi_{n-k}X$ with $n-k\geq 2$ and let 
\begin{equation}\label{:P}
\bar{B}\in\Pi_{n-k,1}Y
\end{equation} 
be an element with
$\partial\bar{B}=b_*(\bar{c})\in\Pi_{n-k,0}Y$. Then $\bar{E}^k\bar{B}\in\Pi_{n,1}E^kY$ satisfies
$$\partial\bar{E}^k\bar{B}=\partial(\bar{B}\wa\imath_{k,0})=(\partial\bar{B})\wa\imath_{k,0}=b_*(\bar{c})\wa\imath_{k,0}=
(E^kb)_*(\bar{c}\wa\imath_{k,0}).$$
Moreover, the track $A$ induces a track in $\C{qpm}$ which is given by a map
$$A_*\colon\Pi_{n,0}E^kX\To\Pi_{n,1}Z$$ with $\partial A_*(x)=(a(E^kb))_*(x)$. Therefore
the element
$$t=a_*(E^k\bar{B})-A_*(\bar{c}\wedge\imath_{k,0})\in\Pi_{n,1}Z$$
satisfies $\partial(t)=0$ and hence $t$ is an element in $h_1\Pi_{n,*}Z$. Recall from (\ref{h0}) that $h_1\Pi_{n,*}Z$
is naturally isomorphic to $\pi_{n+1}Z$ for $n\geq 3$ and to $\pi_3Z/[\pi_2Z,\pi_2Z]$ for $n=2$, where $[-,-]$ is
the Whitehead product.

\begin{lem}\label{:S}
For $n\geq 3$, $a_*\bar{E}^k(\bar{B})-A_*(\bar{c}\wedge\imath_{k,0})\in\set{\alpha,E^k\beta,E^k\gamma}_k$.
Moreover, all elements in $\iota\set{\alpha,E^k\beta,E^k\gamma}_k$ can be obtained in this way. The same equality
holds for $n=2$ $\mod [\pi_2Z,\pi_2Z]$, the image of the Whitehead product.
\end{lem}

\begin{prop}\label{topro}
Let $Y$ be a pointed space and let $r\in\Z$ and $\beta\in\pi_{n-1}Y$ with $r\beta=0$. In the group of homotopy
classes $[EY,EY]$ let $r1_{EY}$ be the $r$-fold sum of the identity $1_{EY}$. Then the Toda bracket
$\set{r1_{EY},E\beta,r\imath_n}_1\subset\pi_{n+1}EY$ is defined and for $n\geq 3$ 
$$\set{r1_{EY},E\beta,r\imath_n}_1\ni\left\{\begin{array}{ll}
0, & \text{if $r$ is odd},\\
&\\
\frac{r}{2}(E\beta)(\S^{n-2}\eta) & \text{if $r$ is even}.
\end{array}\right.$$
Here $\eta\colon S^3\r S^2$ is the Hopf map. For $n=2$ the same formula holds in the quotient  $\pi_3EY/[\pi_2EY,\pi_2EY]$ where
$[-,-]$ is the Whitehead product.
\end{prop}

In \cite{toda} 3.7 Toda proves this result in case $Y$ is a sphere. Toda's proof uses  different methods relying
on the assumption that $Y$ is a sphere.

\begin{proof}[Proof of \ref{topro}]
Let $(\cdot)^r\colon S^1\r S^1$ be the degree $r$ map $z\mapsto z^r$ so that $r1_{EY}$ is represented by
$a=Y\wedge(\cdot)^r$.  Let $b$ be a map representing $\beta$. We choose $\bar{c}\in\Pi_{n-1,0}S^{n-1}$
representing $r\imath_n$ by $\bar{c}=r\imath_{n-1,0}$ and we choose $$\bar{B}\in\Pi_{n-1,1}Y$$ with
$\partial\bar{B}=b_*(\bar{c})=r\, b_*(\imath_{n-1,0})$.
Then we get as in (\ref{:P})
$$\partial\bar{E}(\bar{B})=(Eb)_*(\bar{c}\wa\imath_{1,0})$$
where $\bar{c}\wa\imath_{1,0}=r(\imath_{n-1,0}\wa\imath_{1,0})=r\imath_{n,0}$.
Now we choose $A\colon0\rr a(Eb)$ in such a way that the induced track $A_*$ satisfies
$$A_*(\imath_{n,0})=\bar{E}(\bar{B}).$$ In fact, the boundary of $A_*(\imath_{n,0})$ is
$$\partial A_*(\imath_{n,0})=(a(Eb))_*(\imath_{n,0})=\partial\bar{E}(\bar{B})$$
where $a(Eb)=(Y\wedge(\cdot)^r)(b\wedge S^1)=b\wedge(\cdot)^r=(Eb)(\S^{n-1}(\cdot)^r)$ and
$(\S^{n-1}(\cdot)^r)_*(\imath_{n,0})=r^*\imath_{n,0}=r\imath_{n,0}$. since $H(i_{n,0})=0$.

Now we can compute the element in Lemma \ref{:S}
$$t=a_*\bar{E}(\bar{B})-A_*(\bar{c}\wedge\imath_{1,0})\in\set{r1_{EY},E\beta,r\imath_n}_1.$$
Here $A$ is a track $f_0\rr g_0$ with $f_0=0$ so that
\begin{eqnarray*}
A_*(\bar{c}\wedge\imath_{1,0})&=&A_*(r\imath_{n,0})\\
&=&r\, A_*(\imath_{n,0})\\
&=&r\bar{E}\bar{B}.
\end{eqnarray*}
On the other hand we get 
\begin{eqnarray*}
a_*(\bar{E}\bar{B})&=&(Y\wedge(\cdot)^r)_*(\bar{B}\wa\imath_{1,0})\\
&=&\bar{B}\wa(((\cdot)^r)_*(\imath_{1,0}))\\
&=&\bar{B}\wa(r^*\imath_{1,0})\\
&=&\bar{B}\wa(r\imath_{1,0})\\
&=&r(\bar{B}\wa\imath_{1,0})+t\\
&=&r\bar{E}(\bar{B})+t
\end{eqnarray*}
where
\begin{eqnarray*}
t&=&\binom{r}{2}P(H\partial(\bar{B})\wa(\imath_{1,0}|\imath_{1,0})_H)\\
&=&\binom{r}{2}P(H(r\, b_*(\imath_{n-1,0}))\wa(\imath_{1,0}|\imath_{1,0})_H)\\
&=&{\binom{r}{2}}^2P((b_*(\imath_{n-1,0})|b_*(\imath_{n-1,0}))_H\wa(\imath_{1,0}|\imath_{1,0})_H)\\
&=&{\binom{r}{2}}^2P(b_*(\imath_{n-1,0})\wa\imath_{1,0}|b_*(\imath_{n-1,0})\wa\imath_{1,0})_H.
\end{eqnarray*}
Therefore $t$ represents ${\binom{r}{2}}^2(E\beta)(\S^{n-2}\eta)$. If $r$ is odd we see that $t=0$ since $\beta$
has odd order and $\S\eta$ has even order.
\end{proof}

\part{The construction of the smash product for secondary homotopy groups}\label{III}

In this part we define the smash product operator for secondary homotopy groups and we prove the results
described in Part \ref{II}. A crucial step for this definition will be the construction of canonical tracks 
$$\S(f\sc g)\st{\aleph^{\sc}}\rr f\wedge g\st{\aleph^\sa}\Leftarrow \S(f\sa g)$$
termed the \emph{exterior tracks}, connecting the exterior cup-products and the smash product of two maps. 
Then we use $\aleph^\sa$ and $\aleph^{\sc}$ for the definition of the smash
product operator on the $(1)$-level. Some of the algebraic properties of the smash product are then derived from formulas
concerning the Hopf invariant of the track $(\aleph^{\sc})^\vi\vc\aleph^{\sa}$ 
and of some other tracks between suspensions built out of the exterior tracks.

\section{Exterior cup-products for higher suspensions and tracks}

We begin this section by stating the basic properties of the exterior cup-product operations.

\begin{lem}\label{proep}
We have the following formulas for
suspensions
\begin{enumerate}
\item $f\sa (\S g')=f\wedge g'=f\sc(\S g')$,
\item $(\S f')\sa g=(1\;2)(f'\wedge g)(1\;2)=(\S f')\sc g$,
\end{enumerate}
coproducts
\begin{enumerate}\setcounter{enumi}{2}
\item $(f_1,f_2)\sa g=(f_1\sa g, f_2\sa g)$,
\item $(f_1,f_2)\sc g=(f_1\sc g, f_2\sc g)$,
\item $f\sa (g_1,g_2)=(f\sa g_1,f\sa g_2)$,
\item $f\sc (g_1,g_2)=(f\sc g_1,f\sc g_2)$,
\end{enumerate}
and compositions
\begin{enumerate}\setcounter{enumi}{6}
\item $(f_1f_2)\sa (g_1(\S g_2'))=(f_1\sa g_1)(f_2\sa (\S g_2'))$,
\item $((\S f_1')f_2)\sa (g_1g_2)=((\S f_1')\sa g_1)(f_2\sa g_2)$,
\item $(f_1(\S f_2'))\sc (g_1g_2)=(f_1\sc g_1)((\S f_2')\sc g_2)$,
\item $(f_1f_2)\sc ((\S g_1')g_2)=(f_1\sc (\S g_1'))(f_2\sc g_2)$.
\end{enumerate}
The exterior cup-products are associative
\begin{enumerate}\setcounter{enumi}{10}
\item $f\sa (g\sa h)=(f\sa g)\sa h$,
\item $f\sc (g\sc h)=(f\sc g)\sc h$.
\end{enumerate}
\end{lem}

The proof of this lemma is straightforward.

In order to define the exterior cup-products $$f\sa g,f\sc g\colon\S^{n+m-1}A\wedge Y\r \S^{n+m-1}B\wedge Y$$ 
of maps between higher suspensions $f\colon\S^nA\r\S^nB$, $g\colon \S^mX\r\S^mY$ we take the first spherical
coordinates to the end of the smash product
$$\bar{f}\colon S^1\wedge A\wedge S^{n-1}\cong S^{n-1}\wedge S^1\wedge A\st{f}\To S^{n-1}\wedge S^1\wedge B\cong S^1\wedge
B\wedge S^{n-1},$$
$$\bar{g}\colon S^1\wedge X\wedge S^{m-1}\cong S^{m-1}\wedge S^1\wedge X\st{g}\To S^{m-1}\wedge S^1\wedge Y\cong S^1\wedge
Y\wedge S^{m-1},$$
then we perform the usual exterior cup product on these maps, and we recollect the permuted spherical coordinates
at the beginning of the smash product in an ordered way,
\begin{eqnarray*}f\sa g\colon S^{n-1}\wedge S^{m-1}\wedge S^1\wedge A\wedge X&\cong&
S^1\wedge A\wedge S^{n-1}\wedge X \wedge S^{m-1}\\&\st{\bar{f}\sa \bar{g}}\To&
S^1\wedge B\wedge S^{n-1}\wedge Y \wedge S^{m-1}\\&\cong&
S^{n-1}\wedge S^{m-1}\wedge S^1\wedge B\wedge Y,
\end{eqnarray*}
and the same for $\sc$. These exterior cup-products generalize the classical ones in the following sense. If
$f=\S^{n-1}f'$ and $g=\S^{m-1}g'$ for $f\colon\S A\r \S B$ and $g\colon \S X\r \S Y$ then
\begin{eqnarray*}
(\S^{n-1}f')\sa(\S^{m-1}g')&=&\S^{n+m-1}(f'\sa g'),\\
(\S^{n-1}f')\sc(\S^{m-1}g')&=&\S^{n+m-1}(f'\sc g').
\end{eqnarray*}

The properties of the classical exterior cup-products in Lemma \ref{proep} can be accordingly restated for the
exterior cup-product of maps between higher suspensions.

Let $F\colon f\rr
g$, $G\colon h\rr k$ be now tracks between maps $f,g\colon\S^nA\r\S^nB$, $h,k\colon \S^mX\r\S^mY$. 
The exterior products of a track with a map
$$\begin{array}{c}
F\sa h\colon f\sa h\rr g\sa h,\\
F\sc h\colon f\sc h\rr\ g\sc h,\\
f\sa G\colon f\sa h\rr f\sa k,\\
f\sc G\colon f\sc h\rr f\sc k,
\end{array}$$
are defined by exchanging the interval $I_+$ with the spherical coordinates and using the exterior cup-products of
maps between higher suspensions as defined above. For example the track $F\sa h$ is represented by the homotopy
$$
I_+\wedge S^{n+m-1}\wedge A\wedge X\cong
S^{n+m-1}\wedge I_+\wedge A\wedge X\st{\tilde{F}\sa h}\To
S^{n+m-1}\wedge B\wedge Y,
$$
where $$\tilde{F}\colon S^{n-1}\wedge I_+\wedge A\cong I_+\wedge S^{n-1}\wedge A\st{F}\To S^{n-1}\wedge B$$
is defined from a homotopy $F$ representing the corresponding track.
Now one can define the exterior products of two tracks as the vertical composition
$$
\begin{array}{c}
F\sa G = (g\sa G)\vc(F\sa h)=(F\sa k)\vc(f\sa G),\\
F\sc G = (g\sc G)\vc(F\sc h)=(F\sc k)\vc(f\sc G).
\end{array}
$$

One can also derive from Lemma \ref{proep} analogous properties for the exterior cup-products of tracks.

\section{The exterior tracks}\label{tet}

For any two maps $f\colon\S A\r \S B$ and $g\colon\S X\r\S Y$ the
suspended exterior cup-products $\Sigma (f\sa  g)$ and $\S (f\sc  g)$
are naturally homotopic to the composite $$S^1\wedge S^1\wedge A
\wedge X\st{(2\;3)}\cong S^1\wedge A\wedge S^1\wedge X\st{f\wedge
g}\To S^1\wedge B\wedge S^1\wedge Y\st{(2\;3)}\cong S^1\wedge
S^1\wedge B \wedge Y.$$ 
In order
to construct homotopies we only need to choose a track from the
transposition map $(1\;2)\colon S^1\wedge S^1\r S^1\wedge S^1$ to
$\nu\wedge S^1$, where $\nu\colon S^1\r S^1$ is the co-H-inversion
defined by $\nu(z)=z^{-1}$. Here we use the topological group structure of
$S^1$. The set of all tracks $(1\;2)\rr \nu\wedge S^1$ and
$1_{S^2}\rr 1_{S^2}$ is a group under horizontal composition.
This group is an extension of $\Z/2$ by $\Z$ with the non-trivial action of $\Z/2$, compare
\cite{2hg2} 6.12. Up to isomorphism there is
only one extension of this kind, the trivial extension, given by the infinite dyhedral group $\Z/2*\Z/2$, hence
this group of tracks is  generated by two order $2$ tracks $(1\;2)\rr
\nu\wedge{S^1}$. One of these two generating tracks $\aleph\colon(1\;2)\rr
\nu\wedge{S^1}$ can be constructed as follows. Since $\nu\wedge S^1$
is a homotopy equivalence it is enough to indicate which track is
$\aleph^\vi(\nu\wedge S^1)\colon1_{S^2}\rr(1\;2)(\nu\wedge S^1)$.
The $2$-sphere $S^2=S^1\wedge S^1$ is a quotient of the square
$[-1,1]^2$
by the map $[-1,1]^2\r S^1\wedge S^1 \colon (x,y)\mapsto(\exp \pi
i(1+x),\exp\pi i(1+y))$. The map $(1\;2)(\nu\wedge S^1)$ is induced
by the $90^\circ$ twist (counterclockwise) in the square, so we
obtain $\aleph^\vi(\nu\wedge S^1)$ by using the homeomorphism from
the square to the radius $\sqrt{2}$ circle projecting
from the origin
$$
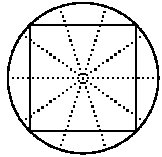
\xy/r1cm/:{\xypolygon4{}}{\xypolygon10{~<{.}~>{}}}\ellipse(1){-}\endxy$$
and twisting continuously the circle $90^\circ$ counterclockwise.

\begin{defn}\label{tracklos}
The \emph{exterior tracks} $$\aleph^{\sa}_{f,g} \colon\Sigma (f\sa
g)\rr(2\;3)(f\wedge g)(2\;3),$$
$${\aleph}^{\sc}_{f,g}\colon \S (f\sc g)\rr(2\;3)(f\wedge
g)(2\;3),$$ are given by the following diagram in $\C{Top}^*$ where
the $2$-cells without a track arrow $\rr$ are strictly commutative
\begin{footnotesize}
\begin{equation*}
\rotatebox{90}{\xymatrix@C=15pt{S^1\wedge S^1\wedge A\wedge
X\ar[rrrrrr]^{\S(f\sa g)}\ar@{=}[dd]&&&&&&
S^1\wedge S^1\wedge B\wedge Y\ar@{=}[dd]\\\\
S^1\wedge S^1\wedge A\wedge X
\ar@/^30pt/[r]^{1                                                     
}_{\;}="a"\ar@/_25pt/[r]^{\;}="b"_{((\nu\wedge S^1)(1\;2))\wedge
A\wedge X}& S^1\wedge S^1\wedge A\wedge X\ar[r]^{S^1\wedge f\wedge
X}& S^1\wedge S^1\wedge B\wedge X
\ar@/^30pt/[r]^{1                                                     
}_{\;}="c"\ar@/_25pt/[r]^{\;}="d"_{((1\;2)(\nu\wedge S^1))\wedge
B\wedge X}& S^1\wedge S^1\wedge B\wedge X\ar[r]^{(2\;3)}& S^1\wedge
B\wedge S^1\wedge X\ar[r]^{S^1\wedge B\wedge g}& S^1\wedge B\wedge
S^1\wedge Y\ar[r]^{(2\;3)}& S^1\wedge S^1\wedge B\wedge Y
\ar@{=>}"a";"b"|<(.29){((\nu\wedge S^1)\aleph^\vi)\wedge A\wedge X}
\ar@{=>}"c";"d"|<(.29){(\aleph^\vi(\nu\wedge S^1))\wedge B\wedge X}\\\\
S^1\wedge S^1\wedge A\wedge X\ar[r]^{(2\;3)}\ar@{=}[uu]\ar@{=}[dd]&
S^1\wedge A\wedge S^1\wedge X\ar[rrrr]^{f\wedge g}&&&&
S^1\wedge B\wedge S^1\wedge Y\ar[r]^{(2\;3)}&S^1\wedge S^1\wedge B\wedge Y\ar@{=}[uu]\ar@{=}[dd]\\\\
S^1\wedge S^1\wedge A\wedge X\ar[r]^{(2\;3)}& S^1\wedge A\wedge
S^1\wedge X\ar[r]^{S^1\wedge A\wedge g}& S^1\wedge A\wedge S^1\wedge
Y\ar[r]^{(2\;3)}& S^1\wedge S^1\wedge A\wedge Y
\ar@/_30pt/[r]_{1                                                     
}^{\;}="e"\ar@/^25pt/[r]_{\;}="f"^{((\nu\wedge S^1)(1\;2))\wedge
A\wedge Y}& S^1\wedge S^1\wedge A\wedge Y\ar[r]^{S^1\wedge f\wedge
Y}& S^1\wedge S^1\wedge B\wedge Y
\ar@/_30pt/[r]_{1                                                     
}^{\;}="g"\ar@/^25pt/[r]_{\;}="h"^{((1\;2)(\nu\wedge S^1))\wedge
B\wedge Y}& S^1\wedge S^1\wedge B\wedge Y
\ar@{=>}"e";"f"|<(.30){((\nu\wedge S^1)\aleph^\vi)\wedge A\wedge Y}
\ar@{=>}"g";"h"|<(.30){(\aleph^\vi(\nu\wedge S^1))\wedge B\wedge Y}\\\\
S^1\wedge S^1\wedge A\wedge X\ar[rrrrrr]^{\S(f\sc
g)}\ar@{=}[uu]&&&&&& S^1\wedge S^1\wedge B\wedge Y\ar@{=}[uu]}}
\end{equation*}
\end{footnotesize}
i. e. the tracks $\aleph^{\sa}_{f,g}$ and ${\aleph}^{\sc}_{f,g}$ are the following composite tracks.
\begin{footnotesize}
\begin{eqnarray*}
\aleph^{\sa}_{f,g}&=&(2\;3)(S^1\wedge B\wedge g)(2\;3)
((\aleph^\vi(\nu\wedge S^1))\wedge B\wedge X)(S^1\wedge f\wedge
X)(((\nu\wedge S^1)\aleph^\vi)\wedge A\wedge X)\\
&=&(2\;3)(S^1\wedge B\wedge g)(2\;3) (\aleph^\vi\wedge B\wedge
X)(S^1\wedge f\wedge
X)(\aleph^\vi\wedge A\wedge X),\\
{\aleph}^{\sc}_{f,g}&=&((\aleph^\vi(\nu\wedge S^1))\wedge B\wedge
Y)(S^1\wedge f\wedge Y) (((\nu\wedge S^1)\aleph^\vi)\wedge A\wedge
Y) (2\;3)(S^1\wedge A\wedge g)(2\;3)
\\&=&(\aleph^\vi\wedge B\wedge Y)(S^1\wedge f\wedge Y)
(\aleph^\vi\wedge A\wedge Y) (2\;3)(S^1\wedge A\wedge g)(2\;3).
\end{eqnarray*}
\end{footnotesize}
\end{defn}

In the next proposition we show elementary properties of the exterior tracks that are
relevant for the definition of the smash product operation on
secondary homotopy groups. They are analogous to the properties of exterior cup-products in Lemma \ref{proep}.

\begin{lem}\label{exteprop}
The exterior tracks satisfy the following formulas for suspensions
\begin{enumerate}
\item $\aleph^\sa_{\S f',g}=0^\vc_{\S((1\;2)(f'\wedge g)(1\;2))}=\aleph^{\sc}_{\S f',g}$,
\item $\aleph^{\sa}_{f,\S g'}=(\aleph^\vi\wedge B\wedge Y)
(S^1\wedge f\wedge g')(\aleph^\vi\wedge A\wedge
X)=\aleph^{\sc}_{f,\S g'}$,
\end{enumerate}
coproducts
\begin{enumerate}\setcounter{enumi}{2}
\item
$\aleph^\sa_{(f_1,f_2),g}=(\aleph^\sa_{f_1,g},\aleph^\sa_{f_2,g})$,
\item
$\aleph^{\sc}_{(f_1,f_2),g}=(\aleph^{\sc}_{f_1,g},\aleph^{\sc}_{f_2,g})$,
\item
$\aleph^{\sa}_{f,(g_1,g_2)}=(\aleph^{\sa}_{f,g_1},\aleph^{\sa}_{f,g_2})$,
\item
$\aleph^{\sc}_{f,(g_1,g_2)}=(\aleph^{\sc}_{f,g_1},\aleph^{\sc}_{f,g_2})$,
\end{enumerate}
and composition of maps
\begin{enumerate}\setcounter{enumi}{6}
\item
$\aleph^\sa_{f_1f_2,g_1(\S
g_2')}=\aleph^\sa_{f_1,g_1}\aleph^\sa_{f_2,\S g_2'}$,
\item $\aleph^\sa_{(\S f_1')f_2,g_1g_2}=\aleph^\sa_{\S f_1',g_1}\aleph^\sa_{f_2,g_2}$,
\item
$\aleph^{\sc}_{f_1(\S
f_2'),g_1g_2}=\aleph^{\sc}_{f_1,g_1}\aleph^{\sc}_{\S f_2',g_2}$,
\item $\aleph^{\sc}_{f_1f_2,(\S g_1')g_2}=\aleph^{\sc}_{f_1,\S g_1'}\aleph^{\sc}_{f_2,g_2}$.
\end{enumerate}
They satisfy the following associativity rules. 
\begin{enumerate}\setcounter{enumi}{10}
\item $(\aleph^{\sa}_{f,g}\wedge h)(\S\aleph^{\sa}_{f\sa g,h})=(f\wedge\aleph^{\sa}_{g,h})(\S\aleph^{\sa}_{f,g\sa
h})\colon \S^2(f\sa g\sa h)\rr f\wedge g\wedge h$,
\item $(\aleph^{\sc}_{f,g}\wedge h)(\S\aleph^{\sc}_{f\sc g,h})=(f\wedge\aleph^{\sc}_{g,h})(\S\aleph^{\sc}_{f,g\sc
h})\colon \S^2(f\sc g\sc h)\rr f\wedge g\wedge h$.
\end{enumerate}
\end{lem}

These properties follow easily from the definition of
exterior tracks above and from the fact that
$\aleph\aleph=0^\vc_{1_{S^2}}$ is the trivial track.
In the right hand side of the equalities (11) and(12) there are some permutations involved that we have omitted.

For suspended maps $\S^{n-1}f\colon\S^nA\r\S^nB$, $\S^{m-1}g\colon\S^mX\r\S^mY$ we define the tracks
$\aleph^{\sa}_{n,f,m,g}$, $\aleph^{\sc}_{n,f,m,g}$ from $\S^{n+m-1}(f\sa g)$,
$\S^{n+m-1}(f\sc g)$, respectively, to
\begin{equation}\label{perara}
S^n\wedge S^m\wedge A\wedge X\cong
S^n\wedge A\wedge S^m\wedge X\st{\S^{n-1} f\wedge \S^{m-1}g}\To S^n\wedge B\wedge S^m\wedge Y
\cong S^n\wedge S^m\wedge B\wedge Y
\end{equation}
as
\begin{equation*}
\xymatrix@R=30pt{{\begin{array}{c}S^{n-1}\wedge S^1\wedge S^{m-1}\wedge S^1\wedge A\wedge X\\\cong^{(2\;3)}\\
S^{n-1}\wedge S^{m-1}\wedge S^1\wedge S^1\wedge A\wedge X\end{array}}
\ar@/^40pt/[d]_{\;}="a"^{\S^{n+m-2}((2\;3)(f\wedge g)(2\; 3))}\ar@/_40pt/[d]^{\;}="b"_{\S^{n+m-1}(f\sa g)}\\
{\begin{array}{c}S^{n-1}\wedge S^{m-1}\wedge S^1\wedge S^1\wedge B\wedge Y\\
\cong^{(2\;3)}\\
S^{n-1}\wedge S^1\wedge S^{m-1}\wedge S^1\wedge B\wedge Y
\end{array}}\ar@{=>}"b";"a"^{\S^{n+m-2}\aleph^{\sa}_{f,g}}}
\end{equation*}
and similarly for $\sc$.
Notice that the last spherical coordinate in these smash products is always the same one. This is relevant in
connection with Lemma \ref{nilson} below.
The tracks $\aleph^{\sa}_{n,f,m,g}$, $\aleph^{\sc}_{n,f,m,g}$ satisfy properties analogous to Lemma
\ref{exteprop} that we do not restate. They also satisfy the following further properties.

\begin{lem}\label{exteprop2}
Given tracks $F\colon \S^{n-1} f_1\rr \S^{n-1} f_2$ and
$G\colon \S^{m-1} g_1\rr \S^{m-1} g_2$ between maps $\S^{n-1}f_i\colon\S^nA\r\S^nB$ and 
$\S^{m-1}g_i\colon\S^mA\r\S^mB$, $i=1,2$, the following equalities are satisfied
\begin{enumerate}
\item $\begin{array}{l}
(F\wedge G)\vc\aleph^{\sa}_{n,f_1,m,g_1}=
\aleph^{\sa}_{n,f_2,m,g_2}\vc(\S(F\sa G))\colon\\
\S^{n+m-1}(f_1\sa g_1)\rr\S^{n+m-2}(f_2\wedge g_2),
\end{array}$

\item $\begin{array}{l}(F\wedge G)\vc\aleph^{\sc}_{n,f_1,m,g_1}=
\aleph^{\sc}_{n,f_2,m,g_2}\vc(\S(F\sc G))\colon\\
\S^{n+m-1}(f_1\sc g_1)\rr\S^{n+m-2}(f_2\wedge g_2).
\end{array}$
\end{enumerate}
\end{lem}

Here the track $F\wedge G$ in (1) and (2) needs to be altered by permutations according to (\ref{perara}). Moreover
$\S(F\sa G)$ and $\S(F\sc G)$ need also to be altered by permutations as follows.
The track $\S(F\sa G)$ should actually be
\begin{equation*}
\xymatrix@R=30pt{{\begin{array}{c}S^{n-1}\wedge S^1\wedge S^{m-1}\wedge S^1\wedge A\wedge X\\\cong^{(1\;2)}\\
S^1\wedge S^{n-1}\wedge S^{m-1}\wedge S^1\wedge A\wedge X\end{array}}
\ar@/^40pt/[d]_{\;}="a"^{\S^{n+m-1}(f_2\sa g_2)}\ar@/_40pt/[d]^{\;}="b"_{\S^{n+m-1}(f_1\sa g_1)}\\
{\begin{array}{c}S^1\wedge S^{n-1}\wedge S^{m-1}\wedge S^1\wedge B\wedge Y\\
\cong^{(1\;2)}\\
S^{n-1}\wedge S^1\wedge S^{m-1}\wedge S^{1}\wedge B\wedge Y
\end{array}}\ar@{=>}"b";"a"^{\S(F\sa G)}}
\end{equation*}
and similarly $\S(F\sc G)$. Notice that the last spherical coordinate remains always in the same place in this 
diagram. This is again relevant in connection with Lemma \ref{nilson} below.

\section{The construction of the smash product operation}

In this section we define the smash product morphism in $\C{qpm}$ which appears in the statement of Theorem
\ref{lamel}. In the last two sections we establish the properties which show that the definition given here is
indeed consistent with the definition of the tensor product of quadratic pair modules.

The secondary homotopy groups $\pi_{n,*}X$, $n\geq 0$, of a pointed
space $X$ were introduced in \cite{2hg1}.

For $n=0$, $\pi_{0,*}X$ is the fundamental pointed groupoid of $X$.
We denote by $\pi_{0,0}X$ to the pointed set of objects, which can
be regarded as the set of pointed maps $S^0\r X$, and by
$\pi_{0,1}X$ to the set of morphisms. Such a morphism $\alpha\colon
x\r y$ is a track $\alpha\colon x\rr y$ between pointed maps
$x,y\colon S^0\r Y$.

For $n=1$, $\pi_{1,*}X$ is a crossed module
$$\partial\colon\pi_{1,1}X\To\pi_{1,0}X=\grupo{\L X}.$$ In particular $\pi_{1,0}X$ acts on the right on $\pi_{1,1}X$.

For $n\geq 2$, $\pi_{n,*}X$ is a reduced quadratic module in the
sense of \cite{ch4c}
\begin{equation}\label{semista}
\otimes^2(\pi_{n,0}X)_\abb\st{\omega}\To\pi_{n,1}X\st{\partial}\To\pi_{n,0}X=\grupo{\L^nX}_\nill,
\end{equation}
which is stable for $n\geq 3$.

For all $n\geq 1$ the elements of $\pi_{n,1}X$ are equivalence
classes $[f,F]$ represented by a map
$$f\colon S^1\To\vee_{\L^nX}S^1$$
and a track
$$\xymatrix{S^n\ar[r]_-{\S^{n-1}f}^<(.98){\;\;\;\;\;}="a"\ar@/^25pt/[rr]^0_{}="b"&\vee_{\L^nX}S^n
\ar[r]_-{ev}&X\ar@{=>}"a";"b"_F}$$ where $ev$ is the obvious
evaluation map. Recall from Section \ref{el1} the notation $S^n_X=\vee_{\L^nX}S^n$. Two elements $[f,F],
[g,G]\in\pi_{n,1}X$ coincide provided there is a diagram like (\ref{algo}) parting to the trivial track
with $\hopf(N)=0$ for $n\geq 2$ and no conditions on $N$ for $n=1$. We refer the reader to \cite{2hg1} for
further details on the construction of the algebraic structure of $\pi_{n,*}X$.


According to the definition of additive secondary homotopy groups of a pointed space given in Remark
\ref{remel} the quadratic pair module $\Pi_{n,*}X$ looks as follows.

\begin{equation}\label{parese}
\Pi_{n,*}X=\left(\begin{array}{c}\xymatrix@C=0pt{&**[r]\Pi_{n,ee}X=\otimes^2\Z[\L^nX]\ar[ld]_{P}&\\\Pi_{n,1}X\ar[rr]_\partial&&**[r]\Pi_{n,0}X=\grupo{\L^nX}_\nill\ar[lu]_H}\end{array}\right)
\end{equation}
Here $H$ is always defined as in (\ref{znil}).



For $n\geq 3$, $\Pi_{n,1}X=\pi_{n,1}X$ and $P=\omega\tau_\otimes$ and $\partial$ in
(\ref{parese}) are given by the homomorphisms in
(\ref{semista}).

For $n=2$ the group $\Pi_{2,1}X$ is the quotient of $\pi_{2,1}X$ by
the relations
$$P(a\otimes b+b\otimes a)=0;\;\; a,b\in(\pi_{2,0}X)_\abb;$$
and $P$ and $\partial$ in (\ref{parese}) are induced by $\omega\tau_\otimes$ and $\partial$ in (\ref{semista}) respectively.

For $n=1$, $\Pi_{1,1}X$ is the quotient of the group
$$\pi_{1,1}X\times(\hat{\otimes}^2(\pi_{1,0}X)_\abb)$$
by the relations
$$(-[f,F]+[f,F]^x,0)=(0,x\hat{\otimes}\partial[f,F]);\;\;[f,F]\in\pi_{1,1}X,\; x\in\pi_{1,0}X;$$
$P(a\otimes b)=(0,a\hat{\otimes} b)$ for $a,b\in (\pi_{1,0}X)_\abb$; and
$\partial([f,F],x\hat{\otimes}y)=\partial[f,F]-y-x+y+x$ in $\grupo{\L X}_\nill$ for
$[f,F]\in\pi_{1,1}X$ and $x,y\in\pi_{1,0}X$.

Finally for $n=0$ the group $\Pi_{0,1}X$ is the quotient of
$$\grupo{\pi_{0,1}X}_\nill\times(\hat{\otimes}^2\Z[\pi_{0,0}X])$$
by the relations
$$([\alpha,\alpha'],0)=(0,(-x'+y')\hat{\otimes}(-x+y))$$ for all morphisms $\alpha\colon x\r y$
and $\alpha'\colon x'\r y'$ in $\pi_{0,*}X$,
$$(\alpha\beta,0)=(\beta+\alpha,0),$$ for all composable morphisms $\bullet\st{\beta}\r
\bullet\st{\alpha}\r\bullet$ in $\pi_{0,*}X$, 
$P(a\otimes
b)=(0,a\hat{\otimes} b)$ for $a,b\in \pi_{0,0}X$; and
$\partial(\alpha,a\hat{\otimes} b)=-x+y+[b,a]$ for $a,b\in
\pi_{0,0}X$ and $\alpha\colon x\r y$ in $\pi_{0,1}X$.

We denote by $\Pi_{n,(0)}X$ and $\Pi_{n,(1)}X$ the square groups
$$\Pi_{n,(0)}X=(\Pi_{n,0}X\mathop{\rightleftarrows}\limits_{\partial P}^{H}\Pi_{n,ee}X),$$
$$\Pi_{n,(1)}X=(\Pi_{n,1}X\mathop{\rightleftarrows}\limits_{P}^{H\partial}\Pi_{n,ee}X),$$
defining the quadratic pair module $\Pi_{n,*}X$.

\begin{defn}\label{ka}
The \emph{smash product operation} for the additive secondary
homotopy groups of two pointed spaces $X, Y$ is given by morphisms
in $\C{qpm}$, $n,m\geq 0$,
\begin{equation}\label{ella}
\Pi_{n,*}X\odot\Pi_{m,*}Y\st{\wedge}\To\Pi_{n+m,*}(X\wedge Y).
\end{equation}
These morphisms are induced by square group morphisms, $n,m\geq 0$,
$0\leq i,j,i+j\leq 1$,
\begin{equation}\label{ellos}
\Pi_{n,(i)}X\odot\Pi_{m,(j)}Y\st{\wedge}\To\Pi_{n+m,(i+j)}(X\wedge
Y).
\end{equation}
defined as follows:

For $i=j=0$ morphism (\ref{ellos}) is the composition
$$\Z_\nill[\L^nX]\odot\Z_\nill[\L^mY]\cong\Z_\nill[(\L^nX)\wedge(\L^mY)]
\st{\Z_\nill[\wedge]}\To\Z_\nill[\L^{n+m}(X\wedge Y)]$$ of the
isomorphism in Lemma \ref{facil} and the morphism induced by the map
between discrete pointed sets
\begin{equation}\label{splash}
\wedge\colon(\L^nX)\wedge(\L^m Y)\To\L^{m+n}(X\wedge
Y),\;\;n,m\geq0,
\end{equation}
defined by
$$(f\colon S^n \r X)\wedge (g\colon S^m\r Y)\mapsto (f\wedge g\colon
S^{n+m}\r X\wedge Y).$$

All morphisms in (\ref{ellos}) coincide in the $ee$-term.

Suppose now that $n,m\geq 1$.

For $i=0$ and $j=1$, an element $g\ca[f,F]$ with $g\in\Pi_{n,0}X$
and $[f,F]\in\Pi_{m,1}Y$ is sent by (\ref{ellos}) to the element
$g\wa[f,F]\in\Pi_{n+m,1}(X\wedge Y)$ represented by the map
$$S^1\st{\bar{g}\sa  f}\To \vee_{(\L^nX)\wedge(\L^mY)}S^1\st{\S\wedge}\To\vee_{\L^{n+m}(X\wedge Y)}S^1,$$
where $\bar{g}\colon S^1\r\vee_{\L^nX} S^1$ is any map with
$(\pi_1\bar{g})_\nill(1)=g$ and the second arrow is the suspension of
$\wedge$ in (\ref{splash}), and by the track
$$\xymatrix{S^{n+m}\ar[rr]_-{\S^{n+m-1}(\bar{g}\sa f)}^{\;}="a"\ar@/^20pt/[rr]_{\;}="b"^<(.7){\S^{n-1}\bar{g}\wedge\S^{m-1}f}
\ar@/^50pt/[rrrrr]_<(.583){\;}="d"^0&&
\vee_{\L^nX\wedge\L^mY}S^{n+m}\ar[rr]_-{\S^{n+m}\wedge}^{\;}="c"&&S^{n+m}_{X\wedge Y}\ar[r]_-{ev}&X\wedge Y
\ar@{=>}"a";"b"_{\aleph^{\sa}_{n,\bar{g},m,f}}
\ar@{=>}"c";"d"_{\bar{g}_{ev}\wedge F}}$$
i. e.
$$g\wa[f,F]=[(\S\wedge)(\bar{g}\sa f), (\bar{g}_{ev}\wedge F)\vc(ev(\S^{n+m}\wedge)\aleph^{\sa}_{n,\bar{g},m,f})].$$ 
Here, and in the following three cases, the smash products of maps and tracks need to be altered by permutations according to (\ref{perara}).

In a
similar way the element $g\cc[f,F]$ is sent by (\ref{ellos}) to
$g{\wc}[f,F]\in\Pi_{n+m,1}(X\wedge Y)$, given by the map
$$S^1\st{\bar{g}\sc  f}\To \vee_{(\L^nX)\wedge(\L^mY)}S^1\st{\S\wedge}\To\vee_{\L^{n+m}(X\wedge Y)}S^1,$$
and by the track
$$\xymatrix{S^{n+m}\ar[rr]_-{\S^{n+m-1}(\bar{g}\sc
f)}^{\;}="a"\ar@/^20pt/[rr]_{\;}="b"^<(.7){\S^{n-1}\bar{g}\wedge\S^{m-1}f}
\ar@/^50pt/[rrrrr]_<(.583){\;}="d"^0&&
\vee_{\L^nX\wedge\L^mY}S^{n+m}\ar[rr]_-{\S^{n+m}\wedge}^{\;}="c"&&S^{n+m}_{X\wedge Y}\ar[r]_-{ev}&X\wedge Y
\ar@{=>}"a";"b"_{\aleph^{\sc}_{n,\bar{g},m,f}}
\ar@{=>}"c";"d"_{\bar{g}_{ev}\wedge F}}$$
i. e.
$$g{\wc}[f,F]=[(\S\wedge)(\bar{g}\sc f), (\bar{g}_{ev}\wedge
F)\vc(ev(\S^{n+m}\wedge)\aleph^{\sc}_{n,\bar{g},m,f})].$$ 

For $i=1$ and $j=0$ the generator $[f,F]\ca g$ with
$[f,F]\in\Pi_{n,1}X$ and $g\in\Pi_{m,0}Y$ is sent by (\ref{ellos})
to the element $[f,F]{\wa}g\in\Pi_{n+m,1}(X\wedge Y)$ represented by
the map
$$S^1\st{f\sa \bar{g}}\To\vee_{(\L^nX)\wedge(\L^m Y)}S^1\st{\S\wedge}\To\vee_{\L^{n+m}(X\wedge Y)}S^1,$$
where $\bar{g}\colon S^1\r\vee_{\L^mY}S^1$ is any map with
$(\pi_1\bar{g})_\nill(1)=g$, and the track
$$\xymatrix{S^{n+m}\ar[rr]_-{\S^{n+m-1}(f\sa
\bar{g})}^{\;}="a"\ar@/^20pt/[rr]_{\;}="b"^<(.7){\S^{n-1}f\wedge\S^{m-1}\bar{g}}
\ar@/^50pt/[rrrrr]_<(.583){\;}="d"^0&&
\vee_{\L^nX\wedge\L^mY}S^{n+m}\ar[rr]_-{\S^{n+m}\wedge}^{\;}="c"&&S^{n+m}_{X\wedge Y}\ar[r]_-{ev}&X\wedge Y
\ar@{=>}"a";"b"_{\aleph^{\sa}_{n,f,m,\bar{g}}}
\ar@{=>}"c";"d"_{F\wedge \bar{g}_{ev}}}$$
i. e.
$$[f,F]\wa g=[(\S\wedge)(f\sa\bar{g}), 
(F\wedge \bar{g}_{ev})\vc(ev(\S^{n+m}\wedge)\aleph^{\sa}_{n,f,m,\bar{g}})].$$  The
element $[f,F]\cc g$ is sent by (\ref{ellos}) to $[f,F]{\wc}
g\in\Pi_{n+m,1}(X\wedge Y)$ given by the map
$$S^1\st{f\sc \bar{g}}\To\vee_{(\L^nX)\wedge(\L^m Y)}S^1\st{\S\wedge}\To\vee_{\L^{n+m}(X\wedge Y)}S^1$$
and the track
$$\xymatrix{S^{n+m}\ar[rr]_-{\S^{n+m-1}(f\sc
\bar{g})}^{\;}="a"\ar@/^20pt/[rr]_{\;}="b"^<(.7){\S^{n-1}f\wedge\S^{m-1}\bar{g}}
\ar@/^50pt/[rrrrr]_<(.583){\;}="d"^0&&
\vee_{\L^nX\wedge\L^mY}S^{n+m}\ar[rr]_-{\S^{n+m}\wedge}^{\;}="c"&&S^{n+m}_{X\wedge Y}\ar[r]_-{ev}&X\wedge Y
\ar@{=>}"a";"b"_{\aleph^{\sc}_{n,f,m,\bar{g}}}
\ar@{=>}"c";"d"_{F\wedge \bar{g}_{ev}}}$$
i. e.
$$[f,F]{\wc} g=[(\S\wedge)(f\sc\bar{g}),
(F\wedge \bar{g}_{ev})\vc(ev(\S^{n+m}\wedge)\aleph^{\sc}_{n,f,m,\bar{g}})].$$  

Suppose now that $n=0$ and $m\geq 0$.

For $i=0$ and $j=1$, an element $g\ca f$ with $g\colon S^0\r X$ and
$f\in\Pi_{m,1}Y$ is sent by (\ref{ellos}) to
$g\wa f=(\Pi_{m,1}(g\wedge Y))(f)\in\Pi_{m,1}(X\wedge Y)$. 

For $i=1$ and $j=0$, an element $F\ca g$ with $F\colon f\rr f'$ a
track between maps $f,f'\colon S^0\r X$ and $g\in\Pi_{m,0}Y$ is sent
by (\ref{ellos}) to $F\wa g=(\Pi_{m,*}(F\wedge Y))(g)\Pi_{m,1}(X\wedge Y)$.

Suppose now that $n\geq0$ and $m=0$.

For $i=0$ and $j=1$, an element $g\cc F$ with $g\in\Pi_{n,0}X$ and
$F\colon f\rr f'$ a track between maps $f,f'\colon S^0\r Y$   is
sent by (\ref{ellos}) to
$g\wc F=(\Pi_{n,*}(X\wedge F))(g)\Pi_{n,1}(X\wedge Y)$. 

For $i=1$ and $j=0$, an element $f\cc g$ with $f\in\Pi_{n,1}X$ and
$g\colon S^0\r Y$  is sent by (\ref{ellos}) to $f\wc
g=(\Pi_{m,1}(X\wedge g))(f)\Pi_{n,1}(X\wedge Y)$.
\end{defn}

\section{The Hopf invariant for tracks and smash products}

In this section we prove two lemmas on the Hopf invariant for tracks defined in \cite{2hg1} 3.3
which will be useful to check the
properties of the smash product operation on secondary homotopy groups.

The first lemma computes the Hopf invariant of a track smashed with a discrete set.

\begin{lem}\label{hopfs}
Let $f,g\colon S^1\wedge A\r S^1\wedge B$ be maps between
suspensions of discrete pointed sets $A, B$; let $F\colon S^1\wedge
f\rr S^1\wedge g$ be a track; and let $X$ be another discrete
pointed set. Then the following equations hold
\begin{enumerate}
\item $\hopf(F\wedge X)=(1\otimes \tau_\otimes\otimes 1)(\hopf(F)\otimes\Delta)$,
\item $\hopf (X\wedge F)=(1\otimes \tau_\otimes\otimes 1)(\Delta\otimes\hopf(F))$.
\end{enumerate}
In particular the smash product of a track with trivial Hopf invariant and a discrete
pointed set has always a trivial Hopf invariant.
\end{lem}

This lemma follows easily from the elementary properties of
the Hopf invariant for tracks in \cite{2hg1} 3.

The second lemma computes the effect of conjugation by an automorphism of a sphere on certain tracks.

\begin{lem}\label{nilson}
Let $f,g\colon S^1\wedge A\r S^1\wedge B$ be maps between
suspensions of discrete pointed sets $A, B$; let $F\colon
S^{n-1}\wedge f\rr S^{n-1}\wedge g$ be a track; and let
$\alpha\colon S^{n-1}\cong S^{n-1}$ be a homeomorphism. Then
$(\alpha\wedge S^1\wedge B)F(\alpha^{-1}\wedge S^1\wedge A)$ is a
track with trivial Hopf invariant if and only if $F$ has trivial Hopf invariant. Moreover, if $\alpha$
has degree $1$ or $n\geq 3$ then
$$F=(\alpha\wedge S^1\wedge B)F(\alpha^{-1}\wedge S^1\wedge A).$$
\end{lem}

\begin{proof}
Let us denote also by $F\colon I_+\wedge S^{n-1}\wedge S^1\wedge A\r
S^{n-1}\wedge S^1\wedge B$ to a map representing the track $F$. The
adjoint map of pairs
$$\begin{array}{l}
ad((\alpha\wedge S^1\wedge B)F(I_+\wedge\alpha^{-1}\wedge S^1\wedge A))
\colon\\ (I_+\wedge S^1\wedge A,(S^1\wedge A)\vee(S^1\wedge A))\To
(\L^{n-1}(S^{n-1}\wedge S^1\wedge B),S^1\wedge B)
\end{array}$$
used to define the Hopf invariant coincides with the composite
$$\xymatrix{(I_+\wedge S^1\wedge A,(S^1\wedge A)\vee(S^1\wedge A))\ar[d]^{ad(F)}\\
(\L^{n-1}(S^{n-1}\wedge S^1\wedge B),S^1\wedge B)
\ar[d]^{\operatorname{map}_*(\alpha^{-1},\alpha\wedge S^1\wedge B)}\\
(\L^{n-1}(S^{n-1}\wedge S^1\wedge B),S^1\wedge B)}$$ The map
$\operatorname{map}_*(\alpha^{-1},\alpha\wedge S^1\wedge B)$ is a
homeomorphism, hence the first part of the lemma follows from the
very definition of the Hopf invariant for tracks.

The homeomorphism $\operatorname{map}_*(\alpha^{-1},\alpha\wedge
S^1\wedge B)$ restricts to the identity on $S^1\wedge B$. Moreover,
if $\alpha$ has degree $1$ or $n\geq 3$ then this homeomorphism is
compatible with the $H$-multiplication of the $(n-1)$-fold loop
space up to homotopy, and therefore with the Pontrjagin product. In
particular by the definition and elementary properties of the Hopf
invariant for tracks
$$\hopf(F)=\hopf((\alpha\wedge S^1\wedge B)F(\alpha^{-1}\wedge S^1\wedge A))$$
and hence
$$F=(\alpha\wedge S^1\wedge B)F(\alpha^{-1}\wedge S^1\wedge A).$$
\end{proof}

\section{Hopf invariant computations related to exterior tracks}

In this section we perform two Hopf invariant computations for tracks. The first computation is 
connected with axiom (7) in the definition
of the tensor product of square groups, see Definition \ref{tpsg}, and the second one is connected 
with the commutativity rule for the smash
product operation on additive secondary homotopy groups, see Theorem \ref{nosym}. Both computations are
crucial steps towards the proof of the main results of this paper stated in Section \ref{LA}. They show that the
algebraic structures described in Parts \ref{I} and \ref{II} are the right algebraic structures to describe the 
smash product operation.

First of all we define a concept which will be useful for computations.

\begin{defn}\label{H}
Let $\C{nil}$ be the category of free groups of nilpotency class $2$ and let $\Phi\colon\C{Ab}\r\C{Ab}$ be a functor.
The  \emph{$\Hg$-group} of $\Phi$ is the class $\Hg(\Phi)$ of all functions
$\chi$ sending two morphisms $f\colon\grupo{A}_\nill\r\grupo{B}_\nill$, 
$g\colon\grupo{X}_\nill\r\grupo{Y}_\nill$ in $\C{nil}$ to a homomorphism
$$\chi(f,g)\colon\Z[A]\otimes\Z[X]\To\Phi(\Z[B]\otimes\Z[Y])$$
in such a way that if $f,f_i,g,g_i$ are morphisms in $\C{nil}$ and $f',f'_i,g',g'_i$ are maps between pointed sets $i=1,2$ then
\begin{enumerate}
\item $\chi (\grupo{f'}_\nill,g)=0$,
\item $\chi (f,\grupo{g'}_\nill)=0$,
\item $\chi ((f_1,f_2),g)=(\chi (f_1,g),\chi (f_2,g))$,
\item $\chi (f,(g_1,g_2))=(\chi (f,g_1),\chi (f,g_2))$,
\item $\chi (f_1\grupo{f_2'}_\nill,g_1\grupo{g_2'}_\nill)=\chi (f_1,g_1)(\Z[f_2']\otimes\Z[g_2'])$,
\item $\chi (\grupo{f_1'}_\nill f_2,\grupo{g_1'}_\nill g_2)=\Phi(\Z[f_1']\otimes\Z[g_1'])\chi (f_2,g_2)$.
\end{enumerate}
A natural transformation $\zeta\colon \Phi\rr\Psi$ between functors $\Phi,\Psi\colon\C{Ab}\r\C{Ab}$  induces a
function $\Hg(\zeta)\colon\Hg(\Phi)\r\Hg(\Psi)$ in the obvious way.

If $\Hg(\Phi)$ is a set then it is an abelian group by addition of abelian group homomorphisms. If
$\Hg(\Psi)$ is also a set then $\Hg(\zeta)$ is an abelian group homomorphism. 

Many functors have a $\Hg$-group which is a set, see for example Lemma \ref{unai} below. Alternatively one can
define $\Hg$-groups by using a small subcategory of $\C{nil}$ to obtain always sets. We therefore
do not care about set theoretic subtleties in what follows.
\end{defn}

The following lemma shows examples of non-trivial elements in the $\Hg$-group of the reduced tensor
square.

\begin{lem}
There are elements, $n,m\geq 0$,
$$\bar{\sigma}\tau_\otimes(1\otimes\tau_\otimes\otimes1)(H\otimes TH),\;\;-\bar{\sigma}\binom{(-1)^{nm}}{2}H(\sc)\in \Hg(\hat{\otimes}^2)$$
which evaluated at $f\colon \grupo{A}_\nill\r\grupo{B}_\nill$, $g\colon\grupo{X}_\nill\r\grupo{Y}_\nill$
send an element $a\otimes x\in\Z[A]\otimes\Z[X]$ with $a\in A$ and $x\in X$ to
$$\bar{\sigma}\tau_\otimes(1\otimes\tau_\otimes\otimes1)(H(f(a))\otimes TH(g(x))),$$
$$-\bar{\sigma}\binom{(-1)^{nm}}{2}H(f(a)\sc g(x)),$$
in $\hat{\otimes}^2(\Z[B]\otimes\Z[Y])$ respectively.
\end{lem}

Properties (1)--(6) in Definition \ref{H} are easy to check in these cases.

The following lemma is left as an exercise for the reader.

\begin{lem}
Given three functors $\Phi,\Psi,\Gamma\colon\C{Ab}\r\C{Ab}$ and a natural exact sequence
$$\Phi\st{\zeta}\hookrightarrow\Psi\st{\xi}\To\Gamma,$$
the sequence
$$\Hg(\Phi)\st{\Hg(\zeta)}\hookrightarrow\Hg(\xi)\st{\Hg(\xi)}\To\Hg(\Gamma)$$
is exact.
\end{lem}

This lemma can be applied to the natural exact sequence
$$A\otimes\Z/2\st{\bar{\tau}}\hookrightarrow\hat{\otimes}^2A\st{q}\twoheadrightarrow\wedge^2A,$$
where $\wedge^2A$ is the exterior square of $A$, $\bar{\tau}(a)=\bar{\sigma}(a\otimes a)$ and
$q\bar{\sigma}(a\otimes b)=a\wedge b$.

Now we define elements in the $\Hg$-group of the reduced tensor square by using the exterior tracks and the Hopf invriant for tracks.

In the rest of this section $A, B, X, Y$ will always be pointed discrete sets.
Given maps $f\colon\S A\r\S B$, $g\colon\S X\r \S Y$, 
we define the following abelian group homomorphism as a Hopf
invariant for tracks
\begin{equation}\label{laK}
K(f,g)=\bar{\sigma}\hopf((\aleph^{\sc}_{f,g})^\vi\vc{\aleph}^{\sa}_{f,g})\colon
\Z[A]\otimes\Z[X]\To\hat{\otimes}^2(\Z[B]\otimes\Z[Y]).
\end{equation}

\begin{prop}\label{Kprop}
The homomorphism $K(f,g)$ defined above only depends on 
$(\pi_1f)_\nill$ and $(\pi_1g)_\nill$. Moreover, $K\in\Hg(\hat{\otimes}^2)$. Furthermore, 
$$\Hg(q)(K)=\Hg(q)(\bar{\sigma}\tau_\otimes(1\otimes\tau_\otimes\otimes1)(H\otimes TH)).$$
\end{prop}

\begin{proof}
By \cite{2hg1} 3.6 (5) and Lemma \ref{nilson}
\begin{equation*}\tag{a}
K(f,g)=\hopf((\aleph^{\sc}_{n,f,m,g})^\vi\vc{\aleph}^{\sa}_{n,f,m,g})
\end{equation*}
for $n,m\geq1$ with $n+m>2$. If $(\pi_1f)_\nill=(\pi_1\bar{f})_\nill$ and
$(\pi_1g)_\nill=(\pi_1\bar{g})_\nill$ then there are tracks $F\colon
\S^{n-1} f\rr\S^{n-1}\bar{f}$, $G\colon \S^{m-1}
g\rr\S^{m-1}\bar{g}$ with trivial Hopf invariant. Moreover, the tracks
$$\S(F\sa G)\colon \S^{n+m-1}(f\sa g)\rr
\S^{n+m-1}(\bar{f}\sa\bar{g}),$$
$$\S(F\sc G)\colon \S^{n+m-1}(f\sc g)\rr
\S^{n+m-1}(\bar{f}\sc\bar{g}),$$
have trivial Hopf invariant by Lemmas \ref{hopfs} and \ref{nilson}, and
\cite{2hg1} 3.6. 

By Lemma \ref{exteprop2} 
we have that
\begin{eqnarray*}
\S(F\sa G)^\vi
\vc(\aleph^{\sc}_{n,\bar{f},m,\bar{g}})^\vi
\vc{\aleph}^{\sa}_{n,\bar{f},m,\bar{g}}
\vc\S(F\sc G)&=&
(\aleph^{\sc}_{n,f,m,g})^\vi\vc(F\wedge G)^\vi
\\&&\vc(F\wedge G)\vc{\aleph}^{\sa}_{n,f,m,g}\\
&=&(\aleph^{\sc}_{n,f,m,g})^\vi\vc{\aleph}^{\sa}_{n,f,m,g}.
\end{eqnarray*}
Hence the first part of the statement follows from the elementary properties of the Hopf invariant for tracks in
\cite{2hg1} 3.

The equation for the images by $q$ follows from \cite{2hg1} 3.6, equation (a),
and the equality
\begin{eqnarray*}
(\pi_1(f\sc g))_\nill(a\wedge x)&=&(\pi_1(f\sa g))_\nill(a\wedge x)\\&&
+\partial\tau_\otimes(1\otimes\tau_\otimes\otimes1)(H((\pi_1f)_\nill(a))\otimes TH((\pi_1g)_\nill(x))).
\end{eqnarray*} 
This last
equality is a consequence of Definition \ref{tpsg} (7) and
Proposition \ref{facil}.

Finally (1)--(6) in Definition \ref{H} follow from Lemma \ref{exteprop} (1)--(10)
and the elementary properties of the Hopf invariant for tracks in
\cite{2hg1} 3.6.
\end{proof}

Let $\hat{\tau}_{n,m}\colon\tau_{n,m}\rr (\cdot)^{(-1)^{nm}}_{n+m}$ be a track between maps 
$$(\cdot)^{(-1)^{nm}}_{n+m},\tau_{n,m}\colon S^{n+m}\r S^{n+m}.$$
This is a lift of the shuffle permutation in (\ref{chuf}) to the symmetric track group $\symt{n+m}$.
Given maps $f\colon\S A\r\S B$, $g\colon\S X\r\S Y$, we define the following Hopf invariants 
\begin{eqnarray}
\label{laL}&&\\
\nonumber L_{1,1}(f,g)&=&
\bar{\sigma}\hopf((\S^2\tau_\wedge)(\aleph^{\sc}_{g,f})^\vi\hat{\tau}_{1,1})\vc
((\hat{\tau}_{1,1}^\vi\wedge B\wedge Y)\aleph^\sa_{f,g})),\\
\nonumber
L_{n,m}(f,g)&=&
\hopf((\S^{m+n}\tau_\wedge)(\aleph^{\sc}_{m,g,n,f})^\vi\hat{\tau}_{n,m})\vc
((\hat{\tau}_{n,m}^\vi\wedge B\wedge Y)\aleph^\sa_{n,f,m,g})).
\end{eqnarray}
Here $n,m\geq 1$ and $n>1$ or $m>1$, and $\tau_\wedge\colon Y\wedge B\r B\wedge Y$ is the symmetry map for the smash
product. Notice that $L_{n,m}(f,g)$ is a homomorphism, $n,m\geq 1$,
$$L_{n,m}(f,g)\colon\Z[A]\otimes\Z[X]\To\hat{\otimes}^2(\Z[B]\otimes\Z[Y]).$$

\begin{prop}\label{Lprop}
For any $n,m\geq1$ the homomorphism $L_{n,m}(f,g)$ defined above does not depend on the choice of
$\hat{\tau}_{n,m}$. Moreover, it only depends on 
$(\pi_1f)_\nill$ and $(\pi_1g)_\nill$. Furthermore, $L_{n,m}\in\Hg(\hat{\otimes}^2)$. In addition, 
$$\Hg(q)(L_{n,m})=\Hg(q)\left(-\bar{\sigma}\binom{(-1)^{nm}}{2}H(\sc)\right).$$
\end{prop}

The proof of this proposition is analogous to Proposition \ref{Kprop}.

We would like to omit $\Hg(q)$ in the equations of Propositions \ref{Kprop} and \ref{Lprop}. For this we prove the following two lemmas.

\begin{lem}\label{unai}
If $\Phi\colon\C{Ab}\r\C{Ab}$ is an additive functor which preserves arbitrary filtered colimits 
then there is an isomorphism
$$\Hg(\Phi)\cong\Phi(\Z).$$
\end{lem}

\begin{proof}
We claim that the isomorphism sends $\chi$ to $\chi(-1,-1)\in\Phi(\Z)$, where $-1$ is the homomorphism $-1\colon\Z\r\Z$. This will be a
consequence of the following formula, that we claim to hold. 
We first notice that
$$\Phi(\Z[B]\otimes\Z[Y])=\Phi(\Z)\otimes\Z[B]\otimes\Z[Y].$$
Given $f\colon \grupo{A}_\nill\r\grupo{B}_\nill$, 
$g\colon\grupo{X}_\nill\r\grupo{Y}_\nill$, $a\in A$ and $x\in X$, if $f_\abb(a)=\sum_i n_ib_i$,
$g_\abb(x)=\sum_jm_jy_j$ are the linear expansions with $b_i\in B$, $y_j\in Y$, $n_i,m_j\in\Z$, then
\begin{equation*}\tag{a}
\chi(f,g)(a\otimes x)=\sum_{i,j}\epsilon(n_i,m_j)\abs{n_im_j}\chi(-1,-1)\otimes b_i\otimes y_j.
\end{equation*}
Here $\epsilon(n_i,m_j)=1$ provided $n_i,m_j<0$ and it is zero otherwise. Conversely any $\chi$
defined by formula (a) with $\chi(-1,-1)$ arbitrarily chosen out of $\Phi(\Z)$ defines an element in $\Hg(\Phi)$. We will just prove the first part of the
claim, the converse is easy.

By (3), (4), (5) and (6) in Definition \ref{H} we have
\begin{equation*}\tag{b}
\chi(f,g)(a\otimes x)=\sum_{i,j}\chi(n_i,m_j)\otimes b_i\otimes y_j.
\end{equation*}
for $n_i,m_j\colon\Z\r\Z$. 
For any $n\in\Z$ we consider the homomorphism $\mu_n\colon\Z\r\grupo{c_1,\dots,c_{\abs{n}}}_\nill$ defined by
$\mu_n(1)=c_1+\cdots+c_{\abs{n}}$ if $n>0$, $\mu_n(1)=-c_1-\cdots-c_{\abs{n}}$ if $n<0$, and $\mu_n(1)=0$ if $n=0$.
By (b) given $n,m\in\Z$ we have that 
\begin{equation*}\tag{c}
\chi(\mu_n,\mu_m)=\sum_{i=1}^{\abs{n}}\sum_{j=1}^{\abs{m}}\chi(\epsilon(n),\epsilon(m))\otimes c_i\otimes c_j.
\end{equation*}
Here $\epsilon$ sends a positive integer to $1$, a negative integer to $-1$, and zero to itself.
Now, again by Definition \ref{H} (6), we have that $\chi(n,m)$ is the sum of all coefficients in (c), i. e. 
$$\chi(n,m)=\abs{nm}\chi(\epsilon(n),\epsilon(m))\in\Phi(\Z).$$
But by Definition \ref{H} (1) and (2) the element $\chi(\epsilon(n),\epsilon(m))\in\Phi(\Z)$ is zero unless
$n,m<0$, hence we are done.
\end{proof}

\begin{lem}\label{papi}
There is a commutative diagram
$$\xymatrix{\Hg(-\otimes\Z/2)\ar@{^{(}->}[r]\ar[d]_\cong&\Hg(\hat{\otimes}^2)\ar@/^10pt/[ld]\\
\Z/2&}$$
\end{lem}

\begin{proof}
Both vertical homomorphisms send an element $\chi$ in the corresponding $\Hg$-group to
$$\chi(-1,-1)\in\Z\otimes\Z/2=\hat{\otimes}^2\Z=\Z/2.$$
The isomorphism is proved in the previous lemma. 
\end{proof}

The following theorem is a key step towards the proof of Theorem 6.1. It is connected to relation (7) in the
definition of the tensor product for square groups, see Definition \ref{tpsg}.

\begin{thm}\label{hopfcal}
$K=\bar{\sigma}\tau_\otimes(1\otimes \tau_\otimes\otimes 1)(H\otimes
TH)\in\Hg(\hat{\otimes}^2).$
\end{thm}

\begin{proof}
By Proposition \ref{Kprop} and Lemma \ref{papi} we only have to check that
$$\bar{\sigma}\tau_\otimes(1\otimes \tau_\otimes\otimes 1)(H(-1)\otimes
TH(-1))=K(\nu,\nu)$$
for $\nu\colon S^1\r S^1\colon z\mapsto z^{-1}$ the complex inversion. It is easy to see from the very definition
of $H$ that 
$$\bar{\sigma}\tau_\otimes(1\otimes \tau_\otimes\otimes 1)(H(-1)\otimes
TH(-1))=1\in\Z/2,$$
hence this theorem follows from Lemma \ref{K-1} below.
\end{proof}

The next theorem is the main step in the proof of the commutativity rule for the smash product operation in Theorem \ref{nosym}.

\begin{thm}\label{hopfcal2}
$L_{n,m}=-\bar{\sigma}\binom{(-1)^{nm}}{2}H(\sc)\in\Hg(\hat{\otimes}^2)$, $n,m\geq1$.
\end{thm}

\begin{proof}
By Proposition \ref{Lprop} and Lemma \ref{papi} we only have to check that
$$-\bar{\sigma}\binom{(-1)^{nm}}{2}H((-1)\sc(-1))=L(\nu,\nu)$$
for $\nu\colon S^1\r S^1\colon z\mapsto z^{-1}$ the complex inversion. But $(-1)\sc(-1)=1$ and
$H(1)=\binom{1}{2}=0$,
hence this theorem follows from Lemma \ref{L-1} below.
\end{proof}

\section{Hopf invariants of tracks between orthogonal transformations}

In this section we are concerned with the computation of Hopf invariants for tracks between self maps of spheres
$S^n$, $n\geq 2$. More concretely, we are interested in tracks between maps $A\colon S^n\r S^n$ which are induced by the left action
of the orthogonal group $O(n)$ on $S^n$, i. e. $A\in O(n)$. The pull-back of this action along the inclusion
$\sym{n}\subset O(n)$ induced by permutation of coordinates in $\Real^n$ yields the action of $S^n$ already
considered in Section \ref{el1}. Let $\det\colon O(n)\r\set{\pm1}$ be the determinant
homomorphism. We consider the group $\set{\pm1}$ embedded in $O(n)$
as 
$$\left(\begin{array}{cccc}
1&&&\text{{\huge 0}}\\
&\ddots&&\\
&&1&\\
\text{{\huge 0}}&&&\pm 1
\end{array}\right).$$
In \cite{2hg1} 6 we compute the group $\widetilde{O}(n)$ of tracks $A\rr\det A$, $A\in O(n)$, with multiplication
given by horizontal composition. Notice that tracks $A\rr B$ between maps $A,B\colon S^n\r S^n$ with $A,B\in O(n)$ coincide with homotopy classes of paths from $A$ to $B$ in the Lie group $O(n)$ since the $J$-homomorphism $\pi_1O(n)\cong\pi_{n+1}S$ is an isomorphism for $n\geq 2$. We can identify the group
$\widetilde{O}(2)$ with the semidirect product $\set{\pm1}\ltimes\Real$. Here $\Real$ is the additive group of
real numbers and $\set{\pm1}$ is the multiplicative group of order $2$ acting on $\Real$ by multiplication. See
\cite{2hg2} 6.12. 

In order to describe $\widetilde{O}(n)$  for $n\geq 3$ we need to recall the definition of the positive Clifford
algebra. 

\begin{defn}\label{clif}
The \emph{positive Clifford algebra} $C_+(n)$ is the unital $\Real$-algebra
generated by $e_i$, $1\leq i\leq n$, with relations
\begin{enumerate}
\item $e_i^2=1$ for $1\leq i \leq n$,
\item $e_ie_j=-e_je_i$ for $1\leq i<j \leq n$.
\end{enumerate}
Clifford algebras are defined for arbitrary quadratic forms on
finite-dimensional vector spaces, see for instance \cite{rclg} 6.1.
The Clifford algebra defined above corresponds to the quadratic form
of the standard positive-definite scalar product in $\Real^n$. We
identify the sphere $S^{n-1}$ with the vectors of Euclidean norm
$1$ in the vector subspace $\Real^n\subset C_+(n)$ spanned by the
generators $e_i$. The vectors in $S^{n-1}$ are units in $C_+(n)$.
Indeed for any $v\in S^{n-1}$ the square $v^2=1$ is the unit element
in $C_+(n)$, so that $v^{-1}=v$.
The group $\widetilde{O}(n)$ can be identified with the
subgroup of units in $C_+(n)$ generated by $S^{n-1}$. This group is also known as the \emph{positive pin group}.
\end{defn}

The group $\widetilde{O}(n)$ is a covering Lie group of $O(n)$ with simply connected components, and with kernel $\Z\subset\Real$ if $n=2$ and
$\set{\pm1}$ if $n\geq 3$. The covering
homomorphism
$$q\colon \widetilde{O}(n)\To O(n),$$
is defined for $n=2$ as
$$q(x,y)=\left(\begin{array}{cc}
\cos 2\pi y&-\sin 2\pi y\\
x\sin 2\pi y&x\cos 2\pi y
\end{array}\right).$$
For $n\geq 3$ the homomorphism $q$ sends an element $v\in S^{n-1}\subset\widetilde{O}(n)\subset C_+(n)$
to the matrix of the reflexion along the plane orthogonal to the unit vector $v$.
An element $x\in \widetilde{O}(n)$ is identified with the track $q(x)\rr\det q(x)$ determined by the push-forward 
along $q$ of the unique track in $\widetilde{O}(n)$ from the point $x$ to the point
$e_n^{\binom{\det q(x)}{2}}$ if $n\geq 3$, or from $x$ to $(\det q(x) ,0)$ if $n=2$.
With the track approach the covering map $q$ sends a track to the source map.
Moreover, a track $\alpha\colon 1_{S^n}\rr 1_{S^n}$ in the kernel of $q$ is identified with $(1,-\hopf(\alpha))$
for $n=2$,
compare \cite{2hg2} 3.4, and with $\hopf(\alpha)\in\Z/2\cong\set{\pm1}$ for $n\geq 3$.

The suspension of tracks defines group inclusions
$$\S\colon\widetilde{O}(n)\hookrightarrow\widetilde{O}(n+1).$$
For $n\geq 3$ this is induced by the algebra inclusion $C_+(n)\hookrightarrow C_+(n+1)$ defined by $e_j\mapsto
e_{j+1}$, $1\leq j\leq n$. For $n=2$ it is given by 
\begin{equation}\label{inklu}
(x,y)\mapsto e_3^{\binom{x+1}{2}}((\sin \pi y)e_2+(\cos \pi y)e_3).
\end{equation}

For $K$ in (\ref{laK}) ane $\nu\colon S^1\r S^1\colon z\mapsto z^{-1}$ we have the following result.

\begin{lem}\label{K-1}
$K(\nu,\nu)=1$.
\end{lem}

\begin{proof}
We are going to prove the following stronger equality:
\begin{equation*}\tag{a}
\hopf((\aleph^{\sc}_{\nu,\nu})^\vi\vc\aleph^{\sa}_{\nu,\nu})=1.
\end{equation*}

From the very definition of $\aleph$ we obtain the following identity in $\widetilde{O}(2)$.
\begin{eqnarray*}
\aleph(\nu\wedge S^1)&=&(1,\frac{1}{4}).
\end{eqnarray*}
Since $(\nu\wedge S^1)\aleph$ is the inverse of $\aleph(\nu\wedge S^1)$ then
\begin{eqnarray*}
(\nu\wedge S^1)\aleph&=&(1,-\frac{1}{4}).
\end{eqnarray*}
Obviously the identity track $0^\vc_{\S\nu}\colon\S\nu\rr\S\nu$ is $(-1,0)$ in $\widetilde{O}(2)$.
By using these equalities we obtain
\begin{eqnarray*}
(\aleph^{\sa}_{\nu,\nu})^\vi&=&0^\vc_{\S\nu}(\aleph (\nu\wedge S^1))0^\vc_{\S\nu}((\nu\wedge S^1)\aleph )\\
&=&(-1,0)(1,\frac{1}{4})(-1,0)(1,-\frac{1}{4})\\
&=&(-1,\frac{1}{4})(-1,-\frac{1}{4})\\
&=&((-1)(-1),-\frac{1}{4}-\frac{1}{4})\\
&=&(1,-\frac{1}{2}),\\
(\aleph^{\sc}_{\nu,\nu})^\vi&=&
(\aleph (\nu\wedge S^1))0^\vc_{\S\nu}((\nu\wedge S^1)\aleph )0^\vc_{\S\nu}\\
&=&(1,\frac{1}{4})(-1,0)(1,-\frac{1}{4})(-1,0)\\
&=&(-1,-\frac{1}{4})(-1,\frac{1}{4})\\
&=&((-1)(-1),\frac{1}{4}+\frac{1}{4})\\
&=&(1,\frac{1}{2}).
\end{eqnarray*}
Let $\beta\colon 1_{S^2}\rr 1_{S^2}$ be a track with $\hopf(\beta)=1$, so that $\beta=(1,-1)$ in $\widetilde{O}(2)$. 
Equation (a) is equivalent
to the following equation in $\widetilde{O}(2)$,
\begin{eqnarray*}
(\aleph^{\sa}_{\nu,\nu})^\vi&=&(\aleph^{\sc}_{\nu,\nu})^\vi\beta,
\end{eqnarray*}
which follows from the equalities above.
\end{proof}

For $L_{n,m}$ in (\ref{laL}) we have the following result.

\begin{lem}\label{L-1}
$L_{n,m}(\nu,\nu)=0$, $n,m\geq1$.
\end{lem}

\begin{proof}
We showed in the proof of Lemma \ref{K-1} that in $\widetilde{O}(2)$
\begin{eqnarray*}
(\aleph^{\sa}_{\nu,\nu})^\vi&=&(1,-\frac{1}{2}),\\
(\aleph^{\sc}_{\nu,\nu})^\vi&=&(1,\frac{1}{2}).
\end{eqnarray*}

For $n=m=1$ we can take $\hat{\tau}_{1,1}=(-1,-\frac{1}{4})$. Then in $\widetilde{O}(2)$  
\begin{eqnarray*}
\hat{\tau}_{1,1}(\aleph^{\sa}_{\nu,\nu})^\vi&=&(-1,-\frac{1}{4})(1,-\frac{1}{2})\\
&=&(-1,-\frac{1}{4}-\frac{1}{2})\\
&=&(-1,-\frac{3}{4}),\\
(\aleph^{\sc}_{\nu,\nu})^\vi\hat{\tau}_{1,1}&=&(1,\frac{1}{2})(-1,-\frac{1}{4})\\
&=&(-1,-\frac{1}{2}-\frac{1}{4})\\
&=&(-1,-\frac{3}{4}).
\end{eqnarray*}
This shows that $L_{1,1}(\nu,\nu)=0$.

Suppose now that $n>1$ or $m>1$. By using formula (\ref{inklu}) one can easily check that for any $k\geq 1$
\begin{eqnarray*}
\S^k(\aleph^{\sa}_{\nu,\nu})^\vi&=& e_{k-1}e_k,\\
\S^k(\aleph^{\sc}_{\nu,\nu})^\vi&=&e_ke_{k-1},
\end{eqnarray*}
in $\widetilde{O}(k+2)$.

For any $k\geq1$ the shuffle permutation $\tau_{1,k-1}$ can be decomposed as
\begin{eqnarray*}
\tau_{1,k-1}&=&(k\; k-1)\cdots(2\; 1),
\end{eqnarray*}
hence for fixed $p,q\geq 0$ with $p+k+q>2$ we can lift $S^p\wedge\tau_{1,k-1}\wedge S^q$ to $\widetilde{O}(p+k+q)$ by
\begin{eqnarray*}
\hat{\tau}_{1,k-1}&=&\frac{1}{2^{\frac{k-1}{2}}}(e_{p+k}-e_{p+k-1})\cdots(e_{p+2}-e_{p+1}),
\end{eqnarray*}
and its inverse in $\widetilde{O}(p+k+q)$ is
\begin{eqnarray*}
\hat{\tau}_{k-1,1}&=&
\frac{1}{2^{\frac{k-1}{2}}}(e_{p+2}-e_{p+1})\cdots(e_{p+k}-e_{p+k-1}).
\end{eqnarray*}
By using these equalities we obtain
\begin{eqnarray*}
(\aleph^\sa_{n,\nu,m,\nu})^\vi&=&(S^{n-1}\wedge\tau_{m-1,1}\wedge S^1)(\S^{m+n-2}(\aleph^{\sa}_{\nu,\nu})^\vi)
(S^{n-1}\wedge\tau_{1,m-1}\wedge S^1)\\
&=&\hat{\tau}_{m-1,1}(\S^{m+n-2}(\aleph^{\sa}_{\nu,\nu})^\vi)
\hat{\tau}_{1,m-1}\\
&=&\frac{1}{2^{m-1}}(e_{n+1}-e_n)\cdots(e_{n+m-1}-e_{n+m-2})e_{n+m-1}e_{n+m}\\
&&\cdot(e_{n+m-1}-e_{n+m-2})\cdots(e_{n+1}-e_n),\\
(\aleph^{\sc}_{m,\nu,n,\nu})^\vi&=&(S^{m-1}\wedge\tau_{n-1,1}\wedge S^1)(\S^{n+n-2}(\aleph^{\sa}_{\nu,\nu})^\vi)
(S^{m-1}\wedge\tau_{1,n-1}\wedge S^1)\\
&=&\hat{\tau}_{n-1,1}(\S^{n+m-2}(\aleph^{\sa}_{\nu,\nu})^\vi)
\hat{\tau}_{1,n-1}\\
&=&\frac{1}{2^{n-1}}(e_{m+1}-e_m)\cdots(e_{n+m-1}-e_{n+m-2})e_{n+m}e_{n+m-1}\\
&&\cdot(e_{n+m-1}-e_{n+m-2})\cdots(e_{m+1}-e_m).
\end{eqnarray*}
For $j>i$ we have the following useful identities.
\begin{equation*}\tag{a}
e_j(e_j-e_i)=-(e_j-e_i)e_i,
\end{equation*}
\begin{equation*}\tag{b}
e_i(e_j-e_i)=-(e_j-e_i)e_j,
\end{equation*}
\begin{equation*}\tag{c}
(e_j-e_i)^2=2.
\end{equation*}
By using (a), (b), (c), and the formulas above one can easily check that
\begin{eqnarray*}
(\aleph^\sa_{n,\nu,m,\nu})^\vi&=&e_ne_{n+m},\\
(\aleph^{\sc}_{m,\nu,n,\nu})^\vi&=&e_{n+m}e_m.
\end{eqnarray*}

The identity $\tau_{n,m}=\tau_{1,n+m-1}^n$ holds, therefore we can also take 
$\hat{\tau}_{n,m}=\hat{\tau}_{1,n+m-1}^n$. With this choice it is not difficult to compute that 
\begin{eqnarray*}
\hat{\tau}_{n,m}(\aleph^{\sa}_{n,\nu,m,\nu})^\vi&=&(\aleph^{\sc}_{m,\nu,n,\nu})^\vi\hat{\tau}_{n,m}
\end{eqnarray*}
in $\widetilde{O}(n+m)$. For this one uses (a) and (b). This last equation is equivalent to $L_{n,m}(\nu,\nu)=0$.
\end{proof}

\section{Properties of the smash product in dimensions $\geq 1$}\label{lasP} 

Restricting to dimensions $\geq 1$ in this section we show a series of properties of the smash product operation for
secondary homotopy groups in Definition \ref{ka} which imply Theorems \ref{lamel}, \ref{lmf}, \ref{nosym} and
\ref{yi} within this range. The case of dimension $0$ is a consequence of the fact that secondary homotopy groups are track functors and of the first technical lemma in the next section.

In this section we will work with the track category $\C{Top}^*$. This track category has a strict zero object $0$ so that the zero morphism $0\colon X\r Y$ is always defined for a pair of pointed spaces. In this situation the \emph{golden rule} says that
\begin{equation*}\tag{GR}
\begin{array}{l}
\text{The composition of a trivial map $0$ with a any track $F$ is always a }\\ \text{trivial track }0F=0^\vc, F0=0^\vc.
\end{array}
\end{equation*}
This is an obvious but crucial property that will be very useful for computations.

Let $n,m\geq 1$ and $0\leq i,j,i+j\leq 1$.

We need to show that the square group morphisms (\ref{ellos}) are
well defined. There is nothing to check in case $i=j=0$. For $i+j=1$
we define operations
\begin{equation}\label{eli}
\wa,\wc\colon\pi_{n,i}X\times\pi_{m,j}Y\To\Pi_{n+m,1}(X\wedge Y)
\end{equation}
as in Definition \ref{ka}. For this we need to choose $\bar{g}$ more
carefully in certain cases. For example, in order to define $g\wa
[f,F]$ and $g\wc [f,F]$ when  $n=1$ and $i=0$
 we need to
take $\bar{g}$ in such a way that $(\pi_1\bar{g})(1)=g\in\grupo{\L X}$.

It is not completely immediate that the operations $\wa$ and $\wc$
do not depend on choices. We check here for instance that $g\wa
[f,F]$ does not depend on the choice of $\bar{g}$, $f$ and $F$.
Let $\bar{g}'$, $f'$ and $F'$ be another choice. Then there are
tracks $N_1\colon \S^{n-1}\bar{g}'\rr\S^{n-1}\bar{g}$ and
$N_2\colon \S^{m-1}f'\rr\S^{m-1}f$ with trivial Hopf invariant such that $F'=F\vc(ev\, N_2)$.
One can now use Lemmas \ref{hopfs}, \ref{exteprop2}, and \ref{nilson} 
together with the golden rule to show that, up to a permutation in the spherical coordinates, the
track $\S(N_1\sa N_2)\colon \S^{n+m-1}(\bar{g}'\sa
f')\rr\S^{n+m-1}(\bar{g}\sa f)$
is a track with trivial Hopf invariant which determines the desired
equality.

In the next lemma we establish the fundamental properties
of (\ref{eli}).

\begin{lem}\label{P1}
\begin{enumerate}
\item The operation $\wa$ is left linear
\begin{eqnarray*}
(x+y)\wa[f,F]&=&x\wa[f,F]+y\wa[f,F],\\
([f,F]+[g,G])\wa x&=&[f,F]\wa x+[g,G]\wa x.
\end{eqnarray*}

\item The operation $\wc$ is right linear
\begin{eqnarray*}
[f,F]\wa(x+y)&=&[f,F]\wc x+[f,F]\wc y,\\
x\wc([f,F]+[g,G])&=&x\wc[f,F]+x\wc[g,G].
\end{eqnarray*}

\item For $m\geq 2$, $i=0$ and $j=1$, given $a\in\otimes^2\Z[\L^mY]$ the  equality 
$x\wa P(a)=P(\otimes^2\Z[\wedge])(1\otimes \tau_\otimes\otimes 1)(\Delta(x)\otimes a)$ holds.

\item For $m\geq 2$, $i=1$ and $j=0$, given $a\in\otimes^2\Z[\L^mY]$ the  equality 
$[f,F]\wa \partial P(a)=P(\otimes^2\Z[\wedge])(1\otimes \tau_\otimes\otimes
1)(\Delta\partial[f,F]\otimes a)$ holds.

\item For $n\geq 2$, $i=1$ and $j=0$, given $a\in\otimes^2\Z[\L^nX]$ the  equality 
$P(a)\wc{x}=P(\otimes^2\Z[\wedge])(1\otimes \tau_\otimes\otimes 1)(a\otimes \Delta(x))$ holds.

\item For $n\geq 2$, $i=0$ and $j=1$, given $a\in\otimes^2\Z[\L^nX]$ the following equality holds
$$\partial P(a)\wc [f,F]=P(\otimes^2\Z[\wedge])(1\otimes \tau_\otimes\otimes 1)(a\otimes
\Delta\partial[f,F]).$$

\item If $m=1$, $i=0$ and $j=1$, the following equality holds
$$g\wa (-[f,F]+[f,F]^x)=P(\otimes^2\Z[\wedge])(1\otimes \tau_\otimes\otimes 1)(\Delta(g)\otimes
\set{x}\otimes\set{\partial[f,F]}).$$ 

\item For $n=1$, $i=1$ and $j=0$, the following equality holds
$$(-[f,F]+[f,F]^x)\wc g=P(\otimes^2\Z[\wedge])(1\otimes \tau_\otimes\otimes
1)(\set{x}\otimes\set{\partial[f,F]}\otimes \Delta(g)).$$ 

\item The equality  
$g\wc[f,F]-g\wa[f,F]=P(\otimes^2\Z[\wedge])(1\otimes\tau_\otimes\otimes1)(H(g)\otimes
TH\partial[f,F])$ holds.

\item The equality 
$[f,F]\wc g-[f,F]\wa
g=P(\otimes^2\Z[\wedge])(1\otimes\tau_\otimes\otimes1)(H\partial[f,F]\otimes
TH(g))$ holds.

\item The  equalities 
$(\partial[f,F])\wa[g,G]=[f,F]\wa(\partial[g,G])$ and
$(\partial[f,F])\wc[g,G]=[f,F]\wc(\partial[g,G])$ hold.
\end{enumerate}
\end{lem}

\begin{proof}
Let us check the first equation in (1).  
\begin{eqnarray*}
x\wa[f,F]+y\wa[f,F]&=&[(\S\wedge)(\bar{x}\sa f,\bar{y}\sa f)\mu,\\
&&((\bar{x}_{ev}\wedge F)\vc(ev(\S^{n+m}\wedge)\aleph^\sa_{n,\bar{x},m,f}),\\&&
(\bar{y}_{ev}\wedge F)\vc(ev(\S^{n+m}\wedge)\aleph^\sa_{n,\bar{y},m,f}))\mu]\\
\text{\ref{proep} (3) and (GR)}&=&[(\S\wedge)((\bar{x},\bar{y})\sa f)\mu,\\
&&((\bar{x}_{ev}\wedge F)\vc(ev(\S^{n+m}\wedge)\aleph^\sa_{n,\bar{x},m,f}),\\&&
(\bar{y}_{ev}\wedge F)\vc(ev(\S^{n+m}\wedge)\aleph^\sa_{n,\bar{y},m,f}))\aleph^\sa_{n,\mu,m,1_{S^1}}]\\
\text{\ref{proep} (7)}&=&[(\S\wedge)(((\bar{x},\bar{y})\mu)\sa f),\\
&&((\bar{x}_{ev}\wedge F,\bar{y}_{ev}\wedge F)(\S^{n-1}\mu\wedge S^m))\vc\\
&&(ev(\S^{n+m}\wedge)(\aleph^\sa_{n,\bar{x},m,f},\aleph^\sa_{n,\bar{y},m,f})\aleph^\sa_{n,\mu,m,1_{S^1}})]\\
\text{\ref{exteprop} (3) and (7)}&=&[(\S\wedge)(((\bar{x},\bar{y})\mu)\sa f),\\
&&(((\bar{x},\bar{y})\mu)_{ev}\wedge F)\vc
(ev(\S^{n+m}\wedge)\aleph^\sa_{n,(\bar{x},\bar{y})\mu,m,f})]\\
&=&(x+y)\wa[f,F].
\end{eqnarray*}
The second equation in (1) and the equations in (2) are analogous. 

In order to check equation (3) we notice that $x\wa P(a)$ is represented by the following diagram
$$\xymatrix{S^{n+m}\ar[rr]_-{\S^{n+m-1}(\bar{x}\sa f)}^{\;}="a"
\ar@/^20pt/[rr]_{\;}="b"|{\S^{n-1}\bar{x}\wedge\S^{m-1}f}^{\;}="c"
\ar@/^70pt/[rr]_{\;}="d"^0&&
\vee_{\L^nX\wedge\L^mY}S^{n+m}\ar[rr]_-{\S^{n+m}\wedge}&&S^{n+m}_{X\wedge Y}\ar[r]_-{ev}&X\wedge Y
\ar@{=>}"a";"b"_{\aleph^\sa_{n,\bar{x},m,f}}
\ar@{=>}"c";"d"|<(.3){\S^{n-1}\bar{x}\wedge (\S^{m-2}F)}}$$
where $F\colon \S f\rr 0$ is any track with $\hopf(F)=-\tau_\otimes(a)$. Here we use claim (*) in the proof of \cite{2hg1} 4.9. By Lemmas \ref{exteprop2} (1), 
\ref{exteprop} (2) and the golden rule this diagram coincides with
$$\xymatrix{S^{n+m}\ar[rr]_-{\S^{n+m-1}(\bar{x}\sa f)}^{\;}="a"
\ar@/^50pt/[rr]_{\;}="d"^0&&
\vee_{\L^nX\wedge\L^mY}S^{n+m}\ar[rr]_-{\S^{n+m}\wedge}&&S^{n+m}_{X\wedge Y}\ar[r]_-{ev}&X\wedge Y
\ar@{=>}"a";"d"|<(.3){\S((\S^{n-1}\bar{x})\sa (\S^{m-2}F))}}$$
Now one can use Lemmas \ref{hopfs} (2) and \ref{nilson} to check that 
$$\hopf(\S((\S^{n-1}\bar{x})\sa (\S^{m-2}F)))=-\bar{\sigma}\tau_\otimes(1\otimes\tau_\otimes\otimes1)(\Delta(x)\otimes a),$$
hence (3) follows form claim (*) in the proof of \cite{2hg1} 4.9.
Equation (5) is analogous. Equations (7) and (8) are unstable versions of (3) and (5). We leave them to the reader.
 
Both sides of the first equation in (11) are represented by
$$\xymatrix{S^{n+m}\ar[rr]_-{\S^{n+m-1}(g\sa f)}^{\;}="a"
\ar@/^20pt/[rr]_{\;}="b"^<(.6){\S^{n-1}g\wedge\S^{m-1}f}
\ar@/^50pt/[rrrrr]_<(.6){\;}="d"^0&&
\vee_{\L^nX\wedge\L^mY}S^{n+m}\ar[rr]_-{\S^{n+m}\wedge}^{\;}="c"&&S^{n+m}_{X\wedge Y}\ar[r]_-{ev}&X\wedge Y
\ar@{=>}"a";"b"_{\aleph^\sa_{n,g,m,f}}
\ar@{=>}"c";"d"_{G\wedge F}}$$
In order to check this fact one only needs to use the golden rule. Similarly for the second equation in (11). 
Now (4) and (6)
follow from (3), (5) and (11). 

Finally (9) and (10) follow from Theorem \ref{hopfcal}, Lemma \ref{nilson}, and claim (*) in the proof of \cite{2hg1} 4.9.
\end{proof}

The equalities in Lemma \ref{P1} can be used to check the following properties of the operations (\ref{eli}).

\begin{lem}\label{P2}
\begin{enumerate}
\item For $n=1$ and $i=0$ the elements $g\wa[f,F]$ and $g\wc[f,F]$ only depend on
$g\in\Pi_{1,0}X$. 

\item For $m=1$ and $j=0$ the elements $[f,F]\wa g$ and $[f,F]\wc g$ only depend on
$g\in\Pi_{1,0}Y$. 

\item For $m=2$ and $j=1$ the elements $g\wa[f,F]$ and $g\wc[f,F]$ only depend on
$[f,F]\in\Pi_{2,1}X$. 

\item For $n=2$ and $i=1$ the elements $[f,F]\wa g$ and $[f,F]\wc g$ only depend on
$[f,F]\in\Pi_{2,1}X$. 
\end{enumerate}
\end{lem}

Lemma \ref{P2} allows us to define the following operations when $m=j=1$ or $n=i=1$.
\begin{equation}\label{eli2}
\begin{array}{c}
\wa,\wc\colon\Pi_{n,0}X\times\left(\pi_{1,1}Y\times(\hat{\otimes}^2\Z[\L Y])\right)\To\Pi_{n+1,1}(X\wedge Y),\\{}\\
\wa,\wc\colon\left(\pi_{1,1}X\times(\hat{\otimes}^2\Z[\L
X])\right)\times\Pi_{m,0}Y\To\Pi_{1+m,1}(X\wedge Y),
\end{array}
\end{equation}
by the formulas
\begin{eqnarray*}
g\wa([f,F],a)&=&g\wa[f,F]+P(\otimes^2\Z[\wedge])(1\otimes\tau_\otimes\otimes1)(\Delta(g)\otimes a),\\
g\wc([f,F],a)&=&g\wc[f,F]+P(\otimes^2\Z[\wedge])(1\otimes\tau_\otimes\otimes1)(\set{g}\otimes\set{g}\otimes a),\\
([f,F],a)\wa g&=&[f,F]\wa g+P(\otimes^2\Z[\wedge])(1\otimes\tau_\otimes\otimes1)(a\otimes\set{g}\otimes\set{g}),\\
([f,F],a)\wc g&=&[f,F]\wc
g+P(\otimes^2\Z[\wedge])(1\otimes\tau_\otimes\otimes1)(a\otimes\Delta(g)).
\end{eqnarray*}

By using Lemma \ref{P1} one can check the following one.

\begin{lem}\label{P3}
\begin{enumerate}
\item For $m=j=1$ the elements $g\wa([f,F],a)$ and $g\wc([f,F],a)$ only depend on
$([f,F],a)\in\Pi_{1,1}Y$ 

\item For $n=i=1$ the elements $([f,F],a)\wa g$ and $([f,F],a)\wc g$ only depend on
$([f,F],a)\in\Pi_{1,1}X$. 
\end{enumerate}
\end{lem}

Moreover, one can extend all properties (1)--(6) and (9)--(11) in Lemma \ref{P1} to the operations in (\ref{eli2}). The explicit statement is left to the reader.

In order to prove the following Lemma one uses Lemma \ref{exteprop} (11) and (12). 
We leave the details as
an exercise.

\begin{lem}\label{P4}
The following associativity rules for the operations (\ref{eli}) hold.
\begin{enumerate}
\item  $g\wa(g'\wa[f,F])=(g\sa g')\wa[f,F]$,

\item $g\wa([f,F]\wa g')=(g\wa[f,F])\wa g'$,

\item $[f,F]\wa(g\sa g')=([f,F]\wa g)\wa g'$,

\item $g\wc(g'\wc[f,F])=(g\sc g')\wc[f,F]$,

\item $g\wc([f,F]\wc g')=(g\wc[f,F])\wc g'$,

\item $[f,F]\wc(g\sc g')=([f,F]\wc g)\wc g'$.
\end{enumerate}
\end{lem}

One can accordingly obtain associativity properties involving also the operations in (\ref{eli2}).

Recall that $\hat{\tau}_{n,m}\colon\tau_{n,m}\rr(\cdot)^{(-1)^{nm}}_{n+m}\in\symt{n+m}$ is a track from the
shuffle permutation in (\ref{chuf}) to its sign, and $\tau_\wedge$ is the symmetry isomorphism for the smash product of pointed spaces.
In the following lemma we establish the commutativity rule for the operations (\ref{eli}). One can
similarly state the commutativity rules for (\ref{eli2}).

\begin{lem}\label{P5}
Given two pointed spaces $X$, $Y$, $[f,F]\in\Pi_{n,1}X$, and $g\in\Pi_{m,0}Y$ we have the equalities
\begin{enumerate}
\item $\grupo{g\wc\partial[f,F],\hat{\tau}_{n,m}}=-(\tau_\wedge)_*([f,F]\wa g)+((-1)^{nm})^* (g\wc[f,F])$,
\item $\grupo{g\wa\partial[f,F],\hat{\tau}_{n,m}}=-(\tau_\wedge)_*([f,F]\wc g)+((-1)^{nm})^* (g\wa[f,F])$.
\end{enumerate}
Moreover if $g\in\Pi_{n,0}X$ and $[f,F]\in\Pi_{m,1}Y$
\begin{enumerate}\setcounter{enumi}{2}
\item $\grupo{\partial[f,F]\wc g,\hat{\tau}_{n,m}}=-(\tau_\wedge)_*(g\wa [f,F])+((-1)^{nm})^* ([f,F]\wc g)$,
\item $\grupo{\partial[f,F]\wa g,\hat{\tau}_{n,m}}=-(\tau_\wedge)_*(g\wc [f,F])+((-1)^{nm})^* ([f,F]\wa g)$.
\end{enumerate}
\end{lem}

\begin{proof}
Here we prove (1). The other equations are analogous. 

The following diagram of pointed sets is commutative.
$$\xymatrix{\L^nX\wedge\L^mY\ar[rr]^\wedge\ar[d]_{\tau_\wedge}&&\L^{n+m}(X\wedge Y)\ar[d]^{\L^{n+m}\tau_\wedge}\\
\L^mY\wedge\L^n X\ar[r]^\wedge&\L^{m+n}(Y\wedge X)\ar[r]^{\tau_{n,m}^*}&\L^{n+m}(Y\wedge X)}$$
Therefore using the definition of induced morphisms for secondary homotopy groups in \cite{2hg1} 4.2 we have
\begin{eqnarray*}
(\tau_\wedge)_*([f,F]\wa g)&=&[(\S\L^{n+m}\tau_\wedge)(\S\wedge)(f\sa \bar{g}),\\
&&\tau_\wedge((F\wedge \bar{g}_{ev})\vc(ev(\S^{n+m}\wedge)\aleph^\sa_{n,f,m,g}))]\\
&=&[(\S\tau_{n,m}^*)(\S\wedge)(\S\tau_\wedge)(f\sa \bar{g}),\\
&&(\tau_\wedge(F\wedge \bar{g}_{ev}))\vc(\tau_\wedge ev(\S^{n+m}\wedge)\aleph^\sa_{n,f,m,g}))]\\
(\ref{tee})&=&[(\S\tau_{n,m}^*)(\S\wedge)(\bar{g}\sc f),\\
&&(\tau_\wedge(F\wedge \bar{g}_{ev}))\vc(ev(\S^{n+m}\wedge)(\tau_{n,m}\wedge\tau_\wedge)\aleph^\sa_{n,f,m,g}))]\\
&=&\text{(a)}.\\
\end{eqnarray*}

On the other hand using claim (*) in the proof of \cite{2hg1} 4.9, if 
$$Q\colon(\S^{n+m}\wedge)((\cdot)^{(-1)^{nm}}_{n+m}\wedge\tau_\wedge)(\S^{n+m-1}(f\sa \bar{g}))\rr(\S^{n+m}\wedge)(\S^{n+m-1}(\bar{g}\sc f))(\cdot)^{(-1)^{nm}}_{n+m}$$
is a track with $$\overline{\hopf}(Q)=-\bar{\sigma}\tau_\otimes\binom{(-1)^{nm}}{2}H(g\wc\partial[f,F])$$
then
\begin{eqnarray*}
((-1)^{nm})^* (g\wc[f,F])&=&[((\cdot)^{(-1)^{nm}}\wedge\L^{m+n}(Y\wedge X))(\S\wedge)(\bar{g}\sc f),\\
&&((\bar{g}_{ev}\wedge
F)(\cdot)^{(-1)^{nm}}_{m+n})\\
&&\vc(ev(\S^{m+n}\wedge)\aleph^{\sc}_{m,\bar{g},n,f}(\cdot)^{(-1)^{nm}}_{m+n})\vc (ev\; Q)]\\
\text{(GR)}&=&[((\cdot)^{(-1)^{nm}}\wedge\L^{m+n}(Y\wedge X))(\S\wedge)(\bar{g}\sc f),\\
&&((\bar{g}_{ev}\wedge F)\tau_{n,m})\\
&&\vc(ev(\S^{m+n}\wedge)\aleph^{\sc}_{m,\bar{g},n,f}\hat{\tau}_{n,m}^\vi)\vc (ev\; Q)]\\
\text{(Theorem \ref{hopfcal2})}&=&[((\cdot)^{(-1)^{nm}}\wedge\L^{m+n}(Y\wedge X))(\S\wedge)(\bar{g}\sc f),\\
&&(\tau_\wedge(F\wedge \bar{g}_{ev}))
\vc(ev(\S^{n+m}\wedge)(\hat{\tau}_{n,m}^\vi\wedge\tau_\wedge)\aleph^\sa_{n,f,m,\bar{g}})]\\
&=&\text{(b)}.\\
\end{eqnarray*}
Finally given a map $\varepsilon\colon S^1\r S^1\vee S^1$ with $(\pi_1\varepsilon)_\nill(1)=-a+b\in\grupo{a,b}_\nill$ and a track $N\colon (1,1)(\S^{n+m-1}\varepsilon)\rr 0$ with $\overline{\hopf}(N)=0$
\begin{eqnarray*}
-\text{(a)}+\text{(b)}&=&[((\S\tau_{n,m}^*)(\S\wedge)(\bar{g}\sc f),
((\cdot)^{(-1)^{nm}}\wedge\L^{m+n}(Y\wedge X)(\S\wedge)(\bar{g}\sc f))\varepsilon,\\
&&(\tau_\wedge(F\wedge \bar{g}_{ev},F\wedge \bar{g}_{ev})(\S^{n+m}\varepsilon))\\
&&\vc(ev(\S^{n+m}\wedge)((\tau_{n,m}\wedge\tau_\wedge)\aleph^\sa_{n,f,m,g},\\
&&(\hat{\tau}_{n,m}^\vi\wedge\tau_\wedge)\aleph^\sa_{n,f,m,g})(\S^{n+m}\varepsilon))]\\
&=&[((\S\tau_{n,m}^*)(\S\wedge)(\bar{g}\sc f),
((\cdot)^{(-1)^{nm}}\wedge\L^{m+n}(Y\wedge X)(\S\wedge)(\bar{g}\sc f))\varepsilon,\\
&&(\tau_\wedge(F\wedge \bar{g}_{ev})(1,1)(\S^{n+m}\varepsilon))\\
&&\vc(ev(\S^{n+m}\wedge)(\hat{\tau}_{n,m}^\vi\wedge\tau_\wedge)\aleph^\sa_{n,f,m,g})(1,1)(\S^{n+m}\varepsilon)) \\
&&\vc(ev(\S^{n+m}\wedge)((\hat{\tau}_{n,m}\wedge\tau_\wedge)(\S^{n+m-1}(f\sa \bar{g})),\\
&&((\cdot)_{n+m}^{(-1)^{nm}}\wedge\tau_\wedge)(\S^{n+m-1}(f\sa \bar{g})))(\S^{n+m}\varepsilon))]\\
\text{(\ref{tee})}&=&[((\S\tau_{n,m}^*)(\S\wedge)(\bar{g}\sc f),
((\cdot)^{(-1)^{nm}}\wedge\L^{m+n}(Y\wedge X)(\S\wedge)(\bar{g}\sc f))\varepsilon,\\
&&(((\tau_\wedge(F\wedge \bar{g}_{ev}))\\
&&\vc(ev(\S^{n+m}\wedge)(\hat{\tau}_{n,m}^\vi\wedge\tau_\wedge)\aleph^\sa_{n,f,m,g}))(1,1)(\S^{n+m}\varepsilon)) \\
&&\vc(ev(\S^{n+m}\wedge)((\hat{\tau}_{n,m}\wedge\L^mY\wedge\L^nX)(\S^{n+m-1}(\bar{g}\sc f)),\\
&&((\cdot)_{n+m}^{(-1)^{nm}}\wedge\L^mY\wedge\L^nX)(\S^{n+m-1}(\bar{g}\sc f)))(\S^{n+m}\varepsilon))]\\
\text{(GR)}&=&[((\S\tau_{n,m}^*)(\S\wedge)(\bar{g}\sc f),
((\cdot)^{(-1)^{nm}}\wedge\L^{m+n}(Y\wedge X))(\S\wedge)(\bar{g}\sc f))\varepsilon,\\
&&(((\tau_\wedge(F\wedge \bar{g}_{ev}))\\
&&\vc(ev(\S^{n+m}\wedge)(\hat{\tau}_{n,m}^\vi\wedge\tau_\wedge)\aleph^\sa_{n,f,m,g}))N) \\
&&\vc(ev(\S^{n+m}\wedge)((\hat{\tau}_{n,m}\wedge\L^mY\wedge\L^nX)(\S^{n+m-1}(\bar{g}\sc f)),\\
&&((\cdot)_{n+m}^{(-1)^{nm}}\wedge\L^mY\wedge\L^nX)(\S^{n+m-1}(\bar{g}\sc f)))(\S^{n+m}\varepsilon))]\\
\text{(GR), (\ref{tee})}&=&[((\S\tau_{n,m}^*)(\S\wedge)(\bar{g}\sc f),
((\cdot)^{(-1)^{nm}}\wedge\L^{m+n}(Y\wedge X))(\S\wedge)(\bar{g}\sc f))\varepsilon,\\
&&(ev(\S^{n+m}\wedge)(S^{n+m-1}\wedge(\cdot)^{(-1)^{nm}}\wedge\L^mY\wedge\L^nX)\\&&(\S^{n+m-1}(\bar{g}\sc f))N) \\
&&\vc(ev(\S^{n+m}\wedge)((\hat{\tau}_{n,m}\wedge\L^mY\wedge\L^nX)(\S^{n+m-1}(\bar{g}\sc f)),\\
&&((\cdot)_{n+m}^{(-1)^{nm}}\wedge\L^mY\wedge\L^nX)(\S^{n+m-1}(\bar{g}\sc f)))(\S^{n+m}\varepsilon))]\\
&=&\grupo{g\wc\partial[f,F],\hat{\tau}_{n,m}}
\end{eqnarray*}
Here we use the definition of the bracket operation $\grupo{-,-}$ in \cite{2hg2} 4.5.
\end{proof}

Finally we prove the compatibility of the smash product operation with the action of the symmetric track group.

\begin{lem}\label{P7}
Let $f\in\Pi_{n,0}X$ and $g\in\Pi_{m,0}Y$. Given $\hat{\sigma}\in\symt{m}$ and $\hat{\gamma}\in\symt{n}$ with $\delta(\hat{\sigma})=\sigma$ and $\delta(\hat{\gamma})=\gamma$ the
following equalities hold.
\begin{enumerate}
\item $\begin{array}{l}
\grupo{f\wa g,S^n\wedge\hat{\sigma}}=f\wa\grupo{g,\hat{\sigma}}\\+P(\otimes^2\Z[\wedge])(1\otimes\tau_\otimes\otimes1)(H(f)\otimes(-\sigma^*g+(\sign\sigma)^*g|\sigma^*g)_H),
\end{array}$

\item $\begin{array}{l}
\grupo{f\wc g,\hat{\tau}_{n,m}^{-1}(S^m\wedge\hat{\gamma})\hat{\tau}_{n,m}}=\grupo{f,\hat{\gamma}}\wc
g\\+P(\otimes^2\Z[\wedge])(1\otimes\tau_\otimes\otimes1)((-\sigma^*f+(\sign \sigma)^*f|\sigma^*f)_H\otimes H(g)).
\end{array}$
\end{enumerate}
\end{lem}

\begin{proof}
Equation (2) follows from (1), Lemma \ref{P5}, and the laws of a sign group action. Let us prove (1). By
Lemma \ref{exteprop} (1), (2), (5), (6), (8) and (9), given $\varepsilon\colon S^1\r S^1\vee S^1$ with $(\pi_1\varepsilon)_\nill(1)=-a+b\in\grupo{a,b}_\nill$ 
$$\begin{array}{l}
\aleph^\sa_{n,\bar{f},m,(\S\sigma^*,(\cdot)^{\sign\sigma}\wedge\L^mY)(\bar{g}\vee\bar{g})\varepsilon}\;\;\;=\\
(\S^{n+m}(\L^nX\wedge\sigma^*),(\cdot)_{n+m}^{\sign\sigma}\wedge\L^nX\wedge\L^mY)
\aleph^\sa_{n,\bar{f},m,(\bar{g}\vee\bar{g})\varepsilon},\\
\aleph^{\sc}_{n,\bar{f},m,(\bar{g}\vee\bar{g})\varepsilon}\;\;\;=\;\;\;
(\aleph^{\sc}_{n,\bar{f},m,\bar{g}}\vee\aleph^{\sc}_{n,\bar{f},m,\bar{g}})(\S^{n+m-1}\varepsilon).
\end{array}$$
Now one can use Theorem \ref{hopfcal}, Lemma \ref{nilson} and the elementary properties of the Hopf invariant for tracks in
\cite{2hg1} 3.6 to check that any track $Q$ from
\begin{small}
$$(\S^{n+m}(\L^nX\wedge\sigma^*),(\cdot)^{\sign\sigma}_{n+m}\wedge\L^nX\wedge\L^mY)
((\S^{n+m-1}\bar{f}\sa\bar{g})\vee(\S^{n+m-1}\bar{f}\sa\bar{g}))(\S^{n+m-1}\varepsilon)$$ 
\end{small}
to
$$\S^{n+m-1}(\bar{f}\sa((\S\sigma^*,(\cdot)^{\sign\sigma}\wedge\L^mY)(\bar{g}\vee\bar{g})\varepsilon))$$
with
$$\overline{\hopf(Q)}=-\bar{\sigma}\tau_\otimes(1\otimes\tau_\otimes\otimes1)(H(f)\otimes(-\sigma^*g+(\sign\sigma)^*g|\sigma^*g)_H),$$
satisfies
\begin{small}
$$\begin{array}{l}
\aleph^\sa_{n,\bar{f},m,(\S\sigma^*,(\cdot)^{\sign\sigma}\wedge\L^mY)(\bar{g}\vee\bar{g})\varepsilon}\vc Q\;\;\;=\\
(\S^{n+m}(\L^nX\wedge\sigma^*),(\cdot)_{n+m}^{\sign\sigma}\wedge\L^nX\wedge\L^mY)
(\aleph^\sa_{n,\bar{f},m,\bar{g}}\vee\aleph^\sa_{n,\bar{f},m,\bar{g}})(\S^{n+m-1}\varepsilon).
\end{array}$$
\end{small}
Then given $N\colon (1,1)(\S^{m-1}\varepsilon)\rr 0$ with $\overline{\hopf}(N)=0$, by using claim (*) in the proof of 
\cite{2hg1} 4.9 and the definition of the bracket $\grupo{-,-}$ in \cite{2hg2} 4.5, the right hand side of equation (1) is
$$\begin{array}{l}
[(\S\wedge)(\S(\L^nX\wedge\sigma^*),(\cdot)^{\sign\sigma}\wedge\L^nX\wedge\L^mY)
((\bar{f}\sa\bar{g})\vee(\bar{f}\sa\bar{g}))\varepsilon,\\
(\bar{f}_{ev}\wedge(ev((\cdot)_n^{\sign\sigma}\wedge\L^mY)(\S^{m-1}\bar{g})N))\\
\vc (\bar{f}_{ev}\wedge(ev(\hat{\sigma}\wedge\L^mY,0^\vc_{(\cdot)_n^{\sign\sigma}\wedge\L^mY})((\S^{m-1}\bar{g})\vee(\S^{m-1}\bar{g}))(\S^{m-1}\varepsilon)))\\
\vc(ev(\S^{n+m}\wedge)(\S^{n+m}(\L^nX\wedge\sigma^*),(\cdot)_{n+m}^{\sign\sigma}\wedge\L^nX\wedge\L^mY)\\
(\aleph^\sa_{n,\bar{f},m,\bar{g}}\vee\aleph^\sa_{n,\bar{f},m,\bar{g}})(\S^{n+m-1}\varepsilon))]=\\

[(\S\wedge)(\S(\L^nX\wedge\sigma^*),(\cdot)^{\sign\sigma}\wedge\L^nX\wedge\L^mY)
((\bar{f}\sa\bar{g})\vee(\bar{f}\sa\bar{g}))\varepsilon,\\
(ev(\S^{n+m}\wedge)((\cdot)_{n+m}^{\sign\sigma}\wedge\L^nX\wedge
\L^mY)(\S^{n-1}\bar{f}\wedge\S^{m-1}\bar{g})(\S^nN))\\
\vc
(ev(\S^{n+m}\wedge)(S^n\wedge\hat{\sigma}\wedge\L^nX\wedge\L^mY,0^\vc_{(\cdot)_{n+m}^{\sign\sigma}\wedge\L^nX\wedge\L^mY})\\
((\S^{n-1}\bar{f}\wedge\S^{m-1}\bar{g})\vee(\S^{n-1}\bar{f}\wedge\S^{m-1}\bar{g}))(\S^{n+m-1}\varepsilon))\\
\vc(ev(\S^{n+m}\wedge)(\S^{n+m}(\L^nX\wedge\sigma^*),(\cdot)_{n+m}^{\sign\sigma}\wedge\L^nX\wedge\L^mY)\\
(\aleph^\sa_{n,\bar{f},m,\bar{g}}\vee\aleph^\sa_{n,\bar{f},m,\bar{g}})(\S^{n+m-1}\varepsilon))]=\\

[(\S\wedge)(\S(\L^nX\wedge\sigma^*),(\cdot)^{\sign\sigma}\wedge\L^nX\wedge\L^mY)
((\bar{f}\sa\bar{g})\vee(\bar{f}\sa\bar{g}))\varepsilon,\\
(ev(\S^{n+m}\wedge)((\cdot)_{n+m}^{\sign\sigma}\wedge\L^nX\wedge
\L^mY)(\S^{n-1}\bar{f}\wedge\S^{m-1}\bar{g})(\S^nN))\\
\vc
(ev(\S^{n+m}\wedge)(S^n\wedge\hat{\sigma}\wedge\L^nX\wedge\L^mY,0^\vc_{(\cdot)_{n+m}^{\sign\sigma}\wedge\L^nX\wedge\L^mY})\\
(\aleph^\sa_{n,\bar{f},m,\bar{g}}\vee\aleph^\sa_{n,\bar{f},m,\bar{g}})(\S^{n+m-1}\varepsilon))]=\\

[(\S\wedge)(\S(\L^nX\wedge\sigma^*),(\cdot)^{\sign\sigma}\wedge\L^nX\wedge\L^mY)
((\bar{f}\sa\bar{g})\vee(\bar{f}\sa\bar{g}))\varepsilon,\\
(ev(\S^{n+m}\wedge)((\cdot)_{n+m}^{\sign\sigma}\wedge\L^nX\wedge
\L^mY)(\aleph^\sa_{n,\bar{f},m,\bar{g}})(\S^nN))\\
\vc
(ev(\S^{n+m}\wedge)(S^n\wedge\hat{\sigma}\wedge\L^nX\wedge\L^mY,0^\vc_{(\cdot)_{n+m}^{\sign\sigma}\wedge\L^nX\wedge\L^mY})\\
((\S^{n+m-1}(\bar{f}\sa\bar{g}))\vee(\S^{n+m-1}(\bar{f}\sa\bar{g})))(\S^{n+m-1}\varepsilon))]=\\

[(\S\wedge)(\S(\L^nX\wedge\sigma^*),(\cdot)^{\sign\sigma}\wedge\L^nX\wedge\L^mY)
((\bar{f}\sa\bar{g})\vee(\bar{f}\sa\bar{g}))\varepsilon,\\
(ev(\S^{n+m}\wedge)((\cdot)_{n+m}^{\sign\sigma}\wedge\L^nX\wedge
\L^mY)(\S^{n+m-1}(\bar{f}\sa\bar{g}))(\S^nN))\\
\vc
(ev(\S^{n+m}\wedge)(S^n\wedge\hat{\sigma}\wedge\L^nX\wedge\L^mY,0^\vc_{(\cdot)_{n+m}^{\sign\sigma}\wedge\L^nX\wedge\L^mY})\\
((\S^{n+m-1}(\bar{f}\sa\bar{g}))\vee(\S^{n+m-1}(\bar{f}\sa\bar{g})))(\S^{n+m-1}\varepsilon))]=\\

\grupo{f\wa
g,S^n\wedge\hat{\sigma}}
\end{array}$$
Here we essentially use the golden rule, concretely for the fourth equation.
\end{proof}

\section{The smash product in dimension $0$}

The first lemma in this section implies that the smash product operation in Definition \ref{ka} is well defined when dimension $0$ is involved and Theorems \ref{lamel}, \ref{lmf}, \ref{nosym}, and \ref{yi} are satisfied also in this case.

Recall that $\C{qpm}$ is a track category. Therefore the morphism set $\hom_\C{qpm}(C,D)$ in $\C{qpm}$ is 
indeed a groupoid. This groupoid is pointed by the zero morphism. If $G_\vc$ is a sign group acting on $C$ and $D$ we can consider the full pointed subgroupoid of $G_\vc$-equivariant morphisms
$$\hom_{G_\vc}(C,D)\subset\hom_\C{qpm}(C,D).$$
If $G_\vc$ is the trivial sign group then this inclusion is always an equality.

A pointed groupoid $\C{G}$ gives rise to a stable quadratic module $\ad_3\ad_2\ad_1\C{G}$, compare Remark \ref{dever}. The low-dimensional group of this stable quadratic module is the free group of nilpotency class $2$ on the pointed set of objects, therefore if we define $H$ as in (\ref{znil}) we obtain a quadratic pair module $\tilde{\C{G}}$ which corresponds to $\ad_3\ad_2\ad_1\C{G}$ by the forgetful functor in Remark \ref{nse}. Compare \cite{2hg2} 1.15.

\begin{lem}\label{Pinf}
Let $C, D$ be quadratic pair modules endowed with an action of the sign group $G_\vc$, let $\C{G}$ be a groupoid, and let
$$\varphi\colon\C{G}\To\hom_{G_\vc}(C,D)$$
be a pointed groupoid morphism. Suppose that $C$ is $0$-good. 
Then there is a well-defined quadratic pair module morphism
$$\tilde{\varphi}\colon\tilde{\C{G}}\odot C\To D$$
given by $$\tilde{\varphi}_i(g\ca x)=\varphi(g)(x),$$
where $g$ is an object in $\C{G}$ and $x\in C_i$ for some $i\in\set{0,1}$, or $g$ is a morphism in $\C{G}$, $x\in
C_0$ and $i=1$; and by
$$\tilde{\varphi}_{ee}((g|g')_H\otimes(x|x')_H)=(\varphi(g)(x)|\varphi(g')(x'))_H,$$
for $g,g'$ objects in $\C{G}$ and $x,x'\in C_0$. This morphism is $G_\vc$-equivariant.
\end{lem}

The proof of the lemma is technical but straightforward. The reader can also check that he construction is natural in $\C{G}$, $G_\vc$, $C$ and $D$.

Finally we show how Lemma \ref{cilin} can be used to obtain the tracks in $\C{qpm}$ induced by the additive
secondary homotopy groups from the smash product operation. For this we notice that there is a unique morphism in
$\C{qpm}$ from the interval quadratic pair module $\I$ in Section \ref{tcqpm} to $\Pi_{0,*}$ of the interval
$I_+$,
$$v\colon\I\To\Pi_{0,*}I_+,$$
sending $i_k\in\I_0$, $k=0,1$, to the map $i_k\colon S^0\r I_+\in\Pi_{0,0}I_+$ corresponding to the inclusion of
$k\in I=[0,1]$.

\begin{lem}\label{inri}
Given a track $F\colon f\rr g$ between maps $f,g\colon X\r Y$ represented by a homotopy $F\colon I_+\wedge X\r Y$ the composite
$$\I\odot\Pi_{n,*}X\st{v\odot1}\To\Pi_{0,*}I_+\odot\Pi_{n,*}X\st{\wedge}\To\Pi_{n,*}(I_+\wedge X)\st{\Pi_{n,*}F}\To\Pi_{n,*}Y$$
corresponds by Lemma \ref{cilin} to the track $\Pi_{n,*}F\colon\Pi_{n,*}f\rr\Pi_{n,*}g$ in $\C{qpm}$, $n\geq 0$.
\end{lem}

This lemma follows easily from the definition of the smash product operation for secondary homotopy groups.

\appendix

\section{Monoidal structures for graded quadratic pair modules and symmetric sequences}

Let $M$ be any additive monoid. The tensor product of two $M$-graded quadratic pair modules $D, E$ is defined as usually
$$(D\odot E)_n=\bigvee_{p+q=n}D_p\odot E_q.$$
However since the tensor product of quadratic pair modules does not preserve coproducts this does not define a
monoidal structure on $\C{qpm}^M$. 
For this we need to restrict to a subcategory of quadratic pair modules where the tensor produc preserves
coproducts. Below we check that the full subcategory of $0$-good quadratic pair modules, already introduced in Definition
\ref{qpm}, is suitable.

The class of good square groups is closed under tensor products and coproducts, see \cite{qaI}
Definition 2 and Section 5.6. The class of $0$-good quadratic pair modules is therefore closed under tensor
products and coproducts as well.
We show in Corollary \ref{pcop} below that the tensor product of $0$-good quadratic pair modules 
is well-behaved with respect to coproducts. As a consequence we obtain the following result.

\begin{prop}
The category $\C{qpm}^M_0$ of $M$-graded $0$-good quadratic pair modules is a monoidal category.
\end{prop}

The tensor product of symmetric sequences $X, Y$ in $\C{qpm}_0^{\symtt}$ is similarly well defined by the formula
$$(X\odot Y)_n=\bigvee_{p+q=n}(X_p\odot Y_q)\odot_{A(\symt{p}\tilde{\times}\symt{q})}A(\symt{n}).$$
Here, as usually, if $M$ and $N$ are a right and a left module over a quadratic pair algebra $R$, respectively, 
then $M\odot_RN$ denotes the coequalizer of the two multiplications $M\odot R\odot N\rightrightarrows M\odot N$.
The left $A(\symt{p}\tilde{\times}\symt{q})$-module structure of $A(\symt{n})$ is given by the sign group
morphism in Proposition \ref{inclusym}. Moreover, we use the
characterization of sign group actions given by Lemma \ref{rein} and the fact that $A$
is strict monoidal, see Proposition \ref{essmon}.

Now we state the technical results of this appendix. For any abelian group $B$ we define following \cite{qaI}
the square group $B^\otimes$ as $B^\otimes_e=B$, $B^\otimes_{ee}=B\oplus B$, $P=(1,1)$ and $H=\binom{1}{1}$.
Clearly $(-)^\otimes$ defines a functor from abelian groups to square groups.

\begin{lem}\label{ll}
Given square groups $L, M, N$ there is a natural push-out diagram
$$\xymatrix{(\coker P_L\otimes\coker P_M\otimes\coker P_L\otimes\coker P_N)^\otimes\ar[d]_\cong\ar@{^{(}->}[r]\ar@{}[rdd]|<(.7){\text{push}}&
(L\odot M)\vee(L\odot N)\ar[dd]^{(L\odot i_M, L\odot i_N)}\\
(\coker P_L\otimes\coker P_L\otimes\coker P_M\otimes\coker P_N)^\otimes\ar[d]_{(-|-)_H\otimes1\otimes1}&\\
(L_{ee}\otimes\coker P_M\otimes\coker P_N)^\otimes\ar@{^{(}->}[r]&L\odot(M\vee N)}$$
\end{lem}

\begin{proof}
The horizontal arrows are the monomorphisms in the exact sequences of square groups in \cite{qaI} Proposition 5 and 5.6. The compatibility of both sequences proves that the square of the statement is indeed a push-out.
\end{proof}

\begin{cor}\label{mm}
If $X$ is a good square group then the functor $X\odot-$ preserves coproducts. 
\end{cor}

\begin{proof}
By Lemma \ref{ll} the functor $X\odot-$ preserves finite coproducts. Since the tensor product $\odot$ of square
groups preserves always filtered colimits, see \cite{qaI} Proposition 4, then $X\odot-$ preserves indeed arbitrary coproducts.

\end{proof}

\begin{lem}
If $\partial\colon C_{(1)}\r C_{(0)}$ is a quadratic pair module such that $C_{(0)}\odot-$ preserves coproducts then for any two quadratic pair module morphisms $f_i\colon M_i\r N_i$, $i=1,2$, the natural morphism 
$$\Phi((\partial\odot f_1)\vee(\partial\odot f_2))\To\Phi(\partial\odot(f_1\vee f_2))$$
in $\C{qpm}$ is an epimorphism. Moreover it is an isomorphism provided $f_i$ are identity morphisms $i=1,2$.
\end{lem}

\begin{proof}
Since $C_{(0)}\odot-$ preserves coproducts the morphism is an isomorphism on the $(0)$-level.
The reflection functor $\Phi$ from pairs of square groups to quadratic pair modules is a left adjoint, so it
preserves colimts. 
In particular by Lemma \ref{ll} the first part of this lemma holds provided the natural morphism
\begin{equation*}\tag{a}
\begin{array}{c}
A=\Phi((\coker P_C\otimes\coker P_C\otimes\coker P_{M_1}\otimes\coker P_{M_2})^\otimes\r C_{(0)}\odot(N_1\vee N_2))\\\downarrow\\
B=\Phi((C_{ee}\otimes\coker P_{M_1}\otimes\coker P_{M_2})^\otimes\r C_{(0)}\odot(N_1\vee N_2))
\end{array}
\end{equation*}
induced by $(-|-)_H\colon\otimes^2\coker P_C\r C_{ee}$ is surjective on the $1$-level. 
On this level $A$ and $B$ are quotient groups, $A_1=\tilde{A}_1/\sim$, $B_1=\tilde{B}_1/\sim$, of the abelian
groups
\begin{equation*}
\begin{array}{rcl}
\tilde{A}_1&=&\coker P_C\otimes\coker P_C\otimes\coker P_{M_1}\otimes\coker P_{M_2}
\oplus D,\\
\tilde{B}_1&=&C_{ee}\otimes\coker P_{M_1}\otimes\coker P_{M_2}\oplus D.
\end{array}
\end{equation*}
and the morphism $A\r B$ is induced by $(-|-)_H\colon\otimes^2\coker P_C\r C_{ee}$ and the identity on $D$.
Here $D$ is the abelian group
$$D=C_{ee}\otimes((N_1)_{ee}\oplus (N_2)_{ee}
\oplus\coker P_{N_1}\otimes\coker P_{N_2}\oplus \coker P_{N_2}\otimes\coker P_{N_1}).$$
The relations $\sim$ in Section \ref{ape} and the definition of $(-)^\otimes$ show that $A_1$ and $B_1$ are generated by image of $D$, and hence (a) is surjective
on the $1$-level.

Now let $f_i=1\colon M_i=N_i$, $i=1,2$. In order to check the second part of the statement it is enough to show that the morphism
\begin{equation*}\tag{b}
\begin{array}{c}
A=\Phi((\coker P_C\otimes\coker P_C\otimes\coker P_{N_1}\otimes\coker P_{N_2})^\otimes\r C_{(0)}\odot(N_1\vee N_2))\\\downarrow\\
\Phi((C_{(1)}\odot N_1)\vee (C_{(1)}\odot N_2)\r C_{(0)}\odot(N_1\vee N_2))
\end{array}
\end{equation*}
given by the upper vertical arrow in Lemma \ref{ll} factors through (a). Both (a) and (b) are the identity on the $(0)$-level, so it is enough to look at the $(1)$-level. The desired factorization is induced by the unique square group morphism
\begin{equation*}
\begin{array}{c}
(C_{ee}\otimes\coker P_{N_1}\otimes\coker P_{N_2})^\otimes\\\downarrow\\
\Phi((C_{(0)}\odot N_1)\vee (C_{(0)}\odot N_2)\r C_{(0)}\odot(N_1\vee N_2))_1
\end{array}
\end{equation*}
which coincides on the $ee$-level with the square group morphism
$$(C_{ee}\otimes\coker P_{N_1}\otimes\coker P_{N_2})^\otimes\To
C_{(0)}\odot(N_1\vee N_2)$$
given by the lower horizontal arrow in Lemma \ref{ll}.
\end{proof}

\begin{lem}
Let $C$ be a quadratic pair module such that $C_{(0)}\odot-$ preserves coproducts. Then $C\odot-$ also preserves coproducts.
\end{lem}

\begin{proof}
Let us fisrt check that $C\odot-$ preserves the coproduct of two qadratic pair modules $D$, $E$. The functor $\Phi$ preserves colimits since it is a left adjoint.
The tensor product of quadratic pair modules is defined by $\Phi$ and the push-out construction  (\ref{tp}) in
the category $\C{SG}$ square groups which is also a push-out in the category $\C{Pair}(\C{SG})$ of pairs of
square groups. Therefore it is enough to check that the natural morphism
\begin{equation*}\tag{a}
\begin{array}{c}
\Phi((C_{(i)}\odot D_{(j)})\vee (C_{(i)}\odot E_{(j)})\r C_{(0)}\odot(D_{(0)}\vee E_{(0)}))\\\downarrow\\
\Phi(C_{(i)}\odot (D_{(j)}\vee E_{(j)})\r C_{(0)}\odot(D_{(0)}\vee E_{(0)}))
\end{array}
\end{equation*}
is an isomorphism if $i=0$ or $j=0$, and an epimorphism if $i=j=1$, and this follos from the previous lemma.

Since the coproduct of two objects is preserved by $C\odot-$ then all finite coproducts are preserved. Moreover,
by \cite{qaI} Proposition 4 the tensor product of square groups preserves filtered colimits. Since $\Phi$
preserves colimits then the tensor product of quadratic pair modules also preserves filtered colimits, hence
$C\odot-$ indeed preserves arbitrary coproducts.
\end{proof}

\begin{cor}\label{pcop}
If $C$ is a $0$-good quadratic pair module then $C\odot-$ preserves coproducts. 
\end{cor}

\bibliographystyle{amsalpha}
\bibliography{Fernando} 
\end{document}